%&amstex         
\input amstex\documentstyle{amsppt}  
\pagewidth{12.5 cm}\pageheight{19 cm}\magnification\magstep1
\topmatter
\title Character sheaves on disconnected groups, VII\endtitle
\author G. Lusztig\endauthor
\address Department of Mathematics, M.I.T., Cambridge, MA 02139\endaddress
\thanks Supported in part by the National Science Foundation.\endthanks
\endtopmatter   
\document
\define\btr{\blacktriangleright}
\define\cce{\che\ce}
\define\asi{\asymp}
\define\sqq{\sqrt q}
\define\pr{\text{\rm pr}}
\define\sgn{\text{\rm sgn}}
\define\part{\partial}
\define\Wb{\WW^\bul}
\define\Bpq{\Bumpeq}

\define\doa{\dot a}
\define\dua{\dot{\un a}}
\define\dw{\dot w}

\define\ds{\dot s}

\define\ua{\un a}

\define\uD{\un D}

\define\uZ{\un Z}
\define\ufs{\un{\fs}}

\define\mpb{\medpagebreak}

\define\hG{\hat G}

\define\hY{\hat Y}
\define\hb{\hat b}
\define\hf{\hat f}
\define\hy{\hat y}
\define\hcl{\hat{\cl}}
\define\hB{\hat B}

\define\hT{\hat T}

\define\hZ{\hat Z}

\define\htt{\hat\t}

\define\dsv{\dashv}

\define\po{\text{\rm pos}}

\define\si{\sim}

\define\sqc{\sqcup}

\define\qua{\quad}

\define\bao{\bar\o}
\define\bcl{\bar{\cl}}

\define\bp{\bar p}

\define\bu{\bar u}

\define\bK{\bar K}

\define\bS{\bar S}

\define\bY{\bar Y}
\define\bZ{\bar Z}

\define\bpi{\bar\p}

\define\bbq{\bar{\QQ}_l}

\define\lb{\linebreak}
\define\eSb{\endSb}

\define\bin{\binom}
\define\op{\oplus}
   
\redefine\sp{\spadesuit}
\define\em{\emptyset}
\define\imp{\implies}
\define\ra{\rangle}
\define\n{\notin}
\define\iy{\infty}
\define\m{\mapsto}
\define\do{\dots}
\define\la{\langle}

\define\lra{\leftrightarrow}

\define\sub{\subset}
\define\bxt{\boxtimes}
\define\T{\times}
\define\ti{\tilde}
\define\nl{\newline}
\redefine\i{^{-1}}
\define\fra{\frac}
\define\un{\underline}
\define\ov{\overline}
\define\ot{\otimes}

\define\Ad{\text{\rm Ad}}
\define\Hom{\text{\rm Hom}}
\define\End{\text{\rm End}}
\define\Aut{\text{\rm Aut}}

\define\Ker{\text{\rm Ker}}

\define\tr{\text{\rm tr}}

\define\supp{\text{\rm supp}}

\define\di{\diamond}

\define\a{\alpha}
\redefine\b{\beta}
\redefine\c{\chi}
\define\g{\gamma}
\redefine\d{\delta}
\define\e{\epsilon}
\define\et{\eta}
\define\io{\iota}
\redefine\o{\omega}
\define\p{\pi}
\define\ph{\phi}
\define\ps{\psi}
\define\r{\rho}
\define\s{\sigma}
\redefine\t{\tau}
\define\th{\theta}
\define\k{\kappa}
\redefine\l{\lambda}
\define\z{\zeta}
\define\x{\xi}

\define\vt{\vartheta}

\redefine\G{\Gamma}
\redefine\D{\Delta}
\redefine\O{\Omega}

\define\Th{\Theta}

\define\Ph{\Phi}
\define\Ps{\Psi}

\redefine\aa{\bold a}

\define\boc{\bold c}
\define\dd{\bold d}

\define\kk{\bold k}

\define\rr{\bold r}
\redefine\ss{\bold s}
\redefine\tt{\bold t}

\define\ww{\bold w}

\define\yy{\bold y}

\redefine\AA{\bold A}

\define\FF{\bold F}

\define\II{\bold I}

\define\NN{\bold N}

\define\QQ{\bold Q}
\define\RR{\bold R}
\define\SS{\bold S}
\define\TT{\bold T}

\define\WW{\bold W}
\define\ZZ{\bold Z}

\define\ca{\Cal A}
\define\cb{\Cal B}
\define\cc{\Cal C}
\define\cd{\Cal D}
\define\ce{\Cal E}
\define\cf{\Cal F}

\define\ch{\Cal H}
\define\ci{\Cal I}
\define\cj{\Cal J}

\define\cl{\Cal L}

\define\co{\Cal O}
\define\cp{\Cal P}

\define\car{\Cal R}
\define\cs{\Cal S}
\define\ct{\Cal T}
\define\cu{\Cal U}

\define\cz{\Cal Z}

\define\fa{\frak a}

\define\fs{\frak s}

\define\fA{\frak A}

\define\fC{\frak C}
\define\fD{\frak D}
\define\fE{\frak E}

\define\fK{\frak K}

\define\fP{\frak P}

\define\fS{\frak S}
\define\fT{\frak T}
\define\fU{\frak U}

\define\fZ{\frak Z}

\define\ta{\ti a}
\define\tb{\ti b}

\define\tf{\ti f}

\define\ts{\ti s}
\define\tit{\ti t}

\define\ty{\ti y}

\define\tE{\ti E}

\define\tG{\ti G}

\define\tR{\ti R}

\define\tT{\ti T}

\define\tV{\ti V}

\define\tY{\ti Y}
\define\tZ{\ti Z}

\define\tcl{\ti\cl}

\define\tid{\ti\d}
\define\tss{\ti\ss}

\define\tce{\ti\ce}

\define\sh{\sharp}

\define\bul{\bullet}
\define\che{\check}
\define\cha{\che\a}
\define\BBD{BBD}
\define\GI{Gi}
\define\KL{KL}
\define\CS{L3}
\define\AD{L9}
\define\YO{Y}
\define\GE{L13}
\define\HA{L12}
\define\MS{MS}

\head Introduction\endhead
Throughout this paper, $G$ denotes a fixed, not necessarily connected, reductive 
algebraic group over an algebraically closed field $\kk$. This paper is a part of a 
series \cite{\AD} which attempts to develop a theory of character sheaves on $G$. 

The usual convolution of class functions on a connected reductive group over a finite 
field makes sense also for complexes in $\cd(G^0)$ and then it preserves (see
\cite{\GI}) in the derived sense the class of character sheaves on $G^0$. In \S32 we 
define, more generally, a natural convolution operation for parabolic character sheaves
(see 32.21(a)). A key role in our study of convolution is played by Theorem 32.6 which
describes explicitly the convolution of two basic complexes of the form 
$\bK^{\ss,\cl}_{J,D}$ in terms of multiplication in some Hecke algebra. Using this we
define a map which to each parabolic character sheaf associates an orbit of a subgroup 
of the Weyl group on the set of isomorphism classes of "tame" local systems of rank $1$
on the torus $\TT$, see 32.25(b); in fact we define a refinement of this map in 
32.25(a). The main result of \S33 is Proposition 33.3 (a generalization of 
\cite{\CS, III, 14.2(b)}). It asserts that (under a cleanness assumption), the 
cohomology sheaves of a character sheaf restricted to an open subset of the support of 
a different character sheaf are disjoint from the local system given by the second 
character sheaf on that open subset. (This plays a key role in the argument in 35.22.) 
In \S34 we study the algebra $H_n$ of 31.2 (or rather an extension $H_n^D$ of it) in 
the spirit of our earlier study \cite{\HA} of a usual Iwahori-Hecke algebra by means of
the asymptotic Hecke algebra. This allows us to construct representations of $H_n^D$ 
starting from representations of $H_n^{D,1}$, the specialization of $H_n^D$ at $v=1$. 
In 34.19 we define some invariants $b_{A,u}^v$ of a character sheaf $A$ which depend 
also on an irreducible representation $E_u$ of $H_n^{D,1}$. These generalize the 
invariants $c_{A,E}$ of \cite{\CS, III,12.10}. From the definition, $b_{A,u}^v$ is a 
rational function in the indeterminate $v$ and one of our goals is to show that 
$b_{A,u}^v$ is in fact a constant. This goal is achieved in \S35 under a cleanness
assumption and a quasi-rationality assumption on $E_u$. (See Theorem 35.23 which is a 
generalization of \cite{\CS, III, 14.9}.) In \S35 we prove an orthogonality formula 
(35.15) for the characteristic functions of complexes of the form $\bK^{\ss,\cl}_{J,D}$
(over a finite field) in the spirit of \cite{\CS, III, 13.5}. A variant of this formula
(see 32.23) can be obtained in an entirely different way as an application of the 
results on convolution in \S32. As an application we associate a sign $\pm 1$ to any 
character sheaf on a connected component of $G$ (see 
35.17), under a cleanness assumption. This generalizes \cite{\CS, III, 13.10}.

\mpb

{\it Erratum to Part V.} In 25.6 replace $R_1^*=R^*\cap R$ by $R_1^*=R^*\cap R_1$.

{\it Erratum to Part VI.} In 28.19 replace $\cl'=(\uD\i)^*\cl$ by 
$\cl'=(\uD\i)^*\che\cl$. In 31.4 replace $\io_D$ by $\uD$.

\head Contents\endhead
 32. Convolution.

 33. Disjointness.

 34. The structure of $H^D_n$.

 35. Functions on $G^{0F}/U$.

\head 32. Convolution\endhead
\subhead 32.1\endsubhead
In this section we define and study the convolution of parabolic character sheaves.

\subhead 32.2\endsubhead
Let $\D$ be a connected component of $G$. Let $J\sub\II$. Let $\ss=(s_1,s_2,\do,s_r)$,
$\ss'=(s'_1,s'_2,\do,s'_{r'})$ be two sequences in $\II$. Let $a,a'\in\WW$. Let 
$$\ww=(s_1,s_2,\do,s_r,a,s'_1,s'_2,\do,s'_{r'},a'),\qua
[\ww]=s_1s_2\do s_ras'_1s'_2\do s'_{r'}a'.$$
Let $\cl\in\fs(\TT)$ be such that $[\ww]\un\D\in\Wb_\cl$. Let 
$$\ct=\{i\in[1,r];s_1\do s_{i-1}s_is_{i-1}\do s_1\in\WW_\cl\},$$ 
$$\ct'=\{j\in[1,r'];
a'{}\i s'_r\do s'_{j+1}s'_js'_{j+1}\do s'_{r'}a'\in\e_\D(\WW_\cl)\},$$ 
$$\align&\uZ^\ww_{\em,J,\D}=\{(B_0,B_1,\do,B_r,B'_0,B'_1,\do,B'_{r'},B,xU_{J,B_0});
B_i\in\cb(i\in[0,r]),\\&B'_j\in\cb(j\in[0,r']),B\in\cb,x\in\D,\po(B_{i-1},B_i)\in
\{s_i,1\}(i\in\ct),\\&\po(B_{i-1},B_i)=s_i(i\in[1,r]-\ct),\po(B'_{j-1},B'_j)\in
\{s'_j,1\}(j\in\ct'),\\&\po(B'_{j-1},B'_j)=s'_j(j\in[1,r']-\ct'),\po(B_r,B'_0)=a,
\po(B'_{r'},B)=a',\\&xB_0x\i=B\}.\endalign$$
Then $Z^\ww_{\em,J,\D}$ (see 28.8) is naturally an open subset of $\uZ^\ww_{\em,J,\D}$.
By 28.8, $\cl$ gives rise to a local system $\tcl$ on $Z^\ww_{\em,J,\D}$. 

(a) {\it $\tcl$ extends uniquely to a local system $\bcl$ on $\uZ^\ww_{\em,J,\D}$.}
\nl
Indeed, let $a=t_1t_2\do t_m,a'=t'_1t'_2\do t'_{m'}$ be reduced expressions for $a,a'$ 
in $\WW$ and let

$\tt=(s_1,s_2,\do,s_r,t_1,t_2,\do,t_m,s'_1,s'_2,\do,s'_{r'},t'_1,t'_2,\do,t'_{m'})$.
\nl
We identify in an obvious way $\uZ^\ww_{\em,J,\D}$ with an open subset of 
$\bZ^\tt_{\em,J,\D}$ (see 28.9) contained in 
$\cup_{\cj\sub\cj_\tt}Z^{\tt_\cj}_{\em,J,\D}$ (notation of 28.9) and we use the fact 
that $\tcl$, regarded as a local system on $Z^\ww_{\em,J,\D}=Z^\tt_{\em,J,\D}$ extends 
to a local system on $\cup_{\cj\sub\cj_\tt}Z^{\tt_\cj}_{\em,J,\D}$, see 28.10. This
extension is unique up to isomorphism since $Z^\ww_{\em,J,\D}$ is open dense in the 
smooth irreducible variety $\uZ^\ww_{\em,J,\D}$ (which is itself open dense in the 
smooth irreducible variety $\bZ^\tt_{\em,J,\D}$).

\subhead 32.3\endsubhead
For any $\aa=(a_0,a_1,\do,a_{r+r'})\in\WW^{r+r'+1}$ let 
$$\align\uZ^{\ww,\aa}_{\em,J,\D}&
=\{(B_0,B_1,\do,B_r,B'_0,B'_1,\do,B'_{r'},B,xU_{J,B_0})
\in\uZ^\ww_{\em,J,\D};\\&\po(B_k,B'_{r'})=a_k(k\in[0,r]),\po(B_r,B'_{r+r'-k})=a_k
(k\in[r,r+r'])\}.\endalign$$
Define $\p_{\ww,\aa}:\uZ^{\ww,\aa}_{\em,J,\D}@>>>Z_{J,\D}$ by
$$(B_0,B_1,\do,B_r,B'_0,B'_1,\do,B'{r'},B,xU_{J,B_0})\m(Q_{J,B_0},Q_{\e_\D(J),B},
xU_{J,B_0})$$
(notation of 28.7). Now $\uZ^{\ww,\aa}_{\em,J,\D}$ is empty unless

(i) $a_{r+r'}=a$,

(ii) $a_k\in\{a_{k-1},s_ka_{k-1}\}$ for $k\in[1,r]$,    

(iii) $a_k\in\{a_{k-1},a_{k-1}s'_{r+r'+1-k}\}$ for $k\in[r+1,r+r']$.
\nl
Indeed, let $(B_0,B_1,\do,B_r,B'_0,B'_1,\do,B'_{r'},B,xU_{J,B_0})\in
\uZ^{\ww,\aa}_{\em,J,\D}$. Clearly, (i) holds. Let $k\in[1,r]$. From 
$\po(B_k,B_{k-1})\in\{1,s_k\}$, $\po(B_{k-1},B'_{r'})=a_{k-1}$, we deduce 
$\po(B_k,B'_{r'})\in\{a_{k-1},s_ka_{k-1}\}$ and (ii) holds. Let $k\in[r+1,r+r']$. From 

$\po(B_r,B'_{r+r'+1-k})=a_{k-1}$, 
$\po(B'_{r+r'+1-k},B'_{r+r'-k})\in\{1,s'_{r+r'+1-k}\}$, 
\nl
we deduce $\po(B_r,B'_{r+r'-k})\in\{a_{k-1},a_{k-1}s'_{r+r'+1-k}\}$ and (iii) holds.

Let $\bcl_\aa$ be the restriction of the local system $\bcl$ from $\uZ^\ww_{\em,J,\D}$ 
to $\uZ^{\ww,\aa}_{\em,J,\D}$. We have a partition 
$\uZ^\ww_{\em,J,\D}=\sqc_\aa\uZ^{\ww,\aa}_{\em,J,\D}$ with $\aa\in\WW^{r+r'+1}$ subject
to (i),(ii),(iii). We set
$$N_\aa=|\{k\in[1,r],a_k>s_ka_k\}|+|\{k\in[r+1,r+r'],a_k>a_ks'_{r+r'+1-k}\}|,$$
$$\align&\ct_\aa=\{i\in\ct;a_{i-1}=a_i<s_ia_i\},\\&
\ct'_\aa=\{j\in\ct';a_{r+r'-j}=a_{r+r'-j+1}<a_{r+r'-j+1}s'_j\}.\endalign$$

\proclaim{Lemma 32.4} Assume that $\aa$ satisfies 32.3(i),(ii),(iii). 

(a) $\uZ^{\ww,\aa}_{\em,J,\D}$ is non-empty if and only if 
$i\in[1,r],a_{i-1}=a_i<s_ia_i\imp i\in\ct$, and
$j\in[1,r'],a_{-j+r+r'}=a_{-j+1+r+r'}<a_{-j+1+r+r'}s'_j\imp j\in\ct'$. 

(b) If $(B_0,B_1,\do,B_r,B'_0,B'_1,\do,B'_{r'},B,xU_{J,B_0})\in
\uZ^{\ww,\aa}_{\em,J,\D}$ then $B_{i-1}=B_i$ for any $i\in\ct_\aa$ and $B'_{j-1}=B'_j$ 
for any $j\in\ct'_\aa$. 

(c) Let ${}^0\uZ^{\ww,\aa}_{\em,J,\D}$ be the subset of $\uZ^{\ww,\aa}_{\em,J,\D}$ 
defined by the following conditions: for $i\in[1,r]$ we have $B_{i-1}=B_i$ if and only 
if $i\in\ct_\aa$; for $j\in[1,r']$ we have $B'_{j-1}=B'_j$ if and only if 
$j\in\ct'_\aa$. If $\uZ^{\ww,\aa}_{\em,J,\D}\ne\em$ then it is smooth, irreducible and 
${}^0\uZ^{\ww,\aa}_{\em,J,\D}$ is open dense in $\uZ^{\ww,\aa}_{\em,J,\D}$.

(d) If $a_k=a_{k-1}$ for some $k\in[1,r]$ with $k\n\ct$ or for some $k\in[r+1,r+r']$ 
with $r+r'+1-k\n\ct'$, then $\p_{\ww,\aa!}\bcl_\aa=0$.

(e) If $a_k\ne a_{k-1}$ for any $k\in[1,r]$ with $k\n\ct$ and for any 
$k\in[r+1,r+r']$ with $r+r'+1-k\n\ct'$, then $a_0a'\un\D\in\Wb_\cl$; moreover, 
$\uZ^{\ww,\aa}_{\em,J,\D}$ is an iterated affine space bundle over 
$Z^{(a_0,b)}_{\em,J,\D}$ with fibres of dimension $N_\aa$ and 
$\p_{\ww,\aa!}\bcl_\aa=K^{(a_0,a'),\cl}_{J,\D}[[-N_\aa]]$.
\endproclaim
We prove (d),(e) by induction on $r+r'$. If $r+r'=0$ then $\ww=(a_0,a'),\aa=\{a_0\}$ 
and 
$$\align\uZ^{\ww,\aa}_{\em,J,\D}=Z^\ww_{\em,J,\D}&=\{(B_0,B'_0,B,xU_{J,B_0});
B_0\in\cb,B'_0\in\cb,B\in\cb,x\in\D,\\&\po(B_0,B'_0)=a_0,\po(B'_0,B)=a',xB_0x\i=B\}.
\endalign$$
Hence $\p_{\ww,\aa!}\bcl_\aa=K^{\ww,\cl}_{J,\D}$. Now assume that $r+r'\ge 1$. Assume 
first that $r'\ge 1$. Let 

$\ww'=(s_1,s_2,\do,s_r,a_{r+r'-1},s'_{r'-1},s'_{r'-2},\do,s'_1,a')$,

$[\ww']=s_1s_2\do s_ra_{r+r'-1}s'_{r'-1}s'_{r'-2}\do s'_1a'$,

$\aa'=(a_0,a_1,\do,a_r,a_{r+1},\do,a_{r+r'-1})$,
$$\align&Y=\{(B_0,B_1,\do,B_r,B'_1,B'_2,\do,B'_{r'},B,xU_{J,B_0});B_i\in\cb
(i\in[0,r]),\\&B'_j\in\cb(j\in[1,r']),B\in\cb,x\in\D,\po(B_{i-1},B_i)\in\{s_i,1\}
(i\in\ct),\\&\po(B_{i-1},B_i)=s_i(i\in[1,r]-\ct),\po(B'_{j-1},B'_j)\in\{s'_j,1\}
(j\in\ct'\cap[2,r']),\\&\po(B'_{j-1},B'_j)=s'_j(j\in[2,r']-\ct'),\po(B_k,B'_{r'})
=a_k(k\in[0,r]),\\&\po(B_r,B'_{r+r'-k})=a_k(k\in[r,r+r'-1]),\po(B'_{r'},B)=a',
xB_0x\i=B'\}.\endalign$$
Define $\p^Y:Y@>>>Z_{J,\D}$, $f:\uZ^{\ww,\aa}_{\em,J,\D}@>>>Y$ by 
$$\p^Y:(B_0,B_1,\do,B_r,B'_1,B'_2,\do,B'_{r'},B,xU_{J,B_0})\m
(Q_{J,B_0},Q_{\e_\D(J),B},xU_{J,B_0}),$$
$$\align&f:(B_0,B_1,\do,B_r,B'_0,B'_1,\do,B'_{r'},B,xU_{J,B_0})\\&\m
(B_0,B_1,\do,B_r,B'_1,B'_2,\do,B'_{r'},xU_{J,B_0}).\endalign$$
The fibre of $f$ at $(B_0,B_1,\do,B_r,B'_1,B'_2,\do,B'_{r'},xU_{J,B_0})\in Y$ may be 
identified with
$$\align\{B'_0\in\cb;&\po(B'_0,B'_1)\in\{s'_1,1\}\text{ if }1\in\ct',\po(B'_0,B'_1)
=s'_1\text{ if }1\n\ct',\\&\po(B_r,B'_0)=a_{r+r'}\}.\endalign$$
If $a_{r+r'}s'_1=a_{r+r'-1}$ then $[\ww']\un\D=[\ww]\un\D\in\Wb_\cl$ and the sets 
analogous to $\ct,\ct'$ (for $\ww'$ instead of $\ww$) are $\ct,\ct'\cap[2,r']$. 
Moreover, $\uZ^{\ww',\aa'}_{\em,J,\D}=Y$ and $f$ is an isomorphism if 
$a_{r+r'}<a_{r+r'}s'_1$ and an affine line bundle if $a_{r+r'}>a_{r+r'}s'_1$. In both 
cases, $\bcl_\aa=f^*(\bcl_{\aa'})$. Hence 
$$f_!\bcl_\aa=\bcl_{\aa'},\p_{\ww,\aa!}\bcl_\aa=\p_{\ww',\aa'!}f_!\bcl_\aa= 
\p_{\ww',\aa'!}\bcl_{\aa'}\text{ if }a_{r+r'}<a_{r+r'}s'_1,$$ 
$$f_!\bcl_\aa=\bcl_{\aa'}[[-1]],\p_{\ww,\aa!}\bcl_\aa=\p_{\ww',\aa'!}f_!\bcl_\aa= 
\p_{\ww',\aa'!}\bcl_{\aa'}[[-1]]\text{ if }a_{r+r'}>a_{r+r'}s'_1;$$
the desired result follows from the induction hypothesis.

If $a_{r+r'-1}=a_{r+r'}<a_{r+r'}s'_1$ and $1\n\ct'$ then $\uZ^{\ww,\aa}_{\em,J,\D}
=\em$. If $a_{r+r'-1}=a_{r+r'}>a_{r+r'}s'_1$ and $1\n\ct'$ then $f$ is a $\kk^*$-bundle
and $f_!\bcl_\aa=0$ (this can be deduced from 28.11) hence 
$\p_{\ww,\aa!}\bcl_\aa=\p^Y_!f_!\bcl_\aa=0$.

Assume now that $a_{r+r'-1}=a_{r+r'}$ and $1\in\ct'$. Then $[\ww']\un\D\in\Wb_\cl$ and 
the sets analogous to $\ct,\ct'$ (for $\ww'$ instead of $\ww$) are 
$\ct,\ct'\cap[2,r']$. Moreover, $\uZ^{\ww',\aa'}_{\em,J,\D}=Y$; also, $f$ is an 
isomorphism if $a_{r+r'}<a_{r+r'}s'_1$ and an affine line bundle if 
$a_{r+r'}>a_{r+r'}s'_1$. In both cases, $\bcl_\aa=f^*(\bcl_{\aa'})$. Hence 
$$f_!\bcl_\aa=\bcl_{\aa'},\p_{\ww,\aa!}\bcl_\aa=\p_{\ww',\aa'!}f_!\bcl_\aa=
\p_{\ww',\aa'!}\bcl_{\aa'}\text{ if }a_{r+r'}<a_{r+r'}s'_1,$$ 
$$f_!\bcl_\aa=\bcl_{\aa'}[[-1]],\p_{\ww,\aa!}\bcl_\aa=\p_{\ww',\aa'!}f_!\bcl_\aa= 
\p_{\ww',\aa'!}\bcl_{\aa'}[[-1]]\text{ if }a_{r+r'}>a_{r+r'}s'_1;$$
the desired result follows from the induction hypothesis.

Assume next that $r'=0$. Then $r\ge 1$. Let 

$\ww''=(s_1,s_2,\do,s_{r-1},a_{r-1},a'),[\ww'']=s_1s_2\do s_{r-1}a_{r-1}a'$,

$\aa''=(a_0,a_1,\do,a_{r-1})$,
$$\align&Y_1=\{(B_0,B_1,\do,B_{r-1},B'_0,B,xU_{J,B_0});B_i\in\cb(i\in[0,r-1]),
B'_0\in\cb,\\&B\in\cb,x\in\D,\po(B_{i-1},B_i)\in\{s_i,1\}(i\in\ct\cap[1,r-1]),\\&
\po(B_{i-1},B_i)=s_i(i\in[1,r-1]-\ct),
\po(B_k,B'_0)=a_k(k\in[0,r-1]),\\&\po(B'_0,B)=a',xB_0x\i=B\}.\endalign$$
Define $\p^{Y_1}:Y_1@>>>Z_{J,\D}$, $f_1:\uZ^{\ww,\aa}_{\em,J,\D}@>>>Y_1$ by 
$$\align&\p^{Y_1}:
(B_0,B_1,\do,B_{r-1},B'_0,B,xU_{J,B_0})\m(Q_{J,B_0},Q_{\e_\D(J),B},xU_{J,B_0}),\\&
f_1:(B_0,B_1,\do,B_r,B'_0,B,xU_{J,B_0})\m(B_0,B_1,\do,B_{r-1},B'_0,B,xU_{J,B_0}).
\endalign$$
The fibre of $f_1$ at $(B_0,B_1,\do,B_{r-1},B'_0,B,xU_{J,B_0})\in Y_1$ may be 
identified with
$$\align\{B_r\in\cb;&\po(B_{r-1},B_r)\in\{s_r,1\}\text{ if }r\in\ct,\po(B_{r-1},B_r)
=s_r\text{ if }r\n\ct,\\&\po(B_r,B'_0)=a_r\}.\endalign$$
If $s_ra_r=a_{r-1}$ then $[\ww'']\un\D=[\ww]\un\D\in\Wb_\cl$ and the set analogous to 
$\ct$ (for $\ww''$ instead of $\ww$) is $\ct\cap[1,r-1]$. Moreover, 
$\uZ^{\ww'',\aa''}_{\em,J,\D}=Y_1$ and $f_1$ is an isomorphism if $a_r<s_ra_r$ and an 
affine line bundle if $a_r>s_ra_r$. In both cases, $\bcl_\aa=f_1^*(\bcl_{\aa''})$. 
Hence 
$$f_{1!}\bcl_\aa=\bcl_{\aa''},\p_{\ww,\aa!}\bcl_\aa=\p_{\ww'',\aa''!}f_{1!}\bcl_\aa=
\p_{\ww'',\aa''!}\bcl_{\aa''}\text{ if }a_r<s_ra_r,$$  
$$f_{1!}\bcl_\aa=\bcl_{\aa''}[[-1]],\p_{\ww,\aa!}\bcl_\aa=\p_{\ww'',\aa''!}
f_{1!}\bcl_\aa=\p_{\ww'',\aa''!}\bcl_{\aa''}[[-1]]\text{ if }a_r>s_ra_r;$$
the desired result follows from the induction hypothesis.

If $a_{r-1}=a_r<s_ra_r$ and $r\n\ct$ then $\uZ^{\ww,\aa}_{\em,J,\D}=\em$. If 
$a_{r-1}=a_r>s_ra_r$ and $r\n\ct$ then $f_1$ is a $\kk^*$-bundle and $f_{1!}\bcl_\aa=0$
(this can be deduced from 28.11) hence 
$\p_{\ww,\aa!}\bcl_\aa=\p^{Y_1}_!f_{1!}\bcl_\aa=0$.

Assume now that $a_{r-1}=a_r$ and $r\in\ct$. Then $[\ww'']\un\D\in\Wb_\cl$ and the set 
analogous to $\ct$ (for $\ww''$ instead of $\ww$) is $\ct\cap[1,r-1]$. Moreover, 
$\uZ^{\ww'',\aa''}_{\em,J,\D}=Y_1$; also, $f_1$ is an isomorphism if $a_r<s_ra_r$ and 
an affine line bundle if $a_r>s_ra_r$. In both cases, $\bcl_\aa=f_1^*(\bcl_{\aa''})$. 
Hence 
$$f_{1!}\bcl_\aa=\bcl_{\aa''},\p_{\ww,\aa!}\bcl_\aa=\p_{\ww'',\aa''!}f_{1!}\bcl_\aa=
\p_{\ww'',\aa''!}\bcl_{\aa''}\text{ if }a_r<s_ra_r,$$ 
$$f_{1!}\bcl_\aa=\bcl_{\aa''}[[-1]],\p_{\ww,\aa!}\bcl_\aa=\p_{\ww'',\aa''!}
f_{1!}\bcl_\aa=\p_{\ww'',\aa''!}\bcl_{\aa''}[[-1]]\text{ if }a_r>s_ra_r;$$
the desired result follows from the induction hypothesis. This completes the proof of 
(d),(e). The previous inductive proof also yields (a),(b),(c). The lemma is proved.

\subhead 32.5\endsubhead
Let $J\sub\II$. Let $D,D',\D$ be three connected components of G with $\D=D'D$. We 
write $\e,\e'$ instead of $\e_D,\e_{D'}:\WW@>>>\WW$. We have a diagram
$$Z_{J,D}\T Z_{\e(J),D'}@<b_1<<Z_0@>b_2>>Z_{J,\D}$$
where 
$$\align Z_0&=\{(Q,Q',Q'',gU_Q,g'U_{Q'});Q\in\cp_J,Q'\in\cp_{\e(J)},
Q''\in\cp_{\e'\e(J)},\\&g\in D,g'\in D', gQg\i=Q',g'Q'g'{}\i=Q''\},\endalign$$
$b_1(Q,Q',Q'',gU_Q,g'U_{Q'})=((Q,Q',gU_Q),(Q',Q'',g'U_{Q'}))$, 

$b_2(Q,Q',Q'',gU_Q,g'U_{Q'})=(Q,Q'',g'gU_Q)$.
\nl
Define a functor ({\it convolution}) $\cd(Z_{J,D})\T\cd(Z_{\e(J),D'})@>>>\cd(Z_{J,\D})$
by 
$$K,K'\m K*K'=b_{2!}b_1^*(K\bxt K').$$
Let $\ss=(s_1,s_2,\do,s_r),\ss'=(s'_1,s'_2,\do,s'_{r'})$ be two sequences in $\II$ and 
let $\cl,\cl'\in\fs(\TT)$ be such that $s_1s_2\do s_r\uD\in\Wb_\cl$, 
$s'_1s'_2\do s'_{r'}\uD'\in\Wb_{\cl'}$. Then $\bK^{\ss,\cl}_{J,D}\in\cd(Z_{J,D})$,
$\bK^{\ss',\cl'}_{\e(J),D'}\in\cd(Z_{\e(J),D'})$ are defined (see 28.12) hence 
$\bK^{\ss,\cl}_{J,D}*\bK^{\ss',\cl'}_{\e(J),D'}\in\cd(Z_{J,\D})$ is defined. 

Let $n\in\NN^*_\kk$ be such that $\cl\in\fs_n,\cl'\in\fs_n$. Let $\l\in\ufs_n$ (resp. 
$\l'\in\ufs_n$) be the isomorphism class of $\cl$ (resp. $\cl'$). 

\proclaim{Theorem 32.6} (a) If $\bK^{\ss,\cl}_{J,D}*\bK^{\ss',\cl'}_{\e(J),D'}\ne 0$ 
then $y\l'=\uD\l$ for some $y\in\WW_{\e(J)}$.

(b) Let $A$ be a simple perverse sheaf on $Z_{J,\D}$. If 
$A\dsv\bK^{\ss,\cl}_{J,D}*\bK^{\ss',\cl'}_{\e(J),D'}$, then $A\in\hZ_{J,\D}^\cl$.

(c) Assume that $\kk,\FF_q,G,F$ are as in 31.7(b), that $A\in\hZ_{J,\D}$ and that 
$\z^A:H_n[\D]@>>>\ca$ is as in 31.7. Then
$$\align&\c^A_v(\bK^{\ss,\cl}_{J,D}*\bK^{\ss',\cl'}_{\e(J),D'})=(v^2-1)^{\dim\TT} 
v^{\dim G-l(w^0_\II w^0_J)}\\&\T\z^A(\sum\Sb y'\in\WW_J\\ y'\uD\i\l'=\l\eSb
v^{2l(w^0_Jy')}C^\ss_{\uD\l}[D]T_{y'}[D\i]C^{\ss'}_{\uD'\l'}[D'][D]T_{y'{}\i}).
\endalign$$
(Notation of 31.5, 31.6. We regard $\cl,\cl'$ as pure of weight $0$ and then 
$\bK^{\ss,\cl}_{J,D}$, $\bK^{\ss',\cl'}_{\e(J),D'}$ and their convolution naturally as 
mixed complexes.) 
\endproclaim
The proof is given in 32.7-32.19.

\subhead 32.7\endsubhead
With notation of 28.7, let
$$\align&V=\{(B_0,B_1,\do,B_r,B'_0,B'_1,\do,B'_{r'},gU_{J,B_0},g'U_{\e(J),B'_0});
B_i\in\cb(i\in[0,r]),\\&B'_j\in\cb(j\in[0,r']),g\in D,g'\in D',gB_0g\i=B_r,
g'B'_0g'{}\i=B'_{r'},\\&
\po(B_{i-1},B_i)\in\{1,s_i\}(i\in[1,r]),\po(B'_{j-1},B'_j)\in\{1,s'_j\}(j\in[1,r']),\\&
\po(B_r,B'_0)\in\WW_{\e(J)}\}.\endalign$$
Let $\bcl$ be the constructible sheaf on $\bZ^\ss_{\em,J,D}$ in 28.10 and let $\bcl'$ 
be the analogous constructible sheaf on $\bZ^{\ss'}_{\em,\e(J),D'}$. The inverse image 
of $\bcl\bxt\bcl'$
under the imbedding $f:V@>>>\bZ^\ss_{\em,J,D}\T\bZ^{\ss'}_{\em,\e(J),D'}$,
$$\align&(B_0,B_1,\do,B_r,B'_0,B'_1,\do,B'_{r'},gU_{J,B_0},g'U_{\e(J),B'_0})\\&\m
((B_0,B_1,\do,B_r,gU_{J,B_0}),(B'_0,B'_1,\do,B'_{r'},g'U_{\e(J),B'_0}))\endalign$$
is a constructible sheaf on $V$ denoted again by $\bcl\bxt\bcl'$. Define 
$\r:V@>>>Z_{J,\D}$ by
$$(B_0,B_1,\do,B_r,B'_0,B'_1,\do,B'_{r'},gU_{J,B_0},g'U_{\e(J),B'_0})\m
(Q_{J,B_0},Q_{\e'\e(J),B'_{r'}},xU_{J,B_0})$$
where $x=g'g\in\D$; this is meaningful since $g'U_{\e(J),B'_0}g=g'U_{\e(J),B_r}g
=g'gU_{J,B_0}$. We have a commutative diagram in which the left square is cartesian
$$\CD
\bZ^\ss_{\em,J,D}\T\bZ^{\ss'}_{\em,\e(J),D'}@<f<<V@>\r>>Z_{J,\D}\\
@VhVV                                             @Vh_0VV     @V1VV\\
Z_{J,D}\T Z_{\e(J),D'}@<b_1<<                   Z_0@>b_2>>Z_{J,\D}
\endCD$$
where 
$$\align&h:((B_0,B_1,\do,B_r,gU_{J,B_0}),(B'_0,B'_1,\do,B'_{r'},g'U_{\e(J),B'_0}))\\&\m
((Q_{J,B_0},Q_{\e(J),B_r},gU_{J,B_0}),(Q_{\e(J),B'_0},Q_{\e'\e(J),B'_{r'}},
g'U_{\e(J),B'_0})),\endalign$$
$$\align&h_0:(B_0,B_1,\do,B_r,B'_0,B'_1,\do,B'_{r'},gU_{J,B_0},g'U_{\e(J),B'_0})\\&\m
((Q_{J,B_0},Q_{\e(J),B_r},Q_{\e'\e(J),B'_{r'}},gU_{J,B_0})),g'U_{\e(J),B_r})).
\endalign$$
Using this commutative diagram and the definitions we have
$$\bK^{\ss,\cl}_{J,D}*\bK^{\ss',\cl'}_{\e(J),D'}=\r_!(\bcl\bxt\bcl').$$
Let
$$\ct=\{i\in[1,r];s_1s_2\do s_i\do s_2s_1\in\WW_\cl\},$$
$$\ct'=\{j\in[1,r'];s'_{r'}\do s'_{j+1}s'_js'_{j+1}\do s'_{r'}\in\e'(\WW_{\cl'})\}.$$
(Thus, $\ct=\cj_\ss$, $\ct'=\cj_{\ss'}$ with the notation of 28.9.) Let
$$\align&V'=\{(B_0,B_1,\do,B_r,B'_0,B'_1,\do,B'_{r'},gU_{J,B_0},g'U_{\e(J),B'_0});
B_i\in\cb(i\in[0,r]),\\&B'_j\in\cb(j\in[0,r']),g\in D,g'\in D',gB_0g\i=B_r,
g'B'_0g'{}\i=B'_{r'},\\&\po(B_{i-1},B_i)\in\{1,s_i\}(i\in\ct),\po(B_{i-1},B_i)=s_i
(i\in[1,r]-\ct),\\&\po(B'_{j-1},B'_j)\in\{1,s'_j\}(j\in\ct'),
\po(B'_{j-1},B'_j)=s'_j(j\in[1,r']-\ct'),\\&\po(B_r,B'_0)\in\WW_{\e(J)}\},\endalign$$
an open subset of $V$. Let $\r':V'@>>>Z_{J,\D}$ be the restriction of 
$\r:V@>>>Z_{J,\D}$. From 28.10 we see that $(\bcl\bxt\bcl')|_{V'}$ is a local system 
and $(\bcl\bxt\bcl')|_{V-V'}=0$. Hence
$$\r_!(V,\bcl\bxt\bcl')=\r'_!(\bcl\bxt\bcl').$$
For any 
$$\aa=(a_0,a_1,\do,a_{r+r'})\in\WW^{r+r'}\T\WW_{\e(J)},\tag a$$
let
$$\align V'_\aa&=\{(B_0,B_1,\do,B_r,B'_0,B'_1,\do,B'_{r'},gU_{J,B_0},
g'U_{\e(J),B'_0})\in V';\\&\po(B_k,B'_{r'})=a_k(k\in[0,r]),\po(B_r,B'_{r+r'-k})=a_k
(k\in[r,r+r'])\}.\endalign$$
Let $\r'_\aa:V'_\aa@>>>Z_{J,\D}$ be the restriction of $\r':V'@>>>Z_{J,\D}$. Then 
$V'=\cup_\aa V'_\aa$ is a partition of $V'$ with $V'_\aa$ locally closed in $V'$ for 
all $\aa$. 

\proclaim{Lemma 32.8} Let $\aa$ as in 32.7(a) be such that 
$(\uD\i)^*\cl\not\cong(a_{r+r'}\i)^*\cl'$. Then $\r'_\aa(\bcl\bxt\bcl')=0$.
\endproclaim
Define $R$, $\ti\r:R@>>>Z_{J,\D}$, $\p:V'_\aa@>>>R$ by
$$\align R&=\{(B_0,B_1,\do,B_r,B'_0,B'_1,\do,B'_{r'},xU_{J,B_0});B_i\in\cb(i\in[0,r]),
\\&B'_j\in\cb(j\in[0,r']),x\in\D,\po(B_r,B'_0)=a_{r+r'},
\po(B'_{r'},xB_0x\i)=\e'(a_{r+r'}\i)\},\endalign$$
$$\ti\r:(B_0,B_1,\do,B_r,B'_0,B'_1,\do,B'_{r'},xU_{J,B_0})\m
(Q_{J,B_0},xQ_{J,B_0},xU_{J,B_0}),$$
$$\align&\p:(B_0,B_1,\do,B_r,B'_0,B'_1,\do,B'_{r'},gU_{J,B_0},g'U_{\e(J),B'_0})\\&\m
(B_0,B_1,\do,B_r,B'_0,B'_1,\do,B'_{r'},g'gU_{J,B_0});\tag a\endalign$$
this is meaningful since $g'U_{\e(J),B'_0}g=g'U_{\e(J),B_r}g=g'gU_{J,B_0}$ and
$$\align&\po(B'_{r'},g'gB_0g\i g'{}\i)=\e'(\po(g'{}\i B'_{r'}g',
gB_0g\i))=\e'(\po(B'_0,B_r))\\&=\e'(a_{r+r'}\i).\endalign$$
Since $\r'_\aa=\ti\r\p$, it suffices to show that $\p_!(\bcl\bxt\bcl')=0$. Hence it 
suffices to show that, for any 
$\x=(B_0,B_1,\do,B_r,B'_0,B'_1,\do,B'_{r'},xU_{J,B_0})\in R$, we have 
$H^e_c(\p\i(\x),\bcl\bxt\bcl')=0$ for all $e$. We may assume that $\p\i(\x)\ne\em$. We 
may identify 
$$\p\i(\x)=\{gU_{J,B_0};g\in D,gB_0g\i=B_r,xg\i B'_0gx\i=B'_{r'}\}$$ 
in an obvious way. We may assume that $B_r=B^*,B'_0=\doa_{r+r'}B^*\doa_{r+r'}\i$ 
(notation of 28.5). Write $B_i=h_iB^*h_i\i,B'_j=h'_jB^*h'_j{}\i$ with 
$h_i,h'_j\in G^0,h_r=1,h'_0=\doa_{r+r'}$. Let 

$w=\po(B_0,B_1)\po(B_1,B_2)\do\po(B_{r-1},B_r)\in\WW$,

$w'=\po(B'_0,B'_1)\po(B'_1,B'_2)\do\po(B'_{r'-1},B'_{r'})\in\WW$.
\nl
Let $T,U^*$ be as in 28.5. Let 
$d\in N_DB^*\cap N_DT,d'\in N_{D'}B^*\cap N_{D'}T$. Then $T$ (a maximal torus of 
$B_r\cap B'_0$) acts freely on $\p_\aa\i(\x)$ by left multiplication and it suffices to
show that for any $T$-orbit $\th$ in $\p_\aa\i(\x)$, we have 
$H^e_c(\th,\bcl\bxt\bcl')=0$ for all $e$. Let $g_0U_{J,B_0}\in\th$. It suffices to show
that the inverse image of $\bcl\bxt\bcl'$ under $t\m tg_0U_{J,B_0}$ (a local system in 
$\fs(T)$) is $\not\cong\bbq$. Using 28.10 and the definitions in 28.5, 28.8, we see 
that this inverse image is just $b^*(\cl)\ot c^*(\cl')$ where $\cl,\cl'$ are regarded 
as local systems on $T$ and $b:T@>>>T,c:T@>>>T$ are given by
$$\align&b(t)=d\i\dw\i n_1n_2\do n_rn_0(n_0\i tn_0),\\&
c(t)=d'{}\i\dw'n'_1n'_2\do n'_{r'}n'_0\doa_{r+r'}\i t\i\doa_{r+r'}\endalign$$
where $n_i\in N_{G^0}T(i\in[1,r])$, $n'_j\in N_{G^0}T(j\in[1,r'])$ are  given by 
$h_{i-1}\i h_i\in U^*n_iU^*$, $h'_{j-1}{}\i h'_j\in U^*n'_jU^*$, and 
$n_0\in N_DB^*\cap N_DT$, $n'_0\in N_{D'}B^*\cap N_{D'}T$ are given by 
$g_0h_0\in U^*n_0$, $h'_{r'}{}\i g_0\i x\doa_{r+r'}\in U^*n'_0$. Since 
$d\i\dw\i n_1n_2\do n_rn_0\in T$, $d'{}\i\dw'{}\i n'_1n'_2\do n'_{r'}n'_0\in T$, 
$n_0\in dT$, we see using 28.1(a) that $b^*\cl=\Ad(d\i)^*\cl$, 
$c^*\cl'=\Ad(\doa_{r+r'}\i)^*\che\cl'$. It then suffices to show that 
$\Ad(d\i)^*\cl\ot\Ad(\doa_{r+r'}\i)^*\che\cl'\not\cong\bbq$. This follows from our 
assumption. The lemma is proved.

\subhead 32.9\endsubhead
Until the end of 32.12 we assume that $\aa$ (as in 32.7(a)) is such that 
$$(\uD\i)^*\cl\cong(a_{r+r'}\i)^*\cl'.\tag a$$
We set $\ua=a_{r+r'},\ua'=\e'(a_{r+r'}\i)$. We show that 
$$s_1s_2\do s_r\ua s'_1s'_2\do s'_{r'}\ua'\un\D\in\Wb_\cl.$$
Since $s_1s_2\do s_r\uD\in\Wb_\cl$, it suffices to show that
$\uD\i\ua s'_1s'_2\do s'_{r'}\uD'\ua\i\uD\in\Wb_\cl$. Since 
$s'_1s'_2\do s'_{r'}\uD'\in\Wb_{\cl'}=\Wb_{\che\cl'}$, it suffices to show that 
$\ua^*(\uD\i)^*\cl\cong\che\cl'$. This holds by our assumption. Let
$$\ww=(s_1,s_2,\do,s_r,\ua,s'_1,s'_2\do,s'_{r'},\ua').$$
Then $\uZ^{\ww,\aa}_{\em,J,\D}$ (see 32.3) is well defined in terms of $\ct,\ct'$ as in
32.7 (or equivalently as in 32.2). As in 32.3, $\cl$ gives rise to a local system 
$\bcl_\aa$ on $\uZ^{\ww,\aa}_{\em,J,\D}$.

\subhead 32.10\endsubhead
We preserve the setup of 32.9. Let $V''_\aa$ be the set of all 

$(B_0,B_1,\do,B_r,B'_0,B'_1,\do,B'_{r'},(U_{B_r}\cap U_{B'_0})g,xU_{J,B_0})$
\nl
where
$(B_0,B_1,\do,B_r,B'_0,B'_1,\do,B'_{r'},xB_0x\i,xU_{J,B_0})\in\uZ^{\ww,\aa}_{\em,J,\D}$
and $g\in D$ satisfies $g\i B_rg=B_0,g\i B'_0g=x\i B'_{r'}x$. (The last equation is 
meaningful. It suffices to show that if $u\in U_{J,B_0}$ then $ug\i B'_0gu\i=g\i B'_0g$
that is, $gug\i\in N_GB'_0$. We have 
$gug\i\in U_{\e_D(J),gB_0g\i}=U_{\e_D(J),B_r}=U_{\e_D(J),B'_0}\sub U_{B'_0}\sub 
N_GB'_0$.) Define $\et:V''_\aa@>>>\uZ^{\ww,\aa}_{\em,J,\D}$, $\k:V'_\aa@>>>V''_\aa$ by 
$$\align&\et:(B_0,B_1,\do,B_r,B'_0,B'_1,\do,B'_{r'},(U_{B_r}\cap U_{B'_0})g,xU_{J,B_0})
\\&\m(B_0,B_1,\do,B_r,B'_0,B'_1,\do,B'_{r'},xB_0x\i,xU_{J,B_0}),\endalign$$ 
$$\align&\k:(B_0,B_1,\do,B_r,B'_0,B'_1,\do,B'_{r'},gU_{J,B_0},g'U_{\e(J),B'_0})\\&\m
(B_0,B_1,\do,B_r,B'_0,B'_1,\do,B'_{r'},(U_{B_r}\cap U_{B'_0})g,g'gU_{J,B_0});\endalign
$$
$\k$ is well defined, by the argument following 32.8(a). Let 
$f_\aa=\et\k:V'_\aa@>>>\uZ^{\ww,\aa}_{\em,J,\D}$. Clearly,

(a) {\it $\k$ is an affine space bundle with fibres of dimension $l(w^0_{\e(J)}\ua)$.}
\nl
Now $\TT$ acts on $V''_\aa$ by
$$\align t:&(B_0,B_1,\do,B_r,B'_0,B'_1,\do,B'_{r'},(U_{B_r}\cap U_{B'_0})g,xU_{J,B_0})
\m\\&(B_0,B_1,\do,B_r,B'_0,B'_1,\do,B'_{r'},(U_{B_r}\cap U_{B'_0})y_tg,xU_{J,B_0})
\endalign$$
where $y_t\in (B_r\cap B'_0)/(U_{B_r}\cap U_{B'_0})$ is defined by the condition that 
its image in $B_r/U_{B_r}$ is the image of $t$ under $\TT@>\si>>B_r/U_{B_r}$. Then

(b) {\it $\et$ is a principal $\TT$-bundle.}
\nl
Let $\x=(B_0,B_1,\do,B_r,B'_0,B'_1,\do,B'_{r'},B,xU_{J,B_0})\in
\uZ^{\ww,\aa}_{\em,J,\D}$. Then $\et\i(\x)$ may be identified with
$$\{(U_{B_r}\cap U_{B'_0})g;g\in D;g\i B_rg=B_0,g\i B'_0g=x\i B'_{r'}x\}.$$
We show only that $\et\i(\x)\cong\TT$. It suffices to show that $\et\i(\x)\ne\em$. For
this it suffices to show that $\po(B'_0,B_r)=\e(\po(x\i B'_{r}x,B_0))$ or that
$\po(B'_0,B_r)=\e\e_\D\i(\po(B'_{r},B))$ or that $\ua\i=\e'{}\i(\ua')$ which is clear.

\subhead 32.11\endsubhead
We preserve the setup of 32.9. Let $w$ be the product of the sequence $s_1,s_2,\do,s_r$
in which the factors $s_i$ with $i\in\ct_\aa$ are replaced by $1$. Let $w'$ be the 
product of the sequence $s'_1,s'_2,\do,s'_{r'}$ in which the factors $s'_j$ with 
$j\in\ct'_\aa$ are replaced by $1$, ($\ct_\aa,\ct'_\aa$ as in 32.3). Let $B^*,U^*,T$ be
as in 28.5. Let $d\in N_DB^*\cap N_DT,d'\in N_{D'}B^*\cap N_{D'}T$. We have a 
commutative diagram     
$$\CD
{}^0V'_\aa@<f_1<<\tV'@>f_2>>T\T T\\
@Vf_3VV          @Vf_4VV     @Vf_5VV\\
{}^0\uZ^{\ww,\aa}_{\em,J,\D}@<f_6<<\tZ @>f_7>> T   
\endCD$$
Here
$$\align&{}^0V'_\aa=\{(B_0,B_1,\do,B_r,B'_0,B'_1,\do,B'_{r'},gU_{J,B_0},
g'U_{\e(J),B'_0})\in V'_\aa;\\&B_{i-1}=B_i(i\in\ct_\aa),
\po(B_{i-1},B_i)=s_i(i\in[1,r]-\ct_\aa),B'_{j-1}=B'_j(j\in\ct'_\aa),\\&
\po(B'_{j-1},B'_j)=s'_j(i\in[1,r']-\ct'_\aa)\},\endalign$$
$$\align&\tV'=\{(h_0,h_1,\do,h_r,h'_0,h'_1,\do,h'_{r'},g,g')\in(G^0)^{r+r'+2}\T D\T D';
\\&h_{i-1}\i h_i\in B^*(i\in\ct_\aa),h_{i-1}\i h_i\in B^*\ds_iB^*(i\in[1,r]-\ct_\aa),
h'_{j-1}{}\i h'_j\in B^*(j\in\ct'_\aa),\\&h'_{j-1}{}\i h'_j\in B^*\ds'_jB^*
(j\in[1,r']-\ct'_\aa),h_r\i gh_0=d,h'_{r'}{}\i g'h'_0=d',\\&h_k\i h'_{r'}\in 
B^*\doa_kB^*(k\in[0,r-1]),
h_r\i h'_{r+r'-k}\in B^*\doa_kB^*(k\in[r,r+r'-1]),\\&h_r\i h'_0=\dua\},\endalign$$
$$\align&\tZ=\{(h_0,h_1,\do,h_r,h'_0,h'_1,\do,h'_{r'},x)\in(G^0)^{r+r'+2}\T D'D;\\&
h_{i-1}\i h_i\in B^*(i\in\ct_\aa),h_{i-1}\i h_i\in B^*\ds_iB^*(i\in[1,r]-\ct_\aa),
h'_{j-1}{}\i h'_j\in B^*(j\in\ct'_\aa),\\&h'_{j-1}{}\i h'_j\in B^*\ds'_jB^*
(j\in[1,r']-\ct'_\aa),h_k\i h'_{r'}\in B^*\doa_kB^*(k\in[0,r-1]),\\&
h_r\i h'_{r+r'-k}\in B^*\doa_kB^*(k\in[r,r+r'-1]),h_r\i h'_0=\dua, \\&
h'_{r'}{}\i xh_0\in B^*d'\dua\i d'{}\i B^*\},\endalign$$        
$$\align&f_1(h_0,h_1,\do,h_r,h'_0,h'_1,\do,h'_{r'},g,g')=
(h_0B^*h_0\i,\do,h_rB^*h_r\i,\\&h'_0B^*h'_0{}\i,\do,h'_{r'}B^*h'_{r'}{}\i,
gU_{J,h_0B^*h_0\i},g'U_{\e_D(J),h'_0B^*h'_0{}\i}),\endalign$$
$$\align&f_2(h_0,h_1,\do,h_r,h'_0,h'_1,\do,h'_{r'},g,g')\\&
=(d\i\dw\i n_1n_2\do n_rd,d'{}\i\dw'{}\i n'_1n'_2\do n'_{r'}d')\endalign$$
with $n_i,n'_j\in N_{G^0}T$ given by $h_{i-1}\i h_i\in U^*n_iU^*(i\in[1,r])$,
$h'_{j-1}{}\i h'_j\in U^*n'_jU^*(j\in[1,r'])$,

$f_3$ is the restriction of $f_\aa:V'_\aa@>>>\uZ^{\ww,\aa}_{\em,J,\D}$, see 32.10,
$$f_4(h_0,h_1,\do,h_r,h'_0,h'_1,\do,h'_{r'},g,g')=(h_0,h_1,\do,h_r,h'_0,h'_1,\do,
h'_{r'},g'g),$$
$$f_5(t,\tit)=d\i\dua d'{}\i\dw'{}\i\dua\i dtd\i\dua\dw'd'\tit\dua\i d
=\Ad(d\i\dua)(\Ad(d'{}\i\dw'{}\i\dua\i d)(t)\tit),$$
$$\align&f_6(h_0,h_1,\do,h_r,h'_0,h'_1,\do,h'_{r'},x)\\&=(h_0B^*h_0\i,\do,h_rB^*h_r\i,
h'_0B^*h'_0{}\i,\do,h'_{r'}B^*h'_{r'}{}\i,xU_{J,h_0B^*h_0\i}),\endalign$$
$$\align&f_7(h_0,h_1,\do,h_r,h'_0,h'_1,\do,h'_{r'},x)\\&=
(d'd)\i d'\dua d'{}\i\dw'{}\i\dua\i\dw\i n_1n_2\do n_r\dua n'_1n'_2\do n'_{r'}md'd,
\endalign$$
where $n_i,n'_j$ are as in the definition of $f_2$, $m\in N_{G^0}T$ is given by 
$h'_{r'}{}\i xh_0d\i d'{}\i\in U^*mU^*$.

\proclaim{Lemma 32.12}We preserve the setup of 32.11. Let 
${}^0\bcl_\aa=\bcl_\aa|_{{}^0\uZ^{\ww,\aa}_{\em,J,\D}}$.

(a) We have $f_3^*({}^0\bcl_\aa)=(\bcl\bxt\bcl')|_{{}^0V'_\aa}$.

(b) We have $f_\aa^*\bcl_\aa=(\bcl\bxt\bcl')|_{V'_\aa}$.
\endproclaim
From the definitions we have $f_6^*({}^0\bcl_\aa)=f_7^*\cl$, 

$f_1^*((\bcl\bxt\bcl')|_{{}^0V'_\aa})=f_2^*(\cl\bxt\Ad(\dua)^*\Ad(d\i)^*\che\cl)$.
\nl
Since $f_1^*$ is a smooth morphism with connected fibres, it suffices to show that 
$f_1^*f_3^*({}^0\bcl_\aa)=f_1^*((\bcl\bxt\bcl')|_{{}^0V'_\aa})$, or that
$f_4^*f_6^*({}^0\bcl_\aa)=f_2^*(\cl\bxt\cl')$, or that
$f_4^*f_7^*\cl=f_2^*(\cl\bxt\cl')$, or that $f_2^*f_5^*\cl=f_2^*(\cl\bxt\cl')$. It 
suffices to show that $f_5^*\cl=\cl\bxt\cl'$.  

Define $f'_5:T\T T@>>>T$ by $f'_5(t,\tit)=t\tit$. Setting $E=\Ad(d\i\dua):T@>>>T$,
$E'=\Ad(d'{}\i\dw'{}\i\dua\i d):T@>>>T$, we have $f_5=Ef'_5 (E'\T 1)$ hence
$$f_5^*\cl=(E'\T 1)^*f'_5{}^*E^*\cl=(E'\T 1)^*(E^*\cl\bxt E^*\cl)
=(EE')^*\cl\bxt E^*\cl.$$
From our assumption we have $\cl'\cong E^*\cl$. Moreover, 
$\cl'\cong\Ad(d'{}\i\dw'{}\i)^*\cl'=(E'E)^*\cl'$. Hence
$E^*\cl\cong(E'E)^*E^*\cl=(EE'E)^*\cl=E^*(EE')^*\cl$. Since $E^*$ is faithful, it 
follows that $(EE')^*\cl\cong\cl$. Thus, $f_5^*\cl\cong\cl\bxt\cl'$. This proves (a).

We prove (b). We may assume that $V'_\aa\ne\em$. From (a) we see that 
$f_\aa^*\bcl_\aa$, $(\bcl\bxt\bcl')|_{V'_\aa}$ are two local systems on $V'_\aa$ with 
the same restriction to the subset ${}^0V'_\aa$. It then suffices to show that $V'_\aa$
is smooth, irreducible and ${}^0V'_\aa$ is open dense in $V'_\aa$. By 32.10(a),(b), 
$f_\aa$ is a fibration with connected smooth fibres and ${}^0V'_\aa$ is the inverse 
image under $f_\aa$ of ${}^0\uZ^{\ww,\aa}_{\em,J,\D}$. Hence it suffices to show that 
$\uZ^{\ww,\aa}_{\em,J,\D}$ is smooth, irreducible and ${}^0\uZ^{\ww,\aa}_{\em,J,\D}$ is
open dense in $\uZ^{\ww,\aa}_{\em,J,\D}$. This follows from 32.4(c). The lemma is 
proved.

\subhead 32.13\endsubhead
Let $\SS$ be the set of all $\aa=(a_0,a_1,\do,a_{r+r'})\in\WW^{r+r'}\T\WW_{\e(J)}$ such
that 

(a) $a_k\in\{a_{k-1},s_ka_{k-1}\}$ for $k\in[1,r]$,    

(b) $a_k\in\{a_{k-1},a_{k-1}s'_{r+r'+1-k}\}$ for $k\in[r+1,r+r']$,

(c) $(\uD\i)^*\cl\cong(a_{r+r'}\i)^*\cl'$,

(d) $i\in[1,r],a_{i-1}=a_i\imp i\in\ct$,

(e) $j\in[1,r'],a_{-j+r+r'}=a_{-j+1+r+r'}\imp j\in\ct'$. 

\proclaim{Lemma 32.14}If $\aa\in(\WW^{r+r'}\T\WW_{\e(J)})-\SS$ then 
$\r'_{\aa!}(\bcl\bxt\bcl')=0$.
\endproclaim
If $\aa$ does not satisfy 32.13(a) or 32.13(b), then $V'_\aa=\em$ and the result is 
trivial. If $\aa$ does not satisfy 32.13(c), the desired result follows from Lemma 
32.8. Assume now that $\aa$ satisfies 32.13(a)-(c) but it does not satisfy 32.13(d) or 
(e). Using $\r'_\aa=\p_{\ww,\aa}f_\aa$ and 32.12(b), we see that it suffices to show 
that $\p_{\ww,\aa!}f_{\aa!}(f_\aa^*\bcl_\aa)=0$, or that 
$\p_{\ww,\aa!}(\bcl_\aa\ot f_{\aa!}\bbq)=0$. It suffices to show that 
$\p_{\ww,\aa!}(\bcl_\aa\ot\ch^e(f_{\aa!}\bbq))=0$ for any $e$. By 32.10 we have 
$f_\aa=\et\k$ and $\k_!\bbq=\bbq[[-c]]$ where $c=l(w^0_{\e(J)}a_{r+r'})$. Hence 
$f_{\aa!}\bbq=\et_!\bbq[[-c]]$. Let $\rr=\dim\TT$. Since $\et$ is a principal 
$\TT$-bundle (see 32.10), the local system $\ch^e(\et_!\bbq)$ admits a filtration whose
associated graded is a direct sum of $\bin{\rr}{2\rr-e}$ copies of $\bbq(\rr-e)$. Since
$\p_{\ww,\aa!}\bcl_\aa=0$ (see 32.4(d)), we see that 
$\p_{\ww,\aa!}(\bcl_\aa\ot\ch^e(f_{\aa!}\bbq))=0$ for any $e$. The lemma is proved.

\subhead 32.15\endsubhead
We now make a short digression. Let $X$ be an algebraic variety over $\kk$. Let 
$C\in\cd(X)$ and let $\{C_n;n\in\ZZ\}$ be a sequence of objects in $\cd(X)$ such that 
$C_n=0$ for all but finitely many $n$. We shall write 
$$C\Bpq\{C_n;n\in\ZZ\}$$
if the following condition is satisfied: there exists a sequence $\{C'_n;n\in\ZZ\}$ of 
objects in $\cd(X)$ such that $C'_n=0$ for $n\ll 0$, $C'_n=C$ for $n\gg 0$ and 
distinguished triangles $(C'_{n-1},C'_n,C_n)$ for $n\in\ZZ$.

If $X,C,C_n$ are as above, $C\Bpq\{C_n;n\in\ZZ\}$ and $X_2@>f_2>>X@>f_1>>X_1$ are 
morphisms of algebraic varieties, we see from definitions that:

$f_{1!}C\Bpq\{f_{1!}C_n;n\in\ZZ\}$,

$f^*_2C\Bpq\{f^*_2C_n;n\in\ZZ\}$.
\nl
Assume now that $C\in\cd(X)$ and that $\{C^u;u\in\cu\}$ are objects of $\cd(X)$ indexed
by a finite set $\cu$. We shall write
$$C\Bpq\{C_u;u\in\cu\}$$
if the following condition is satisfied: there exists a bijection $\cu\lra[0,m]$ such 
that, setting $C_n=C^u$ if $u\lra n\in[0,m]$ and $C_n=0$ for $n\n[0,m]$, we have 
$C\Bpq\{C_n;n\in\ZZ\}$.

For example, if $C\in\cd(X)$, we have $C\Bpq\{{}^pH^nC[-n];n\in\ZZ\}$; in this case we 
can take $C'_n={}^p\t_{\le n}(C)$ (truncation, as in \cite{\BBD}).

Similarly, if $C\in\cd(X)$, we have $C\Bpq\{\ch^n(C)[-n];n\in\ZZ\}$.

As another example, assume that we are given a partition $X=\sqc_{u\in\cu}X^u$ with 
$\cu$ finite, where $X^u$ are locally closed subvarieties of $X$ such that for some 
bijection $\cu\lra[0,m]$, the union $X'_n=X_n\cup X_{n-1}\cup\do\cup X_0$ is open in 
$X$ for any $n\in[0,m]$ (we set $X_n=X^u$ for $u\lra n\in[0,m]$). For any $u\in\cu$ let
$j_u:X^u@>>>X$ be the inclusion and let $C_u=j_{u!}j_u^*C$. We have 
$C\Bpq\{C_u;u\in\cu\}$. Indeed, setting $C_n=C_u$ if $u\lra n\in[0,m]$ and $C_n=0$ for 
$n\n[0,m]$, we have $C\Bpq\{C_n;n\in\ZZ\}$. (We can take $C'_n=0$ for $n<0$, $C'_n=C$ 
for $n>m$, $C'_n=j'_{n!}j'_n{}^*C$ for $n\in[0,m]$, where $j'_n=X'_n@>>>X$ is the 
inclusion.)

\subhead 32.16\endsubhead
Assume that $\aa\in\SS$. As in the proof of 32.14 we have 
$\r'_{\aa!}(\bcl\bxt\bcl')=\p_{\ww,\aa!}(\bcl_\aa\ot f_{\aa!}\bbq)$ and
$f_{\aa!}\bbq=\et_!\bbq[[-c]]$. Moreover, 
$\et_!\bbq\Bpq\{\ch^e(\et_!\bbq)[-e];e\in\ZZ\}$ and for any $e$ we have 
$\ch^e(\et_!\bbq)\Bpq\{C^e_{e'};1\le e'\le\bin{\rr}{2\rr-e}\}$ where 
$C^e_{e'}=\bbq(\rr-e)$. It follows that
$$\r'_{\aa!}(\bcl\bxt\bcl')\Bpq\{\p_{\ww,\aa!}(\bcl_\aa\ot\ch^e(\et_!\bbq))[-e][[-c]];
e\in\ZZ\},\tag a$$
$$\p_{\ww,\aa!}(\bcl_\aa\ot\ch^e(\et_!\bbq)[-e][[-c]]\Bpq\{C'{}^e_{e'};1\le e'\le
\bin{\rr}{2\rr-e}\}\tag b$$
where             
$$C'{}^e_{e'}=\p_{\ww,\aa!}(\bcl_\aa)(\rr-e)[-e][[-c]]
=K^{(a_0,\e'(a_{r+r'}\i)),\cl}_{J,\D}(\rr-e)[-e][[-c]][[-N_\aa]]$$
(see 32.4(e)). By 32.4(e) we have
$$a_0\e'(a_{r+r'}\i)\un\D\in\Wb_\cl.\tag c$$
From (a),(b) we see that, if $A$ is a simple perverse sheaf on $Z_{J,\D}$ such that
$A\dsv\r'_{\aa!}(\bcl\bxt\bcl')$ then $A\dsv K^{(a_0,\e'(a_{r+r'}\i)),\cl}_{J,\D}$.

\subhead 32.17\endsubhead
From the partition $V'=\sqc_\aa V'_\aa$ we get as in 32.15
$$\r_!(\bcl\bxt\bcl')=\r'_!(\bcl\bxt\bcl')\Bpq\{\r'_\aa(\bcl\bxt\bcl');\aa\in\SS\}
\tag a$$
(by 32.14 we can omit the $\aa\n\SS$). Thus, if $A$ is a simple perverse sheaf on 
$Z_{J,D}$ such that $A\dsv\r_!(\bcl\bxt\bcl')$, then for some $\aa\in\SS$ we have 
$A\dsv\r'_{\aa!}(\bcl\bxt\bcl')$ hence, by 32.16, 
$A\dsv K^{(a_0,\e'(a_{r+r'}\i)),\cl}_{J,\D}$, so that $A\in\hZ^\cl_{J,\D}$. We also see
that $\SS\ne\em$; in particular, $(\uD\i)^*\cl\cong y^*\cl'$ for some 
$\ua\in\WW_{\e(J)}$ (see 32.13(c)). Since 
$\bK^{\ss,\cl}_{J,D}*\bK^{\ss',\cl'}_{\e(J),D'}=\r_!(\bcl\bxt\bcl')$ (see 32.7) we see 
that 32.6(b) holds. We also see that 32.6(a) holds since, under the assumption of
32.6(a), we can find an $A$ as above.

\subhead 32.18\endsubhead
In this and the next subsection we place ourselves in the setup of 32.6(c). Then 
$V,V',V'_\aa$ are defined over $\FF_q$ and we can regard $\r'_!(\bcl\bxt\bcl')$ and 
$\r'_{\aa!}(\bcl\bxt\bcl')$ (for any $\aa$) as mixed complexes on $Z_{J,\D}$. Using 
32.17(a), 32.16(a),(b) (or rather their variant in the mixed category) and 31.7(c),(e) 
we see that, with the notation of 31.6, we have
$$\align&\c^A_v(\r_!(\bcl\bxt\bcl'))=\sum_{\aa\in\SS}\c^A_v(\r'_\aa(\bcl\bxt\bcl'))\\&
=\sum_{\aa\in\SS}\sum_{e\in\ZZ}(-1)^ev^{2N_\aa+2l(w^0_{\e(J)}a_{r+r'})-2\rr+2e}
\bin{\rr}{2\rr-e}\c^A_v(K^{(a_0,\e'(a_{r+r'}\i)),\cl}_{J,\D})\\&=(v^2-1)^\rr
\sum_{\aa\in\SS}v^{2N_\aa-2l(w^0_{\e(J)}a_{r+r'})}
\c^A_v(K^{(a_0,\e'(a_{r+r'}\i)),\cl}_{J,\D})
\\&=(v^2-1)^\rr\sum_{\aa\in\SS}v^{2N_\aa
-2l(w^0_{\e(J)}a_{r+r'})}v^{\dim G-l(w^0_\II w^0_J)}\z^A(T_{a_0}T_{\e'(a_{r+r'}\i}
1_{\un\D\l}[\D]).\tag a\endalign$$

\subhead 32.19\endsubhead
Let $h\m h^\flat$ be the antiautomorphism of the algebra $H_n$ defined by
$T_u\m T_{u\i}$ for $u\in\WW$, $1_\l\m 1_\l$ for $\l\in\ufs_n$. We have
$(C^{\tss}_\l)^\flat=C^\ss_{s_r\do s_2s_1\l}$ where $\tss=(s_r,s_{r-1},\do,s_1)$. The 
following identity in the algebra $H_n$ (see 31.2) is a special case of one in 31.11:
$$T_yC^{\ss'}_{\l_1}=\sum_{\yy'}v^{2\d'(\yy')}T_{y'_{r'}}1_{\l_1};\tag a$$
here $y\in\WW,\l_1\in\ufs_n$, the sum is taken over all sequences 
$\yy'=(y'_0,y'_1,\do,y'_{r'})$ in $\WW$ such that 

$y=y'_0$,
 
$y'_i\in\{y'_{i-1},y'_{i-1}s'_i\}$ for $i\in[1,r']$,
 
$i\in[1,r'],y'_{i-1}=y'_i\imp s'_i\in\WW_{s'_{i+1}\do s'_{r'}\l_1}$;
\nl
moreover, $\d'(\yy')=|\{i\in[1,r'];y'_{i-1}s'_i<y'_{i-1}\}|$. Similarly, we have
$$C^\ss_{s_r\do s_2s_1\l_2}T_y=\sum_\yy v^{2\d(\yy)}1_{\l_2}T_{y_0};\tag b$$
here $\l_2\in\ufs_n$, the sum is taken over all sequences $\yy=(y_0,y_1,\do,y_r)$ in 
$\WW$ such that 

$y=y_r$,

$y_i\in\{y_{i-1},s_iy_{i-1}\}$ for $i\in[1,r]$,
 
$i\in[1,r],y_{i-1}=y_i\imp s_i\in\WW_{s_{i-1}\do s_1\l_2}$;
\nl
moreover, $\d(\yy)=|\{i\in[1,r];s_iy_i<y_i\}|$. This can be deduced from (a) using the 
involution $h\m h^\flat$.

Combining (a),(b) we obtain (for $\l_1,\l_2\in\ufs_n$) the identity
$$C^\ss_{s_r\do s_2s_1\l_2}T_yC^{\ss'}_{\l_1}
=\sum_{\yy,\yy'}v^{2\d(\yy)+2\dd(\yy')}1_{\l_2}T_{y_0}1_{\l_1};$$
the sum is taken over the pairs $\yy=(y_0,y_1,\do,y_r)$, $\yy'=(y'_0,y'_1,\do,y'_{r'})$
of sequences in $\WW$ such that 

$y=y'_0,y'_{r'}=y_r$, 

$y_i\in\{y_{i-1},s_iy_{i-1}\}$ for $i\in[1,r]$, 

$y'_i\in\{y'_{i-1},y'_{i-1}s'_i\}$ for $i\in[1,r']$,

$i\in[1,r],y_{i-1}=y_i\imp s_i\in\WW_{s_{i-1}\do s_1\l_2}$,

$i\in[1,r'],y'_{i-1}=y'_i\imp s'_i\in\WW_{s'_{i+1}\do s'_{r'}\l_1}$.
\nl
We have $s_r\do s_2s_1\l=\uD\l$. Take $\l_1=\uD'\l',\l_2=\l$. Take $y\in\WW_{\e(J)}$ 
such that $y\l'=\uD\l$. We replace $(\yy,\yy')$ by $\aa=(a_0,a_1,\do,a_{r+r'})$ where 
$a_k=y_k$ for $k\in[0,r]$, $a_k=y'_{r+r'-k}$ for $k\in[r,r+r']$. Then 
$\d(\yy)+\d'(\yy')=N_\aa$.  We obtain
$$C^\ss_{\uD\l}T_yC^{\ss'}_{\uD'\l'}=\sum_\aa v^{2N_\aa}1_\l T_{a_0}1_{\uD'\l'};$$
the sum is over all $\aa=(a_0,a_1,\do,a_{r+r'})\in\WW^{r+r'}\T\WW_{\e(J)}$ such that

$y=a_{r+r'}$, 

$a_i\in\{a_{i-1},s_ia_{i-1}\}$ for $i\in[1,r]$, 

$a_{r+r'-i}\in\{a_{r+r'-i+1},a_{r+r'-i+1}s'_i\}$ for $i\in[1,r']$,

$i\in[1,r],a_{i-1}=a_i\imp s_i\in\WW_{s_{i-1}\do s_1\l}$,

$i\in[1,r'],a_{r+r'-i+1}=a_{r+r'-i}\imp s'_i\in\WW_{s'_{i+1}\do s'_{r'}\uD'\l'}$.
\nl
Equivalently,
$$C^\ss_{\uD\l}T_yC^{\ss'}_{\uD'\l'}
=\sum_{\aa\in\SS;a_{r+r'}=y}v^{2N_\aa}1_{\l}T_{a_0}1_{\uD'\l'}.\tag c$$
For each $\aa$ in the sum we have $a_0\e'(y\i)\un\D\l=\l$ (see 32.16(a)); combining 
this with $y\l'=\uD\l$ we see that $a_0\uD'\l'=\l$, hence 
$1_{\l}T_{a_0}1_{\uD'\l'}=T_{a_0}1_{\uD'\l'}$. We introduce this in (c), then multiply 
both sides of (c) on the right by 
$$(v^2-1)^\rr v^{\dim G-l(w^0_\II w^0_J)}v^{-2l(y)}T_{\e'(y\i)}1_{\un\D\l}[\D]$$ 
and sum over $y$. We obtain
$$\align&(v^2-1)^\rr v^{\dim G-l(w^0_\II w^0_J)}\sum_{y\in\WW_{\e(J)};y\l'=\uD\l}
v^{-2l(y)}C^\ss_{\uD\l}T_yC^{\ss'}_{\uD'\l'}T_{\e'(y\i)}[\D]\\&=(v^2-1)^\rr 
v^{\dim G-l(w^0_\II w^0_J)}\sum_{\aa\in\SS}v^{2N_\aa-2l(a_{r+r'}}T_{a_0}
T_{\e'(a_{r+r'}\i)}1_{\un\D\l}[\D].\endalign$$
We apply $\z^A_v$ (see 31.7) to both sides and use 32.18(a). We obtain
$$\align&\c^A_v(\r_!(\bcl\bxt\bcl'))=\\&(v^2-1)^\rr v^{\dim G-l(w^0_\II w^0_J)}
\z^A(\sum\Sb y\in\WW_{\e(J)}\\y\l'=\uD\l\eSb v^{2l(w^0_{\e(J)}y)}C^\ss_{\uD\l}T_y
C^{\ss'}_{\uD'\l'}T_{\e'(y\i)}[\D]).\endalign$$
We substitute $y=\e(y'),y'\in\WW_J$ that is, $T_y=[D]T_{y'}[D]\i$. Since \lb
$\c^A_v(\bK^{\ss,\cl}_{J,D}*\bK^{\ss',\cl'}_{\e(J),D'})=\c^A_v(\r_!(\bcl\bxt\bcl'))$, 
32.6(c) follows. This completes the proof of Theorem 32.6.

\subhead 32.20\endsubhead
Let $\cd^{cs}(Z_{J,D})$ (resp. $\cd^\cl(Z_{J,D})$ with $\cl\in\fs(\TT)$) be the 
subcategory of $\cd(Z_{J,D})$ whose objects are those $K\in\cd(Z_{J,D})$ such that for 
any $j$, any simple subquotient of ${}^pH^jK$ is in $\hZ_{J,D}$ (resp. in 
$\hZ^\cl_{J,D}$). We have the following result.

\proclaim{Corollary 32.21}(a) If $K\in\cd^{cs}(Z_{J,D}),K'\in\cd^{cs}(Z_{\e(J),D'})$, 
then $K*K'\in\cd^{cs}(Z_{J,\D})$.

(b) If $\cl\in\fs(\TT)$, $K\in\cd^\cl(Z_{J,D})$, $K'\in\cd^{cs}(Z_{\e(J),D'})$, then
$K*K'\in\cd^\cl(Z_{J,\D})$.
\endproclaim
We prove (a). We may assume that $K\in\hZ_{J,D},K'\in\hZ_{\e(J),D'}$. We can find 
$\ss,\ss',\cl,\cl'$ as in 32.5 and $u,u'\in\ZZ$ such that $K[u]$ is a direct summand of
$\bK^{\ss,\cl}_{J,D}$ and $K'[u']$ is a direct summand of $\bK^{\ss',\cl'}_{\e(J),D'}$.
Then $K*K'[u+u']$ is a direct summand of 
$\bK^{\ss,\cl}_{J,D}*\bK^{\ss',\cl'}_{\e(J),D'}$ which, by 32.6(b), is in 
$\cd^\cl(Z_{J,\D})$. Hence $K*K'[u+u']$ is in $\cd^\cl(Z_{J,\D})$. This proves (a). The
same argument proves (b).

\subhead 32.22\endsubhead
If $E$ is a mixed $\bbq$-vector space (that is, a $\bbq$-vector space which, when 
regarded as a complex over a point, is a mixed complex) we set 

$\c_v(E)=\sum_j\dim(E_j)v^j\in\ca$
\nl
where $E_j$ is the pure subquotient of weight $j$ of $E$. We preserve the setup in 
32.5. Assume that $D'=D\i$ hence $\D=G^0$. Let $\SS$ be as in 32.13. 

Define an $\ca$-linear map $\Ph:H_n@>>>H_n$ by 
$\x\m\fa_D(C^\ss_{\uD\l}\x C^{\ss'}_{\uD'\l'})$ with $\fa_D$ as in 31.4. For any 
$y\in\WW_{\e(J)}$ such that $y\l'=\uD\l$ we have 
$$\Ph(T_y1_{\l'})=\sum_{\aa\in\SS;a_{r+r'}=y}v^{2N_\aa}1_{\uD\l}T_{\e(a_0)}1_{\l'}.$$
(See 32.19(c).) Define an $\ca$-linear map $\Th^J:H_n@>>>H_n$ by 
$\Th^J(T_w1_{\l_1})=T_w1_{\l_1}$ if $w\in\WW_J,\l_1\in\ufs_n$, $\Th^J(T_w1_{\l_1})=0$ 
if $w\in\WW-\WW_J,\l_1\in\ufs_n$. Replacing $J$ by $\e(J)$ we obtain an $\ca$-linear 
map $\Th^{\e(J)}:H_n@>>>H_n$. Define $\Ph':H_n@>>>H_n$ by 
$\Ph'(\x')=\Th^{\e(J)}\Ph(\x')$. Since $H_n$ is a free $\ca$-module and $\Ph'$ is 
$\ca$-linear, $\tr(\Ph',H_n)\in\ca$ is well defined. From the definitions we have
$$\tr(\Ph',H_n)=\sum_{\aa\in\SS;a_{r+r'}=\e(a_0)}v^{2N_\aa}=
\sum_{\aa\in\SS_0}v^{2N_\aa},\tag a$$
where $\SS_0=\{\aa\in\SS;a_0=\e'(a_{r+r'})\}$.

Define an $\ca$-linear map $\Ph'':H_n@>>>H_n$ by 
$\x\m\Th^J(C^{\tss'}_{\l'{}\i}\fa_D(\x)C^{\tss}_{\l\i})$ where 
$\tss=(s_r,\do,s_2,s_1)$, $\tss'=(s'_{r'},\do,s'_2,s'_1)$. Let 
$\part:Z_{J,D}@>>>Z_{\e(J),D'}$ be as in 28.19. Then
$\bK^{\ss,\cl}_{J,D}\ot\part^*\bK^{\ss',\cl'}_{\e(J),D'}\in\cd(Z_{J,D})$ is well 
defined. Let
$$\mu(G^0)=(v^2-1)^\rr\sum_{w\in\WW}v^{2l(w)}\in\ca.$$
The following result is an application of (the proof of) Theorem 32.6.

\proclaim{Corollary 32.23}Assume that $\kk,\FF_q,G,F$ are as in 31.7(b). Then
$$\align&\sum_z(-1)^z\c_v(H^z_c(Z_{J,D},\bK^{\ss,\cl}_{J,D}\ot
\part^*\bK^{\ss',\cl'}_{\e(J),D\i})=v^{2l(w^0_J)}\mu(G^0)\tr(\Ph',H_n)\\&
=v^{2l(w^0_J)}\mu(G^0)\tr(\Ph'',H_n).\tag a\endalign$$
\endproclaim
Let $\fZ=\{(Q,Q',xU_Q)\in Z_{J,\D};Q=Q',x\in U_Q\}$, let $\io:\fZ@>>>Z_{J,\D}$ be the 
inclusion and let $p:\fZ@>>>\text{point}$ be the obvious map.
From the definitions, for any $A\in\cd(\hZ_{J,D}),A'\in\cd(Z_{\e(J),D'})$ we have
$$H^z_c(\text{point},p_!\io^*(A*A'))=H^z_c(Z_{J,D},A\ot\part^*(A'))\tag b$$
for any $z\in\ZZ$. In particular,
$$H^z_c(Z_{J,D},\bK^{\ss,\cl}_{J,D}\ot\part^*\bK^{\ss',\cl'}_{\e(J),D'})=
H^z_c(\text{point},p_!\io^*\r_!(\bcl\bxt\bcl')).$$
Applying $p_!\io^*$ to 32.17(a) gives
$$p_!\io^*\r_!(\bcl\bxt\bcl')\Bpq\{p_!\io^*\r'_{\aa!}(\bcl\bxt\bcl');\aa\in\SS\}.\tag c
$$
Let $\aa\in\SS$. Applying $p_!\io^*$ to 32.16(a), 32.16(b) gives 
$$p_!\io^*\r'_{\aa!}(\bcl\bxt\bcl')\Bpq\{K_e;e\in\ZZ\},\tag d$$
$$K_e\Bpq\{K'{}^e_{e'};1\le e'\le\bin{\rr}{2\rr-e}\}\tag e$$
where $K_e,K'{}^e_{e'}\in\cd(\text{point})$ and 
$K'{}^e_{e'}=p_!\io^*K^{(a_0,\e'(a_{r+r'}\i)),\cl}_{J,\D}(\rr-e)[-e][[-c-N_\aa]]$
(notation of 32.16). Let 
$$X_\aa=\{(B_0,B'_0)\in\cb\T\cb;\po(B_0,B'_0)=a_0,\po(B'_0,B_0)=\e'(a_{r+r'}\i)\}$$
and let $\o:X_\aa@>>>\text{point}$ be the obvious map. From the definitions we see that
$p_!\io^*K^{(a_0,\e'(a_{r+r'}\i)),\cl}_{J,\D}=\o_!\bbq$. If $a_0\ne\e'(a_{r+r'})$ then
$X_\aa=\em$ hence $\o_!\bbq=0$; if $a_0=\e'(a_{r+r'})$ that is, $\aa\in\SS_0$, then
$$H^z(\text{point},\o_!\bbq)=\op_{w\in\WW;2l(w)=z}\bbq(-z-l(a_0)).\tag f$$
Using  (c),(d),(e),(f) (or rather their variant in the mixed category) we see that
$$\align&\sum_z(-1)^z\c_v(H^z(\text{point},p_!\io^*\r_!(\bcl\bxt\bcl')))\\&=\sum_{\aa
\in\SS_0,e\in\ZZ}(-1)^ev^{2N_\aa+2l(w^0_{\e(J)}a_{r+r'})-2\rr+2e}\bin{\rr}{2\rr-e}
\sum_{w\in\WW}v^{2l(w)+2l(a_0)}\\&=v^{2l(w^0_J)}(v^2-1)^\rr\sum_{w\in\WW}v^{2l(w)}
\sum_{\aa\in\SS_0}v^{2N_\aa}=v^{2l(w^0_J)}\mu(G^0)\tr(\Ph',H_n).\endalign$$
It remains to show that $\tr(\Ph',H_n)=\tr(\Ph'',H_n)$. Define $\ca$-linear maps 
$\Ps',\Ps'',\O:H_n@>>>H_n$ by
$$\Ps'(\x)=C^\ss_{\uD\l}\x C^{\ss'}_{\uD'\l'},\Ps''(\x)=C^{\ti\ss'}_{\l'{}\i}\x 
C^{\ti\ss}_{\l\i},\O(T_w1_{\l_1})=1_{\l_1\i}T_{w\i}.$$
One checks that 

$\O\Ps'=\Ps''\O,\fa_D\Th^J=\Th^{\e(J)}\fa_D,\fa_D\Th^J\O=\O\fa_D\Th^J$,
$\Ph'=\Th^{\e(J)}\fa_D\Ps'$,

$\Ph''=\Th^J\Ps''\fa_D$.
\nl
Hence $\Ph'=\fa_D\Th^J\O\i\Ps''\O=\O\i\fa_D\Th^J\Ps''\O$ and
$$\align&\tr(\Ph',H_n)=\tr(\O\i\fa_D\Th^J\Ps''\O,H_n)=\tr(\fa_D\Th^J\Ps'',H_n)\\&
=\tr(\Th^J\Ps''\fa_D,H_n)=\tr(\Ph'',H_n),\endalign$$
as required. The corollary is proved.

\mpb 

For $\l\in\ufs(\TT)$ we set $\WW_\l=\WW_\cl$ where $\cl\in\fs(\TT)$ is in the
isomorphism class $\l$; this agrees with the definition in 31.2 when $\l\in\ufs_n$.

\proclaim{Corollary 32.24}Let $A\in\hZ_{J,D},\cl,\cl''\in\fs(\TT)$. Let $\l$ (resp. 
$\l''$) be the isomorphism class of $\cl$ (resp. $\cl''$). Let $\ss=(s_1,s_2,\do,s_r)$,
$\ss''=(s''_1,s''_2,\do,s''_{r'})$ be sequences in $\II$ such that 
$s_1s_2\do s_r\uD\l=\l$, $s''_1s''_2\do s''_{r'}\uD\l'=\l''$, 
$A\dsv\bK^{\ss,\cl}_{J,D}$ and $A\dsv\bK^{\ss'',\l''}_{J,D}$. Then there exist 
$b\in\WW_{\l''}$,  $a_0\in\WW_J$ such that
$$a_0(\l'')=\l,\qua s_1s_2\do s_r\uD=a_0s''_1s''_2\do s''_{r'}\uD ba_0\i.\tag a$$
\endproclaim
Let $A'=\part_!(\fD(A))$ with $\part$ as in 28.19. Then 
$A'\dsv\bK^{\ss',\cl'}_{\e_D(J),D'}$ where $\ss'=(s'_1,s'_2,\do,s'_{r'})$ is given by 
$s'_k=s''_{r'+1-k}$ and $D'=D\i$, $\cl'=(\uD')^*(\cl'')$, see 28.17, 28.19; hence
$A*A'\in\cd(Z_{J,\D})$ is well defined with $\D=D'D=G^0$. By 32.23(b) we have
$$H^0_c(\text{point},p_!\io^*(A*A'))=H^0_c(Z_{J,D},A\ot\part^*(A'))=
H^0_c(Z_{J,D},A\ot\fD(A)).$$
The last vector space is one dimensional, see \cite{\CS, II,7.4}. It follows that  

$H^0_c(\text{point},p_!\io^*(A*A'))\ne 0$.
\nl
Now some shift of $A*A'$ is a direct summand of 
$\bK^{\ss,\cl}_{J,D}*\bK^{\ss',\cl'}_{\e_D(J),D'}=\r_!(\bcl\bxt\bcl')$ (see 32.6). 
Hence $H^z_c(\text{point},p_!\io^*(\r_!(\bcl\bxt\bcl'))\ne 0$ for some $z\in\ZZ$. Using
this and 32.23(c) we see that there exists $\aa\in\SS$ such that 
$H^z_c(\text{point},p_!\io^*\r'_{\aa!}(\bcl\bxt\bcl'))\ne 0$ for some $z\in\ZZ$. Using 
this and 32.23(d) we see that there exists $e\in\ZZ$ such that
$H^z_c(\text{point},K_e)\ne 0$ for some $z,e\in\ZZ$. Using this and 32.23(e) we see 
that there exists $e'\in\ZZ$ such that $H^z_c(\text{point},K'_{e'}{}^e)\ne 0$ for some
$z,e'\in\ZZ$. As in 32.23 we see that we must have $\aa\in\SS_0$. Thus, there exists a 
sequence $a_0,a_1,\do,a_{r+r'}$ in $\WW^{r+r'}$ such that 

$a_k\in\{a_{k-1},s_ka_{k-1}\}$ for $k\in[1,r]$,    

$a_k\in\{a_{k-1},a_{k-1}(s''_{k-r})$ for $k\in[r+1,r+r']$,

$a_0\in\WW_J$, $a_0(\l'')=\l$, $a_{r+r'}=\e(a_0)$,

$i\in[1,r],a_{i-1}=a_i\imp s_1s_2\do s_i\do s_2s_1\in\WW_\l$

$j\in[1,r'],a_{j+r-1}=a_{j+r}\imp s''_1\do s''_{j-1}s''_js''_{j-1}\do s''_1
\in\WW_{\l''}$.
\nl
For $i\in[1,r]$ we set $t_i=s_1s_2\do s_i\do s_2s_1$ if $a_{i-1}=a_i$ and
$t_i=1$ if $a_{i-1}\n a_i$. Then $t_i\in\WW_\l$ and

$s_1s_2\do s_i=t_is_1s_2\do s_{i-1}a_{i-1}a_i\i$.
\nl
It follows that $s_1s_2\do s_r=t_rt_{r-1}\do t_1a_0a_r\i$. Similarly for $j\in[1,r']$ 
we set $t''_j=s''_1s''_2\do s''_j\do s''_2s''_1$ if $a_{j+r-1}=a_{j+r}$ and 
$t''_j=1$ if $a_{j+r-1}\ne a_{j+r}$. Then $t''_j\in\WW_{\l''}$ and 
$s''_1s''_2\do s''_j=t''_js''_1s''_2\do s''_{j-1}a_{j+r-1}\i a_{j+r}$. It follows that 
$s''_1s''_2\do s''_{r'}=t''_{r'}\do t''_2t''_1a_r\i a_{r+r'}$. Setting 

$\t=t_rt_{r-1}\do t_1,\t''=t''_{r'}\do t''_2t''_1$,
\nl
we have $\t\in\WW_\l$, $\t''\in\WW_{\l''}$, 
$s_1s_2\do s_r=\t a_0a_r\i, s''_1s''_2\do s''_{r'}=\t''a_r\i\e(a_0)$. Let 
$b'=(a_0\i\t a_0)\t''{}\i$. Then $b'\in\WW_{\l''}$ and
$s_1s_2\do s_r\uD=a_0b's''_1s''_2\do s''_{r'}\uD a_0\i$. We set
$b=(s''_1s''_2\do s''_{r'}\uD)\i b's''_1s''_2\do s''_{r'}\uD$. Since 
$s''_1s''_2\do s''_{r'}\uD\l''=\l''$ we have $b\in\WW_{\l''}$. Moreover,
$s_1s_2\do s_r\uD=a_0s''_1s''_2\do s''_{r'}\uD ba_0\i$. The lemma is proved.

\subhead 32.25\endsubhead  
Given $(w,\l),(w',\l')$ in $\Wb\T\ufs(\TT)$ we say that $(w,\l)\asi_J(w',\l')$ if there
exist $a\in\WW_J,b\in\WW_{\l'}$ such that $w=aw'ba\i,\l=a(\l')$. We then have 
$w'=a\i w(ab\i a\i)a$ where $a\i\in\WW_J$, $a\i(\l)=\l'$, 
$ab\i a\i\in\WW_{a\l'}=\WW_\l$ hence $(w',\l')\asi_J(w,\l)$. If, in addition, we 
have $(w',\l')\asi_J(w'',\l'')$ that is, $w'=\ta w''\tb\ta\i,\l'=\ta(\l'')$ with 
$\ta\in\WW_J,\l'=\ta(\l''),\tb\in\WW_{\l''}$, then $w=a\ta w''(\tb\ta\i b\ta)\ta\i a\i$
where $a\ta\in\WW_J$, $\l=a\ta(\l'')$,
$\tb\ta\i b\ta\in\WW_{\l''}\WW_{\ta\i(\l')}=\WW_{\l''}$, hence 
$(w,\l)\asi_J(w'',\l'')$. We see that $\asi_J$ is an equivalence relation on 
$\Wb\T\ufs(\TT)$.

We can now reformulate 32.24 as follows.

(a) {\it To $A\in\hZ_{J,D}$ we can associate an equivalence class $\fE_A$ under 
$\asi_J$ so that the following holds. If $\cl\in\fs(\TT)$, $\l$ is the isomorphism 
class of $\cl$ and $\ss=(s_1,s_2,\do,s_r)$ is a sequence in $\II$ such that
$s_1s_2\do s_r\uD\l=\l$ and $A\dsv\bK^{\ss,\cl}_{J,D}$ then 
$(s_1s_2\do s_r\uD,\l)\in\fE_A$.}

In particular:

(b) {\it To $A\in\hZ_{J,D}$ we can associate a $\WW_J$-orbit $\co_A$ on $\ufs(\TT)$ so 
that the following holds. If $\cl\in\fs(\TT)$, $\l$ is the isomorphism class of $\cl$ 
and $\ss=(s_1,s_2,\do,s_r)$ is a sequence in $\II$ such that 
$s_1s_2\do s_r\uD\l=\l$ and $A\dsv\bK^{\ss,\cl}_{J,D}$ then $\l\in\co_A$.}

\subhead 32.26\endsubhead  
Assume now that $J=\II$. We write $\asi$ instead of $\asi_{\II}$. Thus,
$(w,\l)\asi(w',\l')$ if there exist $a\in\WW,b\in\WW_{\l'}$ such that 
$w=aw'ba\i,\l=a(\l')$. Let $n\in\NN^*_\kk$.

Let $A$ be a character sheaf on $D$. Let $\fE_A$ be the equivalence class in
$\Wb\T\ufs(\TT)$ under $\asi$ defined by $A$ (see 32.25(a)). Let $\z^A:H_n[D]@>>>\ca$ 
be as in 31.7. We show:

(a) {\it If $\ss=(s_1,s_2,\do,s_r)$ is a sequence in $\II$, $\l\in\ufs_n$ and 
$\z^A(C^\ss_{\uD\l}[D])\ne 0$ then $(s_1s_2\do s_r\uD,\l)\in\fE_A$.}
\nl
Indeed, choose $\cl\in\fs(\TT)$ in the isomorphism class $\l$. Our assumption implies 
that $s_1s_2\do s_r\uD\l=\l$ hence $\bK^{\ss,\cl}_D$ is defined. Moreover our 
assumption implies $\sum_j(-v)^jv^{-\dim G}(A:{}^pH^j(\bK^{\ss,\cl}_D))\ne 0$. In 
particular, $A\dsv\bK^{\ss,\cl}_D$. Hence (a) follows from 32.25(a).

We show:

(b) {\it Let $(x,\l)\in\WW\T\ufs_n$ be such that $\z^A(\tT_x1_{\uD\l}[D])\ne 0$. Then
$(x\uD,\l)\in\fE_A$.}
\nl
We argue by induction on $l(x)$. If $x=1$ we have $\tT_x1_{\uD\l}=C^\ss_{\uD\l}$ where 
$\ss$ is the empty sequence and the result follows from (a). Assume now that
$l(x)\ge 1$. From our assumption we have $x\uD\l=\l$. We can find a sequence 
$\ss=(s_1,s_2,\do,s_r)$ in $\II$ with $x=s_1s_2\do s_r$, $r=l(x)$. From the definitions
we have $C^\ss_{\uD\l}=\sum_{y\in\WW_{\uD\l},xy\le x}c_y\tT_{xy}1_{\uD\l}$ with 
$c_y\in\ca,c_1=v^r$. Hence

$\z^A(C^\ss_{\uD\l}[D])=\sum_{y\in\WW_{\uD\l},xy\le x}c_y\z^A(\tT_{xy}1_{\uD\l}[D])$.
\nl
If $\z^A(\tT_{xy}1_{\uD\l}[D])\ne 0$ for some $y\in\WW_{\uD\l},xy<x$ then, by the 
induction hypothesis, we have $(xy\uD,\l)\in\fE_A$; we have $(xy\uD,\l)\asi(x\uD,\l)$ 
so that $(x\uD,\l)\in\fE_A$, as required.

We may therefore assume that $\z^A(\tT_{xy}1_{\uD\l}[D])=0$ for all $y\in\WW_{\uD\l}$ 
such that $xy<x$. Then we have $\z^A(C^\ss_{\uD\l}[D])=v^r\z^A(\tT_x1_{\uD\l}[D])$. 
Hence $\z^A(C^\ss_{\uD\l}[D])\ne 0$. Using (a) we see that $(x\uD,\l)\in\fE_A$, as 
required. This proves (b).

\head 33. Disjointness\endhead
\subhead 33.1\endsubhead
We fix an irreducible component $D$ of $G$. For $(L,S)\in\AA$ with $S\sub D$ and  
$\ce\in\cs(S)$ we define $\fK$ as in 5.6; we regard $\fK$ as a complex on $D$, zero 
outside $\bY_{L,S}$ and we write $(L,S,\ce)\btr_G\fK$.

\proclaim{Lemma 33.2}Let $(L,S)\in\AA,(L',S')\in\AA$ with $S\sub D,S'\sub D$. Let 
$\ce\in\cs(S),\ce'\in\cs(S')$. Let $(L,S,\ce)\btr_G\fK$, $(L',S',\ce')\btr_G\fK'$. 
Assume that $\ce$ (resp. $\ce'$) is strongly cuspidal and clean (see 23.3) relative to 
$N_GL$ (resp. $N_GL'$), that $L=G^0$ hence $Y_{L,S}=S$ and that $Y_{L',S'}\ne S$. Then 
for any $i$, the local systems $\cce$, $\ch^i\fK'|_S$ have no common irreducible direct
summand. 
\endproclaim
If $L'=G^0$ then, since $\ce'$ is clean, we have $\ch^i\fK'|_S=0$. Assume now that 
$L'\ne G^0$. By 23.7 we have $H^j_c(G,\fK\ot\fK')=0$ for all $j$. Since 
$\fK=IC(\bS,\ce)$ and $\ce$ is clean, we have $H^j_c(S,\fK\ot\fK')=0$ for all $j$ hence
$H^i_c(S,\ce\ot\fK')=0$ for all $j$. We must show that the local system 
$\ch^i(\ce\ot(\fK'|_S))$ on $S$ has no direct summand $\bbq$. Assume that 
$\ch^{i_0}(\ce\ot(\fK'|_S))$ has a direct summand $\bbq$ and that $i_0$ is maximum
possible with this property. If $a=\dim S'$, we have 
$H^{2a}_c(S,\ch^{i_0}(\ce\ot(\fK'|_S)))\ne 0$. Hence $E^{2a,i_0}_2\ne 0$ in the 
standard spectral sequence 
$$E^{p,q}_2=H^p_c(S,\ch^q(\ce\ot(\fK'|_S)))\imp H^{p+q}_c(S,\ce\ot(\fK'|_S)).$$
By the proof of 23.5 we have $H^p_c(S,\ce_1)=0$ for any $\ce_1\in\cs(S)$ which has no 
direct summand $\bbq$; in particular, taking $\ce_1=\ch^i(\ce\ot(\fK'|_S))$ with 
$i>i_0$ we see that $E^{p,q}_2=0$ if $q>i_0$. Clearly, $E^{p,q}_2=0$ if $p>2a$, hence 
$E^{2a,i_0}_2=E^{2a,i_0}_3=\do=E^{2a,i_0}_\iy$. Since $E^{2a,i_0}_2\ne 0$, it follows 
that $H^{2a+i_0}_c(S,\ce\ot(\fK'|_S))\ne 0$, a contradiction. The lemma is proved.

\proclaim{Proposition 33.3}Let $(L,S)\in\AA,(L',S')\in\AA$ with $S\sub D,S'\sub D$. 
Let $\ce\in\cs(S),\ce'\in\cs(S')$. Let $(L,S,\ce)\btr_G\fK$, $(L',S',\ce')\btr_G\fK'$. 
Assume that $\ce$ and $\cce$ (resp. $\ce'$) are strongly cuspidal and clean relative to
$N_GL$ (resp. $N_GL'$). Let $A$ (resp. $A'$) be an admissible complex on $D$ (see 6.7) 
which is a direct summand of $\fK$ (resp. of $\fK'$). Assume that $A\not\cong A'$. Let 
$Y=Y_{L,S}$. Let $\cf$ be the local system $A|_Y$. Then for any $i$, $\cf$ is not a 
direct summand of $\ch^i(A')|_Y$ (which is a local system by 25.2).
\endproclaim
Since $\bY_{L',S'}$ is a union of strata of $D$, we have either $Y\cap\bY_{L',S'}=0$ or
$Y=Y_{L',S'}$ or $Y\sub\bY_{L',S'}-Y_{L',S'}$. In the first case we have 
$\ch^i(A')|_Y=0$ and the result is obvious. In the second case we have $\ch^i(A')|_Y=0$
unless $i=0$ and since $A\not\cong A'$, the local system $\ch^0(A)|_Y$ is irreducible, 
non-isomorphic to $\cf$. Thus, we may assume that $Y\sub\bY_{L',S'}-Y_{L',S'}$. It is 
enough to show that for any $i$, $\cf$ is not a direct summand of $\ch^i(\fK')|_Y$ (a 
local system, by 25.2). Let $\d$ be the connected component of $N_GL$ such that 
$S\sub\d$. Let $su=us\in S^*$ with $s$ semisimple, $u$ unipotent. Let $\tid$ be the 
connected component of $Z_G(s)$ such that $u\in\tid$. Since $su$ is isolated in $N_GL$,
we have ${}^\d\cz_L^0={}^{\tid}\cz_{Z_L(s)^0}$; we denote this torus by $\ct$. Let 
$R_1$ be the subvariety of $S$ consisting of all elements of the form $yzsuy\i$ with 
$y\in Z_L(s)^0,z\in\ct$. Since $R_1$ is an orbit of a connected group, it is smooth, 
irreducible. Let $R_1^*=R_1\cap S^*$, an open dense subset of $R_1$ (see 25.4, 25.6).
Let $\p_1:\p\i(R_1^*)@>>>R_1^*$ be the restriction of $\p:\tY_{L,S}@>>>Y$ (as in 3.13).
Let $\tce$ be the local system on $\tY_{L,S}$ defined in 5.6; its restriction to
$\p\i(R_1^*)$ is denoted again by $\tce$. From the definitions, we have
$\fK|_{R_1^*}=\p_{1!}\tce$. By the proof of 3.13(a) we have

$\p\i(R_1^*)=\sqc_{xL\in N(L,S)/L}\{(g,xL);g\in R_1^*\}$
\nl
where $N(L,S)=\{x\in N_{G^0}L,xSx\i=S\}$. Define $\e:s\i R_1@>>>R_1$ by $g\m sg$. We 
see that $\e^*\fK|_{R_1^*}=\op_{xL\in N(L,S)/L}\ce^x|_{s\i R_1^*}$ where $\ce^x$ is the
local system on $s\i R_1$ obtained by taking inverse image of $\ce$ under
$s\i R_1@>>>S,g\m xsgx\i$. Now $s\i R_1$ is an isolated stratum of $Z_G(s)$ contained 
in the connected component $\tid$ (it is the stratum containing $u$). From 23.4 we see 
that $\ce^x$ and $\cce^x$ are strongly cuspidal and clean with respect to $Z_G(s)$. By 
16.12 we can find complexes $\fK'_j(j\in[1,m])$ on $\tid$ of the same type as $\fK'$
and an open subset $\cu$ of $\tid$ containing all unipotents in $\tid$ such that  

(a) $\e'{}^*(\fK'|_{s\cu})\cong\op_j\fK'_j|_\cu$,
\nl
where $\e':\cu@>>>s\cu$ is $g\m sg$. Note that $R_1^*\cap s\cu$ contains $su$ hence is
non-empty. Since $s\cu$ is open in $s\tid$, and $R_1$ is an irreducible subset of
$s\tid$, we see that $R_1\cap s\cu$ is an open dense subset of $R_1$. Since $R_1^*$ is
another open dense subset of $R_1$ we see that 

(b) {\it $R_1^*\cap s\cu=(R_1\cap s\cu)\cap R_1^*$ is open dense in $R_1$.}
\nl
It suffices to show that the local systems $\ch^i(\fK')|_{R_1^*\cap s\cu}$, 
$\fK|_{R_1^*\cap s\cu}$  have no common irreducible direct summand. Using (a) we see 
that it suffices to show that for any $j\in[1,m],x\in N(L,S)$, 

(c) {\it the local systems $\ch^i(\fK'_j)|_{s\i R_1^*\cap\cu}$, 
$\ce^x_{s\i R_1^*\cap\cu}$ have no common irreducible direct summand.}
\nl
Since $s\i R_1$ is an isolated stratum of $Z_G(s)$, $\ch^i(\fK'_j)|_{s\i R_1}$ is a 
local system. Using (b) we see that (c) would follow from the following statement:

(d) {\it the local systems $\ch^i(\fK'_j)|_{s\i R_1}$, $\ce^x$ have no common 
irreducible direct \lb summand.}
\nl
By 16.12(b) we may assume that there exists $x'\in G^0$ such that $x'{}\i sx'\in S'_s$ 
and the following holds. Let $L''=x'L'x'{}\i,S''=x'S'x'{}\i$, $L'_0=L''\cap Z_G(s)^0$, 

$S'_0$ is a stratum of $N_GL'_0$ contained in $\tid$, containing unipotent elements and
such that $S'_0\sub s\i S''$, 

$\ce'_0$ is the local system on $S'_0$, inverse image of $\ce'$ under $S'_0@>>>S'$,
$g\m x'{}\i sgx'$,
\nl
$(L'_0,S'_0,\ce'_0)\btr_{Z_G(s)}\fK'_j$.
\nl
From 23.4 we see that $\ce'_0$ is strongly cuspidal and clean with respect to 
$N_{Z_G(s)}(L'_0)$. We see that (d) follows from 33.2 (applied to $Z_G(s),\fK'_j,
\cce^x$ instead of $G,\fK',\ce$) provided we can show that 

(e) {\it $Y_{L'_0,S'_0}$ (defined in terms of $Z_G(s)$) is not equal to $s\i R_1$.}
\nl
Assume that $Y_{L'_0,S'_0}=s\i R_1$. Since $s\i R_1$ is an isolated stratum of 
$Z_G(s)$, it follows that $L'_0=Z_G(s)^0$ and $S'_0=s\i R_1$ hence $Z_G(s)^0\sub L''$ 
and $u\in S'_0$. Since $sS'_0\sub S''$, we have $su\in S''$. We can find a parabolic 
$P'$ of $G^0$ with Levi $L'$ such that $S'\sub N_GP'\cap N_GL'$ hence
$S''\sub N_G(x'P'x'{}\i)\cap N_G(x'L'x'{}\i)$. We see that
$su\in N_G(x'P'x'{}\i)\cap N_G(x'L'x'{}\i)$. Using 2.1(c) with $g=su,Q=x'P'x'{}\i$, we
see that $L(su)\sub x'L'x'{}\i=L''$ where $L(su)$ is defined as in 2.1. We can find a 
parabolic $P$ of $G^0$ with Levi $L$ such that $S\sub N_GP\cap N_GL$ hence 
$su\in N_GP\cap N_GL$. Moreover, $su$ is isolated in $N_GP\cap N_GL$. From 3.8(a) we 
see that $L\sub L(su)$. Combining with $L(su)\sub L''$, we see that $L\sub L''$. Since 
$Y\sub\bY_{L',S'}-Y_{L',S'}$, we have $\bY\sub\bY_{L',S'}$. Taking images under the map
$\s:D@>>>D//G^0$ (see 7.1) we obtain $\dim\s(\bY)\le\dim\s(\bY_{L',S'})$. Using 7.3(b) 
we can rewrite the last inequality in the form 
$\dim({}^\d\cz_L^0)\le\dim({}^{\d'}\cz_{L'}^0)$ where $\d'$ is the connected component 
of $N_GL'$ that contains $S'$. Equivalently,

(f) $\dim({}^\d\cz_L^0)\le\dim({}^{\d''}\cz_{L''}^0)$
\nl
where $\d''$ is the connected component of $N_GL''$ that contains $S''$. From 
$L\sub L''$ we deduce $\cz_{L''}\sub\cz_L$. Intersecting both sides with $Z_G(su)$ and 
noting that $su\in\d,su\in\d''$ we see that ${}^{\d''}\cz_{L''}\sub{}^\d\cz_L$. Taking 
identity components we have ${}^{\d''}\cz_{L''}^0\sub{}^\d\cz_L^0$. Using (f) we deduce
${}^{\d''}\cz_{L''}^0={}^\d\cz_L^0$. Taking the centralizer of both sides in $G^0$ and 
using 1.10(a) we obtain $L=L''$. Now $S$ and $S''$ are strata of $N_GL=N_GL''$ which 
contain a common point, $su$. Hence $S=S''$. Since $(L,S)=(L'',S'')$ we have 
$Y_{L,S}=Y_{L'',S''}$ hence $Y=Y_{L',S'}$. This contradicts 
$Y\sub\bY_{L',S'}-Y_{L',S'}$ and proves (e). The proposition is proved.

\subhead 33.4\endsubhead
Let $\ci$ be a finite collection of mutually non-isomorphic character sheaves on $D$
and let $A\in\ci$. Let $Y=Y_{L,S}$ be the stratum of $D$ such that $\supp(A)=\bY$. Let
$\tY=\{(g,xL)\in D\T G^0/L;x\i gx\in S^*\}$ (see 3.13). Define $\p_1:\tY@>>>Y$ by
$\p_1(g,xL)=g$. By 25.2, for any $A'\in\ci$ and any $i\in\ZZ$ there exists a local 
system $\ce\in\cs(S)$ such that $\ch^i(A')|_Y$ is a local system isomorphic to a direct
summand of $\p_{1!}\tce$ with $\tce$ as in 5.6. Replacing $\ce$ by the direct sum of
the local systems $\ce$ (for various $j,i$ as above) we see that we may assume that 
$\ce$ is the same for all $A',i$. We can find $n'\in\NN^*_\kk$ such that 
$\ce\in\cs_{n'}(S)$. Let $\d$ be the connected component of $N_GL$ that contains $S$. 
Let $g_1\in S$. Let $H=\{(z_1,l_1)\in{}^\d\cz^0_L\T L;l_1z_1^{n'}g_1l_1\i=g_1\}$. Let 

$V=\{(g,x,z,l)\in D\T G^0\T{}^\d\cz^0_L\T L;x\i gx=lz^{n'}g_1l\i\in S^*\}$.
\nl
Now $V$ is irreducible; it is isomorphic to the product of $G^0$ with an open dense 
subset of ${}^\d\cz^0_L\T L$. We have a commutative diagram with cartesian squares
$$\CD
{}@.\tY'@<a'<<Z'@>b'>>S'\\
@.@V\p_2VV @V\p_3VV @V\p_4 VV\\
Y@<\p_1<<\tY@<a<<Z@>b>>S
\endCD$$
where 

$S'$ is the space of $H^0$-orbits on ${}^\d\cz^0_L\T L$ for the free $H^0$-action by 
right translation,

$Z=\{(g,x)\in D\T G^0;x\i gx\in S^*\}$,

$\tY'$ is the space of $(L\T H^0)$-orbits on $V$ for the free $L\T H^0$-action 
$(l_0,(z_1,l_1)):(g,x,z,l)\m(g,xl_0\i,zz_1\i,l_0ll_1\i)$,

$Z'$ is the space of $H^0$-orbits on $V$ for the free $H^0$-action 
$(z_1,l_1)):(g,x,z,l)\m(g,x,zz_1\i,ll_1\i)$,

$a(g,x)=(g,xL),b(g,x)=x\i gx$, $a'$ is the obvious map, $b'(g,x,z,l)\m(z,l)$,
$\p_2(g,x,z,l)=(g,xL),\p_3(g,x,z,l)=(g,x)$, $\p_4(z_1,l_1)=l_1\i z_1^{n'}g_1l_1$.
\nl
Now $\tY'$ is irreducible since $V$ is irreducible; $\tY$ is irreducible since it
equals $\p_2(\tY')$. Since $\ce\in\cs_{n'}(S)$, the local system $\p_4^*\ce$ on $S'$ is
$({}^\d\cz^0_L\T L)$-equivariant (for the action by left translation). Since this 
action is transitive with connected isotropy groups, we see that $\p_4^*\ce\cong\bbq^e$
for some integer $e\ge 1$. Hence $\p_3^*b^*\ce=b'{}^*\p_4^*\ce\cong\bbq^e$. By 
definition, $a^*\tce=b^*\ce$. Hence 
$a'{}^*\p_2^*\tce=\p_3^*a^*\tce=\p_3^*b^*\ce\cong\bbq^e$. Since $a'$ is a principal
$L$-bundle, it follows that $\p_2^*\tce\cong\bbq^e$. Now $\p_0:=\p_1\p_2:\tY'@>>>Y$ is 
a composition of two (finite) principal coverings ($\p_1$ is a principal 
$H/H^0$-covering since $\p_4$ is a principal $H/H^0$-covering; $\p_1$ is a principal 
covering by 3.13(a)) hence it is a not necessarily principal, finite unramified 
covering. Let $N=|\p_0\i(y)|$ for some/any $y\in Y$. Let $Y''$ be the set of all pairs 
$(y,f)$ where $y\in Y$ and $f:\{1,2,\do,N\}@>>>\p_0\i(y)$ is a bijection. Then $Y''$ is
an algebraic variety and $\p'_0:Y''@>>>Y,(y,f)\m y$ is a principal covering whose group
is the symmetric group $\fS_N$. Moreover, $\p'_0$ factors as $Y''@>\t>>\tY'@>\p_0>>Y$ 
where $\t(y,f)=f(1)$. Let $\hY$ be a connected component of $Y''$. Then 
$\t_0:\hY@>>>\tY'$ (restriction of $\t$) is a finite unramified covering. Let 
$\p:\hY@>>>Y$ be the restriction of $\p'_0$. Then $\p$ is a (finite) principal bundle 
whose group is the group $\G$ consisting of all elements of $\fS_N$ which map $\hY$ 
into itself. Moreover, $\p$ factors as $\hY@>\t_1>>\tY@>\p_1>>Y$ where $\t_1=\p_1\t_0$ 
is a finite unramified covering. Since $\p_2^*\tce\cong\bbq^e$, we have 
$\t_1^*\tce\cong\bbq^e$. Hence any irreducible direct summand of the local system 
$\tce$ is a direct summand of $\t_{1!}\bbq$. Now let $\ce_1$ be an irreducible 
local system on $Y$ which is a direct summand of $\p_{1!}\tce$. We can find an 
irreducible direct summand $\ce_2$ of $\tce$ such that $\ce_1$ is a direct summand of 
$\p_{1!}\ce_2$. Then $\ce_2$ is a direct summand of $\t_{1!}\bbq$, hence $\p_{1!}\ce_2$
is a direct summand of $\p_{1!}\t_{1!}\bbq=\p_!\bbq$. Since $\ce_1$ is a direct summand
of $\p_{1!}\ce_2$ it follows that 

(a) {\it $\ce_1$ is a direct summand of $\p_!\bbq$.}
\nl
Let $\cc$ be the category whose objects are local systems on $Y$ which are direct sums
of irreducible direct summands of $\p_!\bbq$. Let $\cc_\G$ be the category of
$\bbq[\G]$-modules of finite dimension over $\bbq$. We have an equivalence of 
categories $\cc_\G@>>>\cc$: it attaches to an object $M$ of $\cc_\G$ the local system 
$[M]=(M^*\ot\p_!\bbq)^\G$ in $\cc$; here $\p_!\bbq$ is regarded naturally as a local 
system with $\G$-action, $M^*$ is the dual of $M$  and the superscript denotes 
$\G$-invariants. Using (a) and the definition of $\ce$, we see that, for any 
$A'\in\ci,i\in\ZZ$, we have $\ch^i(A')|_Y\in\cc$. Hence $\ch^i(A')|_Y\cong[M_{A',i}]$ 
for some $M_{A',i}\in\cc_\G$, well defined up to isomorphism. Let $e=\dim Y$. Then 
$M_{A,-e}$ is an irreducible object of $\cc_\G$. 

In the remainder of this section we assume that:

(b) {\it $D$ is clean in the sense that, for any parabolic subgroup $P$ of $G^0$ such 
that $N_DP\ne\em$, any cuspidal character sheaf of $N_DP/U_P$ is $0$ on the complement 
of some isolated stratum of $N_DP/U_P$.}
\nl
We show:

(c) {\it if $A'\in\ci,i\in\ZZ$ and $A'\ne A$ then $M_{A',i}$ contains no direct summand
isomorphic to $M_{A,-e}$.}
\nl
This follows from 33.3 which is applicable in view of (b), the admissibility of 
character sheaves (30.6), the strong cuspidality of cuspidal character sheaves (31.15)
and the fact that $\fD(A)$ is a character sheaf (28.18).

In the remainder of this section we assume that $\kk$ is an algebraic closure of a 
finite field $\FF_q$ and that $G$ has a fixed $\FF_q$-rational structure whose 
Frobenius map $F$ induces the identity map on $G/G^0$. Replacing $\FF_q$ by a finite
extension, we may assume that $F(Y)=Y$, that $\hY$ and $\p:\hY@>>>Y$ are defined over
$\FF_q$, that the Frobenius map $F:\hY@>>>\hY$ satisfies $F(\g\hy)=\g F(\hy)$ for all 
$\g\in\G,\hy\in\hY$, that $F^*A'\cong A'$ for all $A'\in\ci$ and that for any $\g\in\G$
and any integer $m\ge 1$ there exists $\hy_{\g,m}\in\hY$ such that 
$F^m(\hy_{\g,m})=\g\hy_{\g,m}$. (We then set $y_{\g,m}=\p(\hy_{\g,m})$.)

Let $M\in\cc_\G$. The stalk of $[M]$ at $y\in Y$ is the vector space 
$$[M]_y=\{f:\p\i(y)@>>>M^*;f(\hy)=\g(f(\g\i\hy))\text{ for all }\g\in\G,\hy\in\hY\}.$$
Let $m\ge 1$. For any $R\in\Aut_{\cc_\G}(M^*)$ there is a unique isomorphism of local 
systems $\tR:F^{m*}[M]@>\si>>[M]$ such that for any $y\in Y$, $\tR$ induces on stalks
the linear map $\tR_y:[M]_{F^m(y)}@>>>[M]_y$ which to a function 
$f:\p\i(F^m(y))@>>>M^*$ associates the function $f':\p\i(y)@>>>M^*$ given by 
$f'(\hy)=R(f(F^{-m}(\hy)))$. Clearly, any isomorphism $F^{m*}[M]@>\si>>[M]$ is of the 
form $\tR$ for a unique $R$ as above.

For $\g\in\G$ we have an isomorphism 

(d) $[M]_{y_{\g,m}}@>\si>>M^*,f\m f(\hy_{\g,m})$.
\nl
If $R$ is as above then $\tR_{y_{\g,m}}$ maps $[M]_{y_{\g,m}}$ into itself (since 
$F^m(y_{\g,m})=y_{\g,m}$) and it corresponds under (d) to the automorphism 
$\g\i R=R\g\i:M^*@>>>M^*$. Hence

(e) $\tr(\tR_{y_{\g,m}},[M]_{y_{\g,m}})=\tr(\g\i R,M^*)=\tr({}^tR\g,M)$.

\subhead 33.5\endsubhead
Let $V$ be an algebraic variety defined over $\FF_q$ with Frobenius map $F:V@>>>V$. Let
$K\in\cd(V)$ and let $\ph:F^*K@>\si>>K$ be an isomorphism. For any integer $m\ge 1$ we 
denote by $\ph^{(m)}:F^{m*}K@>\si>>K$ the composition
$$(F^m)^*K@>(F^{m-1})^*\ph>>(F^{m-1})^*K@>(F^{m-2})^*>>\do@>F^*\ph>>F^*K@>\ph>>K.$$

\subhead 33.6\endsubhead
For each $A'\in\ci$ we choose an isomorphism $\k_{A'}:F^*A'@>\si>>A'$. Let 
$\k'_A:F^*\fD(A)@>\si>>\fD(A)$ be the isomorphism such that for any $y\in Y$ the 
isomorphism $\ch^{-e}\fD(A)_{F(y)}@>\si>>\ch^{-e}\fD(A)_y$ (that is, 
$\ch^{-e}(A)\che{}_{F(y)}@>\si>>\ch^{-e}(A)\che{}_y$) induced by $\k'_A$ is 
$q^{\dim D-e}$ times the transpose inverse of the isomorphism 
$\ch^{-e}(A)_{F(y)}@>\si>>\ch^{-e}(A)_y$ induced by $\k_A$.

\proclaim{Proposition 33.7} Let $A'\in\ci$. For any integer $m\ge 1$ we have
$$q^{-(\dim D-e)m}|\G|\i\sum_{\g\in\G}\c_{A',\k_{A'}^{(m)}}(y_{\g,m})
\c_{\fD(A),\k'_A{}^{(m)}}(y_{\g,m})=\d_{A,A'}.\tag a$$
\endproclaim
Under an isomorphism $\ch^i(A')|_Y\cong[M_{A',i}]$, the isomorphism \lb
$F^{*m}\ch^i(A')@>\si>>\ch^i(A')$ induced by $\k_{A'}^{(m)}:F^{*m}A'@>\si>>A'$ 
corresponds to an isomorphism $F^{*m}[M_{A',i}]@>\si>>[M_{A',i}]$ which must be of the 
form $\tR$ for some $R=R_{m,A',i}\in\Aut_{\cc_\G}(M_{A',i}^*)$ hence
$$\tr(\k_{A'}^{(m)},\ch^i(A')_{y_{\g,m}})=\tr({}^tR_{m,A',i}\g,M_{A',i}).$$
Next we have
$$\align&\tr(\k'_A{}^{(m)},\ch^{-e}(\fD(A))_{y_{\g,m}})
=q^{(\dim D-e)m}\tr((\k_A^{(m)})\i,\ch^{-e}(A)_{y_{\g,m}})\\&
=q^{(\dim D-e)m}\tr({}^tR_{m,A,-e}\i\g\i,M_{A,-e}).\endalign$$
Hence the left hand side of (a) equals
$$\sum_i(-1)^{i+e}|\G|\i\sum_{\g\in\G}\tr({}^tR_{m,A',i}\g,M_{A',i})
\tr({}^tR_{m,A,-e}\i\g\i,M_{A,-e})\tag b$$
that is,
$$\sum_i(-1)^{i+e}|\G|\i\tr(({}^tR_{m,A',i}\ot {}^tR_{m,A,-e})
\sum_{\g\in\G}(\g\ot\g\i),M_{A',i}\ot M_{A,-e}).$$
Assume first that $A'\ne A$. To show that (b) is zero it is enough to show
that for any $i$, $\sum_{\g\in\G}(\g\ot\g\i)$ acts as $0$ on $M_{A',i}\ot M_{A,-e}$.
This follows from the fact that the $\G$-invariant part of the $\G$-module
$M_{A',i}\ot M_{A,-e}^*$ is zero (see 33.4(c)).

Assume next that $A'=A$. Then we have $M_{A',i}=0$ unless $i=-e$. We must show that
$$|\G|\i\sum_{\g\in\G}\tr({}^tR_{m,A,-e}\g,M_{A,-e})\tr({}^tR_{m,A,-e}\i\g\i,M_{A,-e})
=1.$$
Since $M_{A,-e}$ is an irreducible $\G$-module, ${}^tR_{m,A,-e}$ acts as on it as a 
scalar, hence the desired equality follows from the Schur orthogonality relations for 
irreducible characters of $\G$.

\head 34. The structure of $H^D_n$\endhead
\subhead 34.1\endsubhead
We give (a variant of) some definitions in \cite{\GE, 1}. Let $\car$ be a commutative 
ring with $1$. Let $\fA$ be an associative $\car$-algebra with $1$ with a given finite 
basis $B$ as an $\car$-module. We assume that $1$ is compatible with $B$ in the 
following sense: $1=\sum_\l1_\l$ where $1_\l\in B$ are distinct, 
$1_\l1_{\l'}=\d_{\l,\l'}1_\l$ and any $b\in B$ satisfies $1_\l b1_{\l'}=b$ for some 
(uniquely determined) $\l,\l'$. For $b,b'\in B$ we have 
$bb'=\sum_{b''\in B}r_{b,b'}^{b''}b''$ where $r_{b,b'}^{b''}\in\car$. We say that 
$b'\preceq b$ if $b'\in\cap_{K\in\cf;b\in K}K$ where $\cf$ is the collection of all
subsets $K$ of $B$ such that $\sum_{b_1\in K}\car b_1$ is a two-sided ideal of $\fA$; 
we say that $b\si b'$ if $b'\preceq b$ and $b\preceq b'$. This is an equivalence 
relation on $B$ and the equivalence classes are the {\it two-sided cells}. (Replacing
two-sided ideals by left ideals in the definition of $\preceq$ and of two-sided cells 
we obtain the notion of left cells. The left cells form a partition of $B$ finer than 
that given by two-sided cells.) We say that $b'\prec b$ if $b'\preceq b$ and 
$b'\not\si b$. For any $b\in B$ let $\fA_{\prec b}=\op_{b'\in B;b'\prec b}\car b$.

Assume now that $\car=\ca=\ZZ[v,v\i]$. Let $b\in B$. We can find an integer $m\ge 0$ 
such that $v^{-m}r_{b,b'}^{b''}\in\ZZ[v\i]$ for any $b',b''$ in the two-sided cell of 
$b$. The smallest such $m$ is denoted by $a(b)$. We say that $B$ satisfies $P_1$ if 
$a(b)=a(b')$ whenever $b,b'$ are in the same two-sided cell. Assume that $B$ satisfies 
$P_1$. For $b\in B$ we set $\hb=v^{-a(b)}b\in\fA$. Let 
$\fA^-=\sum_{b\in B}\ZZ[v\i]\hb\sub\fA$. Then $\fA^-$ is an associative 
$\ZZ[v\i]$-algebra for the multiplication 
$\hb*\hb'=\sum_{b''\in B;b''\si b}v^{-a(b)}r_{b,b'}^{b''}\hb''$ if $b\si b'$, 
$\hb*\hb'=0$ if $b\not\si b'$. Let $\fA^\iy=\fA^-/v\i\fA^-$ and let 
$t_b=\hb+v\i\fA^-\in\fA^\iy$. Then $\fA^\iy$ is a ring with $\ZZ$-basis 
$\{t_b;b\in B\}$ and with multiplication defined by 
$t_bt_{b'}=\sum_{b''\in B}\g_{b,b'}^{b''}t_{b''}$ where $\g_{b,b'}^{b''}\in\ZZ$ is 
given by $v^{-a(b)}r_{b,b'}^{b''}=\g_{b,b'}^{b''}\mod v\i\ZZ[v\i]$ if $b,b',b''$ are in
the same two-sided cell and $\g_{b,b'}^{b''}=0$, otherwise. We say that $B$ satisfies 
$P_2$ if $\fA^\iy$ has a unit element compatible with the basis $\{t_b;b\in B\}$. We 
say that $B$ satisfies $P_3$ if for any $b_1,b_2,b_3,b_4\in B$ such that $b_2\si b_4$ 
we have
$$\sum_{\b\in B;\b\si b_2}r^\b_{b_1,b_2}(v)r^{b_4}_{\b,b_3}(v')
=\sum_{\b\in B;\b\si b_2}r^{\b_4}_{b_1,\b}(v)r^\b_{\b_2,b_3}(v')$$
where $v'$ is an indeterminate independent of $v$. In this case, assuming also that
$b_2\si b_3\si b_4$, we pick the coefficient of $v'{}^{a(b_2)}=v'{}^{a(b_4)}$ in both 
sides and we obtain
$$\sum_{\b\in B;\b\si b_2}r^\b_{b_1,b_2}\g^{b_4}_{\b,b_3}
=\sum_{\b\in B;\b\si b_2}r^{\b_4}_{b_1,\b}\g^\b_{\b_2,b_3}.\tag a$$
Assume that $B$ satisfies $P_1,P_2,P_3$. The unit element of $\fA^\iy$ is of the form 
$\sum_{b\in\cd}t_b$ where $\cd\sub B$. We say that $\cd$ is the set of 
{\it distinguished} elements of $B$.

Let $\fA^\iy_\ca=\ca\ot\fA^\iy$. We define an 
$\ca$-linear map $\Ph:\fA@>>>\fA^\iy_\ca$ by 
$$\Ph(b)=\sum_{b_1\in\cd,b_2\in B;b_1\si b_2}r^{b_2}_{b,b_1}t_{b_2}$$
for $b\in B$. Then $\Ph$ is an $\ca$-algebra homomorphism taking $1$ to $1$. If we 
identify $\fA,\fA^\iy_\ca$ as $\ca$-modules via $b\lra t_b$, the obvious left 
$\fA^\iy_\ca$-module structure on $\fA^\iy_\ca$ becomes the left $\fA^\iy_\ca$-module 
structure on $\fA$ given by $t_b*b'=\sum_{b''\in B}\g_{b,b'}^{b''}b''$. For 
$x\in\fA,b\in B$ we have 

(b) $xb=\Ph(x)*b\mod\fA_{\prec b}$.
\nl
Indeed, we may assume that $x\in B$. Using (a), we have
$$\align&\Ph(x)*b=\sum_{b_1\in\cd,b_2\in B;b_1\si b_2}r^{b_2}_{x,b_1}t_{b_2}*b
=\sum_{b_1\in\cd,b_2,b''\in B;b_1\si b_2}r^{b_2}_{x,b_1}\g_{b_2,b}^{b''}b''\\&
=\sum_{b_1\in\cd,b_2,b''\in B;b_1\si b\si b''}r^{b_2}_{x,b_1}\g_{b_2,b}^{b''}b''
=\sum_{b_1\in\cd,b'_1,b''\in B;b_1\si b\si b''}r^{b''}_{x,b'_1}\g_{b_1,b}^{b'_1}b''\\&
=\sum_{b_1\in\cd,b''\in B;b_1\si b\si b''}r^{b''}_{x,b}\g_{b_1,b}^bb''
=\sum_{b''\in B;b\si b''}r^{b''}_{x,b}b''=xb\mod\fA_{\prec b},\endalign$$
as required.

Let $K$ be a field and let $\ca@>>>K$ be a homomorphism of rings with $1$. Let 
$\fA_K=K\ot_\ca\fA,\fA^\iy_K=K\ot_\ZZ\fA^\iy,\fA_{K,\prec b}=K\ot_\ca\fA_{\prec b}$
($b\in B$). Then $\Ph$ induces a $K$-algebra homomorphism $\Ph_K:\fA_K@>>>\fA^\iy_K$. 
We show:

(c) {\it If $\fA_K$ is a semisimple algebra then $\Ph_K$ is an (algebra) isomorphism.}
\nl
Since $\fA_K,\fA^\iy_K$ have the same (finite) dimension, it suffices to show that
$\Ph_K$ is injective. The $\fA^\iy_\ca$-module structure on $\fA$ extends to an 
$\fA^\iy_K$-module structure on $\fA_K$ denoted again by $*$. From (b) we deduce that
$xb=\Ph_K(x)*b\mod\fA_{K,\prec b}$ for any $x\in\fA_K,b\in B$. In particular, if
$x\in\Ker\Ph_K,b\in B$ then $xb\in\fA_{K,\prec b}$. Applying this repeatedly, we see 
that for any $m\ge 1$, any $x_1,x_2,\do,x_m$ in $\Ker\Ph_K$ and any $b\in B$,
$x_1x_2\do x_mb$ is a $K$-linear combination of elements $b'\in B$ such that
$b'=b_m\prec b_{m-1}\prec\do\prec b_0=b$ (with $b_i\in B$). If $m$ is large enough, no
such $b'$ exists. Thus for large enough $m$ we have $x_1x_2\do x_mb=0$ for all $b\in B$
hence $x_1x_2\do x_m=0$. We see that $\Ker\Ph_K$ is a nilpotent two-sided ideal of
$\fA_K$. Hence it is $0$. Thus $\Ph_K$ is injective and (c) is proved.

\subhead 34.2\endsubhead
Let $D$ be a connected component of $G^0$. Let $\WW^D$ be the subgroup of 
$\WW^\bul\sub\Aut(\TT)$ generated by $\WW$ and by $\uD$. Now $\WW$ is a normal subgroup
of $\WW^D$ and $\WW^D/\WW$ is a finite cyclic group.

We fix $n\in\NN^*_\kk$. Let $\l\in\ufs_n$. We write $R_\l$ instead of $R_\cl$ (see 
28.3) where $\l$ is the isomorphism class of $\cl\in\fs_n$. Then $R_\l$ is a root 
system and $R_\l^+=R_\l\cap R^+$ is a set of positive roots for $R_\l$. Let $\Pi_\l$ be
the unique set of simple roots for $R_\l$ such that $\Pi_\l\sub R_\l^+$. Recall that 
$\WW_\l$, the subgroup of $\WW$ generated by $\{s_\a;\a\in R_\l\}$ is the Weyl group of
the root system $R_\l$. Let $\II_\l=\{s_\a;\a\in\Pi_\l\}\sub\WW_\l$. Then 
$(\WW_\l,\II_\l)$ is a Coxeter group. 
Let $\WW^D_\l=\{w\in\WW^D;w\l=\l\}$. Let $\O^D_\l=\{w\in\WW^D_\l;w(R_\l^+)=R_\l^+\}$. 
(Here $\WW^D$ acts on $R$ by $w:\a\m w\a,(w\a)(t)=\a(w\i t)$ for $t\in\TT$.) Then 
$\WW_\l$ is a normal subgroup of $\WW^D_\l$, $\O^D_\l$ is a subgroup of $\WW^D_\l$ and 
$\WW^D_\l$ is the semidirect product of $\WW_\l$ and $\O^D_\l$. Define $l:\WW^D@>>>\NN$
by $l(w)=|\{\a\in R^+;w(\a)\in R^-\}|$. This extends the length function $\WW@>>>\NN$. 
Define $l_\l:\WW^D@>>>\NN$ by $l_\l(w)=|\{\a\in R_\l^+;w(\a)\in R^-\}|$. Then 
$\O^D_\l=\{w\in\WW^D_\l;l_\l(w)=0\}$, $\II_\l=\{w\in\WW_\l;l_\l(w)=1\}$. The standard 
partial order $\le_\l$ of the Coxeter group $\WW_\l$ is extended to a partial order 
$\le_\l$ on $\WW^D_\l$ as follows: if $w_1,w'_1\in\O^D_\l,w_2,w'_2\in\WW_\l$, we say 
that $w_1w_2\le_\l w'_1w'_2$ if $w_1=w'_1$ and $w_2\le_\l w'_2$. 

Let $H^D_\l$ be the $\ca$-algebra defined by the generators $\tT_w^\l(w\in\WW^D_\l)$ 
and relations

(a) $\tT_w^\l\tT_{w'}^\l=\tT_{ww'}^\l$ if $w,w'\in\WW^D_\l$, 
$l_\l(ww')=l_\l(w)+l_\l(w')$,

(b) $(\tT_\s^\l)^2=\tT_1^\l+(v-v\i)\tT_\s^\l$ for $\s\in\II_\l$.
\nl
Then $\{\tT_w^\l;w\in\WW^D_\l\}$ is an $\ca$-basis of $H^D_\l$. Let $H_\l$ be the 
$\ca$-submodule of $H^D_\l$ with $\ca$-basis $\{\tT_w^\l;w\in\WW_\l\}$. This is an 
$\ca$-subalgebra of $H^D_\l$. Let $\bar{}:H^D_\l@>>>H^D_\l$ be the unique ring 
homomorphism such that $\ov{v^m\tT_w^\l}=v^{-m}(\tT_{w\i}^\l)\i$ for all 
$w\in\WW^D_\l$, $m\in\ZZ$. From the definitions, for any $w\in\WW^D_\l$ we have 
$\ov{\tT_w^\l}-\tT_w^\l\in\sum_{y\in\WW^D_\l;y\le_\l w,y\ne w}\ca\tT_y^\l$. By an
argument similar to one in \cite{\HA, 5.2} we see that for any $w\in\WW^D_\l$ there is 
a unique element $c_w^\l\in H^D_\l$ such that $\ov{c_w^\l}=c_w^\l$ and 
$c_w^\l-\tT_w^\l\in\sum_{y\in\WW^D_\l;y\ne w}v\i\ZZ[v\i]\tT_y^\l$. Also, 
$\{c_w^\l;w\in\WW^D_\l\}$ is an $\ca$-basis of $H^D_\l$ and $\{c_w^\l;w\in\WW_\l\}$ is 
an $\ca$-basis of $H_\l$ (as in \cite{\KL}). 

\proclaim{Lemma 34.3}The $\ca$-algebra $H^D_\l$ with its basis 
$(c_w^\l)_{w\in\WW^D_\l}$ satisfies $P_1,P_2,P_3$ in 34.1.
\endproclaim
The analogous statement where $H^D_\l,\WW^D_\l$ are replaced by $H_\l,\WW_\l$ holds by
\cite{\HA, \S15}. The proof of the lemma is entirely similar; alternatively, it can be
reduced to the case of $H_\l$ using the identities 

(a) $c_{w_1w_2}^\l=\tT_{w_1}^\l c_{w_2}^\l$, $c_{w_2w_1}^\l=c_{w_2}^\l\tT_{w_1}^\l$ for
$w_1\in\O^D_\l,w_2\in\WW_\l$,

(b) $\tT_{w_1}^\l\tT_{w'_1}^\l=\tT_{w_1w'_1}^\l$ for $w_1\in\O^D_\l,w'_1\in\O^D_\l$. 

\mpb

The function $a:\{c_w^\l;w\in\WW^D_\l\}@>>>\NN$ (see 34.1) is determined by the 
analogous function $a:\{c_w^\l;w\in\WW_\l\}@>>>\NN$ (defined in terms of $H_\l$) by
$a(c_{w_1w_2}^\l)=a(c_{w_2w_1}^\l)=a(c_{w_2}^\l)$ for $w_1\in\O^D_\l,w_2\in\WW_\l$. The
two-sided cells of $\{c_w^\l;w\in\WW^D_\l\}$ are the sets of the form 
$\tT_{w_1}^\l\boc\tT_{w'_1}^\l$ where $w_1,w'_1$ run through $\O^D_\l$ and $\boc$ is a 
two-sided cell of $\{c_w^\l;w\in\WW_\l\}$. We show:

(c) {\it If $c_w^\l (w\in\WW^D_\l)$ is a distinguished basis element of $H^D_\l$ (see
34.1) then $w\in\WW_\l$ and $w^2=1$.}
\nl
By \cite{\HA} any left cell of $\WW_\l$ contains a unique distinguished basis element.
By the same argument, any left cell of $\WW^D_\l$ contains a unique distinguished basis
element. Let $\G$ be the left cell of $\WW^D_\l$ that contains $c_w^\l$. (See 3.1.) 
Write $w=w_1w_2$ with $w_1\in\O^D_\l,w_2\in\WW_\l$. From (a) we have
$c_w^\l=\tT_{w_1}^\l c_{w_2}^\l,c_{w_2}^\l=\tT_{w_1\i}^\l c_{w^\l}$. Hence
$c_{w_2}^\l\in\G$. We see that, if $\G'$ is the left cell of $\WW_\l$ that contains 
$c_{w_2}^\l$, then $\G'\sub\G$. Let $c_{w_3}^\l$ be the unique distinguished basis
element of $\WW_\l$ that is contained in $\G'$. Then $c_{w_3}^\l$ is also a 
distinguished basis element of $\WW^D_\l$ contained in $\G$ hence, by uniqueness, we
have $c_{w_3}^\l=c_w^\l$. We see that $w=w_3\in\WW_\l$. The fact that $w^2=1$ also
follows from \cite{\HA}.

\subhead 34.4\endsubhead
Let $H_n^D$ be the $\ca$-algebra with $1$ defined by the generators 
$\tT_w(w\in\WW^D)$, $1_\l(\l\in\ufs_n)$ and the relations

$1_\l1_\l=1_\l$ for $\l\in\ufs_n$,  $1_\l1_{\l'}=0$ for $\l\ne\l'$ in $\ufs_n$,

$\tT_w\tT_{w'}=\tT_{ww'}$ for $w,w'\in\WW^D$ with $l(ww')=l(w)+l(w')$,

$\tT_w1_\l=1_{w\l}\tT_w$ for $w\in\WW^D,\l\in\ufs_n$,

$\tT_s^2=\tT_1+(v-v\i)\sum_{\l;s\in\WW_\l}\tT_s1_\l$ for $s\in\II$,

$\tT_1=\sum_\l1_\l$.
\nl
We identify $H_n$ (see 31.2) with the subalgebra of $H_n^D$ generated by
$\tT_w(w\in\WW),1_\l(\l\in\ufs_n)$ by $T_w\m v^{l(w)}\tT_w(w\in\WW)$,
$1_\l\m 1_\l(\l\in\ufs_n)$.
There is a unique ring homomorphism $\bar{}:H^D_n@>>>H^D_n$ such that 
$\ov{\tT_w}=\tT_{w\i}\i$ for all $w\in\WW^D$, $\ov{v^m1_\l}=v^{-m}1_\l$ for all $\l$ 
and all $m\in\ZZ$. It has square $1$. Its restriction to $H_n$ is the involution
$\bar{}:H_n@>>>H_n$ described in 31.3. From the definitions, for any $w\in\WW^D$,
$\l\in\ufs_n$ we have 
$\ov{\tT_w1_\l}-\tT_w1_\l\in\sum_{y\in\WW^D;y\le w,y\ne w}\ca\tT_y1_\l$. By an argument
similar to one in \cite{\HA, 5.2} we see that for any $w\in\WW^D$, $\l\in\ufs_n$, there 
is a unique element $c_{w,\l}\in H^D_n$ such that $\ov{c_{w,\l}}=c_{w,\l}$ and 
$c_{w,\l}-\tT_w1_\l\in\sum_{y\in\WW^D;y\ne w}v\i\ZZ[v\i]\tT_y1_\l$. We have 
$c_{w,\l}\in 1_{w\l}H^D_n1_\l$. Also, $\{c_{w,\l};w\in\WW^D,\l\in\ufs_n\}$ is an 
$\ca$-basis of $H^D_n$. 

\proclaim{Proposition 34.5}The $\ca$-algebra $H^D_n$ with its basis 
$(c_{w,\l})_{(w,\l)\in\WW^D\T\ufs_n}$ satisfies $P_1,P_2,P_3$ in 34.1.
\endproclaim
The proof is given in 34.10.

\subhead 34.6\endsubhead
In the setup of 34.2, the $\ca$-algebra $1_\l H^D_n1_\l$ (a subalgebra of $H^D_n$) has 
a unit element $1_\l$, an $\ca$-basis $\{\tT_w1_\l;w\in\WW^D_\l\}$ and an $\ca$-basis 
$\{c_{w,\l};w\in\WW^D_\l\}$. 

\proclaim{Lemma 34.7}The $\ca$-algebra $1_\l H^D_n1_\l$ with its basis 
$\{c_{w,\l};w\in\WW^D_\l\}$ satisfies $P_1,P_2,P_3$ in 34.1.
\endproclaim
Define $\vt_\l:H^D_\l@>>>1_\l H^D_n1_\l$ by $\tT^\l_w\m\tT_w1_\l$ (an isomorphism of 
$\ca$-modules). Using 34.3, we see that it suffices to show that $\vt_\l$ is an 
isomorphism of $\ca$-algebras carrying $c_w^\l$ to $c_{w,\l}$ for any $w\in\WW^D_\l$. 
We use the notation in 34.2. We show:

(a) {\it Let $w\in\WW^D,\a\in\Pi_\l,\s=s_\a\in\II_\l$. Then 
$\tT_w\tT_\s1_\l=\tT_{w\s}1_\l+\d(v-v\i)\tT_w1_\l$ (in $H^D_n$) with $\d\in\{0,1\}$. If
in addition $w\in\WW^D_\l$ then $\d=0$ if $l_\l(w\s)>l_\l(w)$ and $\d=1$, otherwise.}
\nl
The proof has some common features with one in \cite{\MS, 3.3.5}. We have 
$\s=s_1s_2\do s_r$ with $s_i\in\II,r=l(\s)$. By \cite{\CS, I, 5.3}, there exists 
$j\in[1,r]$ such that $s_r\do s_{j+1}s_js_{j+1}\do s_r\in\WW_\l$ and 
$s_r\do s_{i+1}s_is_{i+1}\do s_r\n\WW_\l$ for $i\in[1,r]-\{j\}$; by \cite{\CS, I, 5.6},
we have $\s=s_r\do s_{j+1}s_js_{j+1}\do s_r$. Hence 
$s_1s_2\do s_{j-1}s_{j+1}\do s_r=1$. From the relations of $H^D_n$ we have:
$$\tT_{ws_1s_2\do s_{j-1}}\tT_{s_j}1_{s_{j+1}\do s_r\l}=\tT_{ws_1s_2\do s_j}1_{s_{j+1}\do
s_r\l}+\d'(v-v\i)\tT_{ws_1s_2\do s_{j-1}}1_{s_{j+1}\do s_r\l}$$
where $\d'=0$ if $l(ws_1s_2\do s_j)>l(ws_1s_2\do s_{j-1})$ and $\d'=1$ otherwise,
$$\tT_{ws_1s_2\do s_{i-1}}\tT_{s_i}1_{s_{i+1}\do s_r\l}
=\tT_{ws_1s_2\do s_i}1_{s_{i+1}\do s_r\l}$$
for $i\in[1,r]-\{j\}$,
$$\tT_{ws_1s_2\do s_{j-1}s_{j+1}\do s_{i-1}}\tT_{s_i}1_{s_{i+1}\do s_r\l}=
\tT_{ws_1s_2\do s_{j-1}s_{j+1}\do s_i}1_{s_{i+1}\do s_r\l}$$ 
for $i\in[j+1,r]$. From these identities we see that
$$\align&\tT_w\tT_\s1_\l=\tT_w\tT_{s_1}\tT_{s_2}\do\tT_{s_r}1_\l=\tT_{ws_1s_2\do 
s_{j-1}}\tT_{s_j}\tT_{s_{j+1}}\do\tT_{s_r}1_\l\\&=\tT_{ws_1s_2\do s_j}\tT_{s_{j+1}}\do
\tT_{s_r}1_\l+\d'(v-v\i)\tT_{ws_1s_2\do s_{j-1}}\tT_{s_{j+1}}\do\tT_{s_r}1_\l\\&
=\tT_{ws_1s_2\do s_js_{j+1}\do s_r}1_\l+\d'(v-v\i)\tT_{ws_1s_2\do s_{j-1}s_{j+1}\do 
s_r}1_\l\\&=\tT_{w\s}1_\l+\d'(v-v\i)\tT_w1_\l.\endalign$$
Assume now that $w\in\WW^D_\l$. We show that $\d=\d'$. The condition that $\d=0$ is 
equivalent to the condition that $w(\a)\in\RR^+_\l$. The condition that $\d'=0$ is 
equivalent to the condition that $ws_1s_2\do s_{j-1}(\a_j)\in R^+$ where $\a_j\in R^+$ 
is defined by $s_j=s_{\a_j}$. Since $\a=s_1s_2\do s_{j-1}(\a_j)$, this completes the 
proof of (a).

We show:

(b) {\it Let $w\in\WW^D_\l,w'\in\O^D_\l$. Then 
$\tT_w\tT_{w'}1_\l=\tT_{ww'}1_\l\in H^D_n$.}
\nl
We write $w'=s_1s_2\do s_r$ with $s_i\in\II,r=l(w')$. Using \cite{\CS, I, 5.3} we see 
that $s_r\do s_{i+1}s_is_{i+1}\do s_r\n\WW_\l$ for all $i\in[1,r]$. From the relations 
of $H^D_n$ we have
$$\tT_{ws_1s_2\do s_{i-1}}\tT_{s_i}1_{s_{i+1}\do s_r\l}
=\tT_{ws_1s_2\do s_i}1_{s_{i+1}\do s_r\l}$$
for $i\in[1,r]$. Using these identities we see that
$\tT_w\tT_{w'}1_\l=\tT_w\tT_{s_1}\tT_{s_2}\do\tT_{s_r}1_\l=\tT_{ws_1s_2\do s_r}1_\l$ 
and (b) follows.

We show that $\vt^\l$ is an algebra homomorphism. We must check that 
$(\tT_\s1_\l)^2=1_\l+(v-v\i)\tT_\s1_\l$ for $s\in\II_\l$. This is a special case of (a)
(take $w=\s$). We must also check that $(\tT_w1_\l)(\tT_{w'}1_\l)=\tT_{ww'}1_\l$ if 
$w,w'\in\WW^D_\l$, $l_\l(ww')=l_\l(w)+l_\l(w')$. If $w,w'\in\WW_\l$, this is proved by 
induction on $l_\l(w')$, the induction step being provided by (a). The general case can
be reduced to this special case using (b). We see that $\vt_\l$ is an $\ca$-algebra
isomorphism. We show that
 
$\ov{\vt_\l(h)}=\vt_\l(\ov h)$ for $h\in H^D_\l$.
\nl
Assume first that $h=\tT^\l_w$ where $w\in\O^D_\l$. Then 
$\ov h=(\tT^\l_{w\i})\i=\tT^\l_w$. Hence
$$\ov{\vt_\l(h)}=\ov{\tT_w1_\l}=\tT_{w\i}\i1_\l=\tT_w\i1_\l=\vt_\l(\tT^l_w)
=\vt_\l(\ov h),$$
as required. Assume next that $h=\tT^\l_\s$ where $\s\in\II_\l$. Then 
$$\align&\vt_\l(\ov h)=\vt_\l((\tT^\l_\s)\i)=\vt_\l(\tT^\l_\s+(v\i-v)\tT^\l_1)=
\tT_\s1_\l+(v\i-v)\tT_11_\l=\tT_\s\i1_\l\\&=\ov{\tT_\s1_\l}=\ov{\vt_\l(h)},\endalign$$
as required.

We see that for $w\in\WW^D_\l$ we have $\ov{\vt_\l(c_w^\l)}=\vt_\l(c_w^\l)$. Hence
$\vt_\l(c_w^\l)$ satisfies the defining properties of $c_{w,\l}$ so that 
$\vt_\l(c_w^\l)=c_{w,\l}$. The lemma is proved.

\mpb

Using now 34.3(c) we see that

(c) {\it If $c_{w,\l} (w\in\WW^D_\l)$ is a distinguished basis element of 
$1_\l H^D_n1_\l$ (see 34.1) then $w\in\WW_\l$ and $w^2=1$.}

\subhead 34.8\endsubhead
Let $\ufs'_n$ be a set of representatives for the $\WW^D$-orbits in $\ufs_n$. For 
$\l\in\ufs_n$ define $\l^0\in\ufs'_n$ by $\l^0\in\WW^D\l$ (the $\WW^D$-orbit of $\l$).
Let 

$\G=\{(\l_1,\l_2)\in\ufs_n\T\ufs_n;\WW^D\l_1=\WW^D\l_2\}$.
\nl
Let $E^D_n$ be the set of all formal sums $x=\sum_{(\l_1,\l_2)\in\G}x_{\l_1,\l_2}$
where $x_{\l_1,\l_2}\in 1_{\l_1^0}H^D_n1_{\l_2^0}$. Then $E^D_n$ is naturally an 
$\ca$-module and an associative $\ca$-algebra where the product $xy$ of $x,y\in E^D_n$ 
is given by $(xy)_{\l_1,\l_2}=\sum_{\ti\l\in\WW^D\l_1}x_{\l_1,\ti\l}y_{\ti\l,\l_2}$. 
This algebra has a unit element, namely the element $1$ such that 
$1_{\l_1,\l_2}=\d_{\l_1,\l_2}1_{\l_1}$ for $(\l_1,\l_2)\in\G$. Define a ring involution
$\bar{}:E^D_n@>>>E^D_n$ by $x\m\ov x$ where $\ov x_{\l_1,\l_2}=\ov{x_{\l_1,\l_2}}$. 
(Note that $\bar{}:H^D_n@>>>H^D_n$ maps $1_{\l_1^0}H^D_n1_{\l_2^0}$ onto itself.)

Let $C=\{(\l_1,\l_2,w)\in\ufs_n\T\ufs_n\T\WW^D;w\l_1^0=\l_1^0=\l_2^0\}$. For
$(\l_1,\l_2,w)\in C$ define $x^{\l_1,\l_2,w}\in E^D_n$ by 

$x^{\l_1,\l_2,w}_{\l'_1,\l'_2}=\d_{(\l'_1,\l'_2),(\l_1,\l_2)}\tT_w1_{\l_1^0}$.
\nl
Then $\{x^{\l_1,\l_2,w};(\l_1,\l_2,w)\in C\}$ is an $\ca$-basis of $H''_n$. From the 
definitions, for $(\l_1,\l_2,w)\in C$ we have

$\ov{x^{\l_1,\l_2,w}}-x^{\l_1,\l_2,w}\in
\sum_{y\in\WW^D;y\le w,y\ne w,y\l_1^0=\l_1^0}\ca x^{\l_1,\l_2,y}$ 
\nl
By an argument similar to one in \cite{\HA, 5.2} we see that for any \lb
$(\l_1,\l_2,w)\in C$ there is a unique element $c^{\l_1,\l_2,w}\in E^D_n$ such that 
$\ov{c^{\l_1,\l_2,w}}=c^{\l_1,\l_2,w}$ and 
$$c^{\l_1,\l_2,w}-x^{\l_1,\l_2,w}\in
\sum_{y\in\WW^D,y\l_1^0=\l_1^0,y\ne w}v\i\ZZ[v\i]x^{\l_1,\l_2,y}.$$
Also, $\{c^{\l_1,\l_2,w};(\l_1,\l_2,w)\in C\}$ is an $\ca$-basis of $E^D_n$. 

\proclaim{Lemma 34.9}The $\ca$-algebra $E^D_n$ with its basis 
$\{c^{\l_1,\l_2,w};(\l_1,\l_2,w)\in C\}$ satisfies $P_1,P_2,P_3$ in 34.1.
\endproclaim
For $\l\in\ufs'_n$ let $N_\l=|\WW^D\l|$ and let $M_{N_\l}(1_\l H^D_n1_\l)$ be the 
algebra of $N_\l\T N_\l$ matrices with entries in $1_\l H^D_n1_\l$. From the 
definitions we have a decomposition

$E^D_n=\op_{\l\in\ufs'_n}M_{N_\l}(1_\l H^D_n1_\l)$
\nl
compatible with the algebra structures and with the natural bases. Using this, the 
lemma is reduced to the similar statement for $1_\l H^D_n1_\l$ where it is known by 
34.7.

\mpb

The function $a:\{c^{\l_1,\l_2,w};(\l_1,\l_2,w)\in C\}$ (see 34.1) is given 
by $a(c^{\l_1,\l_2,w})=a(c_{w,\l_1^0})$ where $a(c_{w,\l_1^0})$ is defined as in 34.1 
in terms of $1_{\l_1^0}H^D_n1_{\l_1^0}$. The two-sided cells of 
$\{c^{\l_1,\l_2,w};(\l_1,\l_2,w)\in C\}$ are the sets of the form $\{c^{\l_1,\l_2,w}\}$
where $\l_1,\l_2$ run through $\WW^D\l$ (with $\l\in\ufs'_n$ fixed) and $w$ running 
through a subset $X$ of $\WW^D_\l$ such that $\{c_{w,\l};w\in X\}$ is a two-sided cell 
of $\{c_{w,\l};w\in\WW^D_\l\}$ (see 34.7).

Using 34.7(c) we obtain:

(a) {\it If $c^{\l_1,\l_2,w},(\l_1,\l_2,w)\in C)$ is a distinguished basis element of
$E^D_n$ then $\l_1=\l_2,w\in\WW_{\l_1^0},w^2=1$.} 

\subhead 34.10\endsubhead
We prove Proposition 34.5. It is enough to construct an algebra isomorphism 
$H^D_n@>\si>>E^D_n$ which carries the basis $(c_{w,\l})$ onto the basis 
$(c^{\l_1,\l_2,w})$.

For each $\l\in\ufs_n$ we choose a sequence $\ss_\l=(s_1,s_2,\do,s_r)$ where, for 
$i\in[1,r]$, $s_i$ is either in $\II$ or is a power of $\uD$ and 
$\l^0=s_1s_2\do s_r\l\ne s_2\do s_r\l\ne\do\ne s_r\l\ne\l$ or, equivalently,
$\l=s_r\i\do s_2\i s_1\i\l^0\ne s_{r-1}\i\do s_1\i\l^0\ne\do\ne s_1\i\l^0\ne\l^0$. Let 
$[\ss_\l]=s_1s_2\do s_r$. We set 

$\t_\l=\tT_{s_1}\tT_{s_2}\do\tT_{s_r}\in H^D_n$, 
$\t'_\l=\tT_{s_r\i}\do\tT_{s_2\i}\tT_{s_1\i}\in H^D_n$.
\nl
We show:

(a) $1_{\l^0}\t_\l\t'_\l=1_{\l^0},1_\l\t'_\l\t_\l=1_\l$.
\nl
We have 
$$1_{\l^0}\t_\l\t'_\l=1_{\l^0}\tT_{s_1}\tT_{s_2}\do\tT_{s_r}\tT_{s_r\i}\do\tT_{s_2\i}
\tT_{s_1\i}=\tT_{s_1}\tT_{s_2}\do\tT_{s_r}1_\l\tT_{s_r\i}\do\tT_{s_2\i}\tT_{s_1\i}.$$
Since $s_r\l\ne\l$, we can replace $\tT_{s_r}1_\l\tT_{s_r\i}$ by $1_{s_r\l}$ and we 
obtain 

$\tT_{s_1}\tT_{s_2}\do\tT_{s_{r-1}}1_{s_r\l}\tT_{s_{r-1}\i}\do\tT_{s_2\i}\tT_{s_1\i}$.
\nl
Since $s_{r-1}s_r\l\ne s_r\l$, we can replace $\tT_{s_{r-1}}1_{s_r\l}\tT_{s_{r-1}\i}$ 
by $1_{s_{r-1}s_r\l}$ and we obtain 

$\tT_{s_1}\tT_{s_2}\do\tT_{s_{r-2}}1_{s_{r-1}s_r\l}\tT_{s_{r-2}\i}\do\tT_{s_2\i}
\tT_{s_1\i}.$
\nl
Continuing in this way we find $1_{s_1\do s_{r-1}s_r\l}=1_{\l^0}$. This proves the 
first identity in (a). The second identity is proved in a similar way. 

We have 

(b) $\ov{\t_\l1_\l}=\t_\l1_\l,\ov{1_\l\t'_\l}=1_\l\t'_\l$.
\nl
The first identity in (b) is equivalent to 
$\tT_{s_1\i}\i\tT_{s_2\i}\i\do\tT_{s_r\i}\i1_\l=\t_\l1_\l$ or to 
$\t'_\l{}\i1_\l=\t_\l1_\l$, which follows from (a). Similarly, the second identity in 
(b) follows from (a).

We define an $\ca$-linear map $\Ps:H^D_n@>>>E^D_n$ by 

$\Ps(h)_{\l_1,\l_2}=\t_{\l_1}1_{\l_1}h1_{\l_2}\t'_{\l_2}\in 1_{\l_1^0}H^D_n1_{\l_2^0}$.
\nl
We show that $\Ps$ is a ring homomorphism. Let $h,h'\in H^D_n$, $x=\Ps(h),y=\Ps(h')$,
$z=\Ps(hh')$, $z'=\Ps(h)\Ps(h')$. We have
$$\align&(\Ps(h)\Ps(h'))_{\l_1,\l_2}=
\sum_{\ti\l\in\WW^D\l_1}\Ps(h)_{\l_1,\ti\l}\Ps(h')_{\ti\l,\l_2}\\&
=\sum_{\ti\l\in\WW^D\l_1}\t_{\l_1}1_{\l_1}h1_{\ti\l}\t'_{\ti\l}\t_{\ti\l}1_{\ti\l}
h'1_{\l_2}\t'_{\l_2}=
\sum_{\ti\l\in\WW^D\l_1}\t_{\l_1}1_{\l_1}h1_{\ti\l}h'1_{\l_2}\t'_{\l_2}\endalign$$ 
where the last equality comes from (a). Since $1_{\l_1}h1_{\ti\l}=0$ if 
$\ti\l\in\ufs_n-\WW^D\l_1$, we see that
$$(\Ps(h)\Ps(h'))_{\l_1,\l_2}=\t_{\l_1}1_{\l_1}h\sum_{\ti\l\in\ufs_n}1_{\ti\l}h'
1_{\l_2}\t'_{\l_2}=\t_{\l_1}1_{\l_1}hh'1_{\l_2}\t'_{\l_2}=\Ps(hh')_{\l_1,\l_2}.$$
Thus $\Ps(h)\Ps(h')=\Ps(hh')$, as required.

We show that 

(c) $\ov{\Ps(h)}=\Ps(\ov h)$ for $h\in H^D_n$.
\nl
We have 

$(\ov{\Ps(h)})_{\l_1,\l_2}=\ov{\t_{\l_1}1_{\l_1}h1_{\l_2}\t'_{\l_2}}$,
$(\Ps(\ov h))_{\l_1,\l_2}=\t_{\l_1}1_{\l_1}\bar h1_{\l_2}\t'_{\l_2}$.
\nl
It suffices to show that $\ov{\t_{\l_1}1_{\l_1}}=\t_{\l_1}1_{\l_1}$, 
$\ov{1_{\l_2}\t'_{\l_2}}=1_{\l_2}\t'_{\l_2}$. This follows from (b).

We show that

(d) $\Ps(\tT_w1_\l)=x^{w\l,\l,[\ss_{w\l}]w[\ss_\l]\i}$ for $w\in\WW^D,\l\in\ufs_n$.
\nl
Indeed, $(\Ps(\tT_w1_\l))_{\l_1,\l_2}=\t_{\l_1}1_{\l_1}\tT_w1_\l1_{\l_2}\t'_{\l_2}$. 
This is $0$ unless $\l_2=\l,\l_1=w\l$. If $\l_2=\l,\l_1=w\l$ we see as in the proof of 
(a) that 

$\t_{\l_1}1_{\l_1}\tT_w1_\l1_{\l_2}\t'_{\l_2}=
\t_{w\l}\tT_w1_\l\t'_\l=\tT_{[\ss_{w\l}]w[\ss_\l]\i}1_{\l^0}$.
\nl
This proves (d).

(d) shows that $\Ps$ is induced by a map 

$\Ps_0:\WW^D\T\ufs_n@>>>C$, $(w,\l)\m(w\l,\l,[\ss_{w\l}]w[\ss_\l]\i)$. 
\nl
This is a bijection with inverse $(\l_1,\l_2,w)\m([\ss_{\l_1}]\i w[\ss_{\l_2}],\l_2)$).
It follows that $\Ps$ is an isomorphism.

Let $w\in\WW^D,\l\in\ufs_n$. Since $\Ps\i$ is compatible with $\bar{}$, we have

$\ov{\Ps\i(c^{w\l,\l,[\ss_{w\l}]w[\ss_\l]\i})}
=\Ps\i(c^{w\l,\l,[\ss_{w\l}]w[\ss_\l]\i})$.
\nl
From (d) we see that 

$\Ps\i(c^{w\l,\l,[\ss_{w\l}]w[\ss_\l]\i})-\tT_w1_\l\in
\sum_{y\in\WW^D,y\l=w\l;y\ne w}v\i\ZZ[v\i]\tT_y1_\l$.
\nl
Since these properties characterize $c_{w,\l}$, we see that 
$\Ps\i(c^{w\l,\l,[\ss_{w\l}]w[\ss_\l]\i})=c_{w,\l}$ that is,
$\Ps(c_{w,\l})=c^{w\l,\l,[\ss_{w\l}]w[\ss_\l]\i}$. Thus, $\Ps$ restricts to a bijection
between the basis $(c_{w,\l})$ of $H^D_n$ and the basis $(c^{\l_1,\l_2,w})$ of 
$E^D_n$, induced by the bijection $\Ps_0:\WW^D\T\ufs_n@>\si>>C$. Proposition 34.5 is 
proved.

\mpb

We show:

(e) {\it If $c_{w,\l} (w\in\WW^D,\l\in\ufs_n)$ is a distinguished basis element of
$H^D_n$ then $w\in\WW_\l,w^2=1$.}
\nl
Note that, with the notation above, $c^{w\l,\l,[\ss_{w\l}]w[\ss_\l]\i}$ is a 
distinguished basis element of $E^D_n$ hence, by 34.9(a), we have $w\l=\l$ and 
$[\ss_{w\l}]w[\ss_\l]\i\in\WW_{\l^0}$ has square $1$. Thus,
$[\ss_\l]w[\ss_\l]\i\in\WW_{\l^0}$ has square $1$. It follows that
$w\in\WW_{[\ss_\l]\i\l^0}=\WW_\l$ has square $1$. This proves (e).

\subhead 34.11\endsubhead
The algebra $H_n^{G^0}$ defined as in 34.4 with $D$ replaced by $G^0$ is the same as
$H_n$; we identify it in an obvious way with a subalgebra of $H_n^D$ (see 34.4). The 
$\ca$-basis of $H_n=H_n^{G^0}$ analogous to the $\ca$-basis 
$\{c_{w,\l};w\in\WW^D,\l\in\ufs_n\}$ of $H_n^D$ is the subset of the last basis given 
by $\{c_{w,\l};w\in\WW,\l\in\ufs_n\}$. The analogue of Proposition 34.5 holds: the 
$\ca$-algebra $H_n$ with its basis $(c_{w,\l})_{(w,\l)\in\WW\T\ufs_n}$ satisfies 
$P_1,P_2,P_3$ in 34.1. 

\subhead 34.12\endsubhead
Let $K$ be a field of characteristic $0$ and let $\ca@>>>K$ be a homomorphism of rings 
with $1$ which carries $v\in\ca$ to $v_0\in K^*$. We show:

(a) {\it If $\l\in\ufs_n$ and $\sum_{w\in\WW_\l}v_0^{2l_\l(w)}\ne 0$, then the 
$K$-algebra $H^D_{\l,K}=K\ot_\ca H^D_\l$ is semisimple.}
\nl
Let $M$ be an $H^D_{\l,K}$-module of finite dimension over $K$ and let $M'$ be an
$H^D_{\l,K}$-submodule of $M$. It is enough to show that there exists an 
$H^D_{\l,K}$-submodule of $M$ complementary to $M'$. It is well known that under our 
assumption, the $K$-algebra $H_{\l,K}=K\ot_\ca H_\l$ is semisimple. Hence there exists 
an $H_{\l,K}$-submodule of $M$ complementary to $M'$ that is, there exists an 
$H_{\l,K}$-linear map $f:M@>>>M'$ such that $f(x)=x$ for all $x\in M'$. Define
$\tf:M@>>>M'$ by $\tf(x)=|\O^D_\l|\i\sum_{z\in\O^D_\l}\tT_{z\i}f(\tT_zx)$. For
$w\in\WW_\l,x\in M$ we have
$$\align&|\O^D_\l|\tf(\tT_wx)=\sum_{z\in\O^D_\l}\tT_{z\i}f(\tT_z\tT_wx)=
\sum_{z\in\O^D_\l}\tT_{z\i}f(\tT_{zwz\i}\tT_zx)\\&=\sum_{z\in\O^D_\l}\tT_{z\i}
\tT_{zwz\i}f(\tT_zx)=\sum_{z\in\O^D_\l}\tT_w\tT_{z\i}f(\tT_zx)=|\O^D_\l|\tT_w\tf(x).
\endalign$$
(The third equality holds since $f$ is $H_{\l,K}$-linear.) We see that $\tf$ is 
$H_{\l,K}$-linear. Since $\tf(\tT_yx)=\tT_y\tf(x)$ for $x\in M$, $y\in\O^D_\l$, we see 
that $\tf$ is $H^D_{\l,K}$-linear. If $x\in M'$ we have $\tT_zx\in M'$ for 
$z\in\O^D_\l$ hence $\tf(x)=|\O^D_\l|\i\sum_{z\in\O^D_\l}\tT_{z\i}\tT_zx=x$. It follows
that the $\ker(\tf)$ is an $H^D_{\l,K}$-submodule of $M$ complementary to $M'$. This
proves (a).

{\it Let $\fU$ be the subfield of $\bbq$ generated by the roots of $1$.}
\nl
For any $\k\in\fU^*,\l\in\ufs_n$, let 

$H_n^{D,\k}=\fU\ot_\ca H^D_n$, $H_n^\k=\fU\ot_\ca H_n$, 
$H_\l{}^{D,\k}=\fU\ot_\ca H^D_\l$,
\nl
where $\fU$ is regarded as an $\ca$-algebra via the ring homomorphism 
$\ca@>>>\fU,v\m\k$. Let 

$H_n^{D,v}=\fU(v)\ot_\ca H^D_n$, $H_n^v=\fU(v)\ot_\ca H_n$ 
\nl
where $\fU(v)$ (field of rational functions in $v$ with coefficients in $\fU$) is 
regarded as an $\ca$-algebra via the ring homomorphism $\ca@>>>\fU(v),v\m v$. Now 
$\Ph:H^D_n@>>>\ca\ot_\ZZ H^{D,\iy}_n$ ( where $H^{D,\iy}_n=(H^D_n)^\iy$ is defined as 
in 34.1 in terms of the basis $(c_{w,\l})$ in 34.4) induces algebra homomorphisms 
$\Ph^\k:H_n^{D,\k}@>>>\fU\ot_\ZZ H^{D,\iy}_n$ for $\k\in\fU^*$ and 
$\Ph^v:H_n^{D,v}@>>>\fU(v)\ot_\ZZ H^{D,\iy}_n$. 

Similarly, $\Ph:H_n@>>>\ca\ot_\ZZ H_n^\iy$ (defined as in 34.1 in terms of the basis
$(c_{w,\l})$ in 34.11) induces algebra homomorphisms $H_n^\k@>>>\fU\ot_\ZZ H_n^\iy$ for
$\k\in\fU^*$ and $H_n^v@>>>\fU(v)\ot_\ZZ H_n^\iy$, denoted again by $\Ph^\k,\Ph^v$. 
From the definitions we see that $H_n^\iy$ may be identified with the subgroup of
$(H_n^D)^\iy$ spanned by the basis elements of $(H_n^D)^\iy$ indexed by 
$\{(w,\l);w\in\WW,\l\in\ufs_n\}$ and $\Ph:H_n@>>>\ca\ot_\ZZ H_n^\iy$ becomes the 
restriction of $\Ph:H^D_n@>>>\ca\ot_\ZZ H^{D,\iy}_n$. We show:

(b) {\it If $\k\in\fU^*,\sum_{w\in\WW}\k^{2l(w)}\ne 0$ then $H_n^{D,\k}$ is a 
semisimple $\fU$-algebra and $\Ph^\k:H_n^{D,\k}@>>>\fU\ot_\ZZ H^{D,\iy}_n$ is an 
algebra isomorphism. Moreover, $H_n^\k$ is a semisimple $\fU$-algebra and 
$\Ph^\k:H_n^\k@>>>\fU\ot_\ZZ H_n^\iy$ is an algebra isomorphism.}

(c) {\it $H_n^{D,v}$ is a split semisimple $\fU(v)$-algebra and 
$\Ph^v:H_n^{D,v}@>>>\fU(v)\ot_\ZZ H^{D,\iy}_n$ is an algebra isomorphism. Moreover, 
$H_n^v$ is a split semisimple $\fU(v)$-algebra and 
$\Ph^v:H_n^v@>>>\fU(v)\ot_\ZZ H_n^\iy$ is an algebra isomorphism.}
\nl
We prove (b). The following statement is easily verified:

{\it If $\l\in\ufs_n$ then $\sum_{w\in\WW}v^{2l(w)}=Q\sum_{w\in\WW_\l}v^{2l_\l(w)}$ for
some $Q\in\ZZ[v^2]$.}
\nl
We see that, if $\k$ is as in (b) and $\l\in\ufs_n$ then 
$\sum_{w\in\WW_\l}\k^{2l_\l(w)}\ne 0$; hence, by (a), $H_\l^{D,\k}$ is a semisimple 
algebra. By the arguments in 34.7-34.10, $H_n^{D,\k}$ is a direct sum of matrix rings 
over rings of the form $H_\l^{D,\k}$. Hence $H_n^{D,\k}$ is a semisimple algebra. Using
34.1(c) we see that $\Ph^\k:H_n^{D,\k}@>>>\fU\ot_\ZZ H^{D,\iy}_n$ is an algebra 
isomorphism. It remains to show that $H_n^{D,\k}$ is split over $\fU$. Since $\Ph^\k$ 
is an isomorphism it is enough to show that $\fU\ot_\ZZ H^{D,\iy}_n$ is split over 
$\fU$. Since $\Ph^1$ is an isomorphism it is enough to show that $H_n^{D,1}$ is split 
over $\fU$. As above, $H_n^{D,1}$ is a direct sum of matrix rings over rings of the 
form $H_\l^{D,1}$. Since $H_\l^{D,1}$ is the group algebra of a finite group with 
coefficients in $\fU$, it is split over $\fU$ by Brauer's theorem. This proves the 
first sentence in (b). The second sentence in (b) is obtained from the first by 
replacing $D$ by $G^0$.

Now the proof of (c) is just like that of (b) except for the splitness assertion. By 
(b), $\fU\ot_\ZZ H^{D,\iy}_n$ is a split semisimple $\fU$-algebra hence
$\fU(v)\ot_\ZZ H^{D,\iy}_n$ is a split semisimple $\fU(v)$-algebra. Since $\Ph^v$ is an
isomorphism, it follows that $H_n^{D,v}$ is split over $\fU(v)$. Similarly, $H_n^v$ is 
split over $\fU(v)$. This proves (c).

\subhead 34.13\endsubhead
Define an $\ca$-linear map $\t:H_n@>>>\ca$ by $\t(\tT_w1_\l)=\d_{w,1}$ for all 
$w\in\WW,\l\in\ufs_n$. Define a bilinear form $(,):H_n\T H_n@>>>\ca$ by 
$(x,x')=\t(xx')$. We show that

(a) $(\tT_w1_\l,\tT_{w'}1_{\l'})=\d_{w\i,w'}\d_{\l,w'\l'}$
\nl
for $w,w'\in\WW,\l,\l'\in\ufs_n$. (This shows that $(,)$ is symmetric; indeed, we have 
$\d_{w\i,w'}\d_{\l,w'\l'}=\d_{w'{}\i,w}\d_{\l',w\l}$.) To prove (a) it suffices to show
that, for $w,w'\in\WW,\l\in\ufs_n$, we have $\t(\tT_w\tT_{w'}1_\l)=\d_{ww',1}$. We argue
by induction on $l(w)$. If $l(w)=0$ then $\tT_w\tT_{w'}=\tT_{ww'}$ and the result is 
clear. Assume now that $l(w)\ge 1$. We can find $s\in\II$ such that $l(w)=l(ws)+1$. 
Then $\t(\tT_w\tT_{w'}1_\l)=\t(\tT_{ws}\tT_s\tT_{w'}1_\l)$. If $l(sw')=l(w')+1$ then, 
by the induction hypothesis, 

$\t(\tT_{ws}\tT_s\tT_{w'}1_\l)=\t(\tT_{ws}\tT_{sw'}1_\l)=\d_{wssw',1}=\d_{ww',1}$,
\nl
as required. Assume now that $l(sw')=l(w')-1$. We have 
$$\align&\t(\tT_{ws}\tT_s\tT_{w'}1_\l)=\t(\tT_{ws}\tT_s\tT_s\tT_{sw'}1_\l)\\&=
\t(\tT_{ws}\tT_{sw'}1_\l)+(v-v\i)\sum_{\l';s\in\WW_{\l'}}
\t(\tT_{ws}\tT_s1_{\l'}\tT_{sw'}1_\l).\endalign$$
If $s\n\WW_{sw'\l}$ this equals (by the induction hypothesis) 
$\d_{wssw',1}=\d_{ww',1}$, as required; if $s\in\WW_{sw'\l}$ this equals (by the 
induction hypothesis)

$\d_{wssw',1}+(v-v\i)\t(\tT_{ws}\tT_{w'}1_\l)=\d_{ww',1}+(v-v\i)\d_{wsw',1}$.
\nl
It remains to note that $wsw'\ne 1$ whenever $l(w)=l(ws)+1$, $l(sw')=l(w')-1$. This 
proves (a).

\subhead 34.14\endsubhead
Let $\fC$ be a finite dimensional semisimple split (associative) algebra with $1$ over 
a field $K$. Let $\{E_u;u\in\cu\}$ be a set or representatives for the isomorphism 
classes of simple $\fC$-modules. Let $\fa:\fC@>>>\fC$ be an algebra automorphism. For 
$u\in\cu$, $c:e\m\fa(c)e$ defines a $\fC$-module structure on the $K$-vector space 
$E_u$ which is isomorphic to $E_{\bu}$ for a unique $\bu\in\cu$. Then $u\m\bu$ is a 
permutation of $\cu$. Let $\cu^\fa=\{u\in\cu;u=\bu\}$. For $u\in\cu^\fa$ we can find a 
$K$-linear isomorphism $\fa_u:E_u@>>>E_u$ such that $\fa_u(ce)=\fa(c)\fa_ue$ for all 
$c\in\fC,e\in E_u$. Note that $\fa_u$ is uniquely determined up to multiplication by an
element in $K^*$. We show:

(a) {\it If $c,c'\in\fC$, then the trace of the $K$-linear map $\fC@>>>\fC$,
$c_1\m c\fa(c_1)c'$ equals} $\sum_{u\in\cu^\fa}\tr(c\fa_u,E_u)\tr(\fa_u\i c',E_u)$.
\nl
Under the algebra isomorphism 

(b) $\fC@>\si>>\op_{u\in\cu}\End_K(E_u)$, $c\m[e\m ce,e\in E_u]$, 
\nl
the linear map $\fC@>>>\fC$ in (a) corresponds to an endomorphism of 
$\op_{u\in\cu}\End_K(E_u)$ which permutes the summands according to $u\m\bu$ and whose 
restriction to a summand with $u=\bu$ is $\End_K(E_u)\m\End_K(E_u)$, 
$f\m c\fa_uf\fa_u\i c'$. From this (a) follows easily. (Compare 20.3(b).)

We show:

(c) {\it If $y:\fC@>>>K$ is $K$-linear and $y(cc')=y(c'\fa(c))$ for all $c,c'\in\fC$
then there exist $b_u\in K (u\in\cu^\fa)$ such that 
$y(c)=\sum_{u\in\cu^\fa}b_u\tr(c\fa_u,E_u)$ for all $c\in\fC$.}
\nl
Let $\fC_u$ be the inverse image of the summand $\End_K(E_u)$ under (b). Then 
$\fC=\op_u\fC_u$ and $\fa(\fC_u)=\fC_{\bu}$ for $u\in\cu$. Let $y_u:\fC_u@>>>K$ be the 
restriction of $y$ to $\fC_u$. Let $c\in\fC_u$ where $u\ne\bu$. Let $c'\in\fC_u$ be the
projection of $1\in\fC$ onto $\fC_u$. We have $\fa(c)\in\fC_{\bu}$ hence $c'\fa(c)=0$. 
Also, $cc'=c\in\fC_u$. Thus, $y_u(c)=y_u(cc')=y(cc')=y(c'\fa(c))=0$. We see that 
$y_u=0$. We are reduced to the case where $\cu=\cu^\fa$ consists of a single element 
$u$ and $\fC=\End_K(E_u)$. We can find $h\in\fC$, invertible, such that $\fa(c)=hch\i$ 
for all $c\in\fC$. We can assume that $\fa_ue=he$ for all $e\in E_u$. We have 
$y(cc')=y(c'hch\i)$ for all $c,c'\in\fC$. Define $\ty:\fC@>>>K$ by $\ty(c)=y(ch\i)$. We
have $y(cc'h\i)=y(c'h\i hch\i)=y(c'ch\i)$ hence $\ty(cc')=\ty(c'c)$ for all 
$y,y'\in\fC$. Thus there exists $b\in K$ such that $\ty(c)=b\tr(c,E_u)$. Then 
$y(c)=b\tr(ch:E_u@>>>E_u)=b\tr(c\fa_u:E_u@>>>E_u)$. This proves (c).

Assume now that we are given a $K$-linear map $z:\fC@>>>K$ such that 
$(c,c')=z(cc')=z(c'c)$ is a nondegenerate (symmetric) $K$-bilinear form 
$\fC\T\fC@>>>K$. Let $(c_i)_{i\in I}$ be a $K$-basis of $\fC$. Define a $K$-basis 
$(c'_i)_{i\in I}$ of $\fC$ by $(c_i,c'_j)=\d_{ij}$. 

For $u\in\cu$ and $c\in\fC$ invertible, $\sum_{i\in I}\tr(c_ic,E_u)c\i c'_i$ is in the 
centre of $\fC$; if $u'\in\cu,u'\ne u$, this sum acts on $E_{u'}$ as $0$ and on 
$E_u$ as $f_u$ times the identity, where $f_u\in K^*$ is independent of 
$(c_i),(c'_i),c$. (We apply \cite{\HA, 19.2} to the dual bases $(c_ic),(c\i c'_i)$: we
have $(c_ic,c\i c'_j)=(c_i,cc\i c'_j)=(c_i,c'_j)=\d_{ij}$.) We see that for 
$u,u'\in\cu$ we have
$$\sum_{i\in I}\tr(c_ic,E_u)\tr(c\i c'_i,E_{u'})=\d_{u,u'}f_u\dim E_u.\tag d$$
Now assume that $u,u'\in\cu^\fa$. We can pick $c\in\fC$ invertible such that $c$ acts 
on $E_u$ as $\fa_u$ and on $E_{u'}$ as $\fa_{u'}$. From (d) we deduce
$$\sum_{i\in I}\tr(c_i\fa_u,E_u)\tr(\fa_{u'}\i c'_i,E_{u'})=\d_{u,u'}f_u\dim E_u.\tag e
$$

\subhead 34.15\endsubhead
We write $\fa:H_n@>>>H_n$ instead of $\fa_D:H_n@>>>H_n$ (see 31.4); this is the algebra
automorphism given by $h\m\tT_{\uD}h\tT_{\uD\i}$ (product in $H_n^D$) for $h\in H_n$.
The same formula defines an algebra automorphism of $H_n^\k$ or $H_n^v$ denoted again 
by $\fa$. From the definitions we see that $\fa:H_n@>>>H_n$ takes $c_{w,\l}$ to 
$c_{\e_D(w),\uD\l}$ for $w\in\WW,\l\in\ufs_n$. Hence it induces a ring automorphism 
$H_n^\iy@>>>H_n^\iy$ denoted again by $\fa$. It also induces algebra automorphisms 
$H_n^\k@>>>H_n^\k,\k\in\fU^*$ and $H_n^v@>>>H_n^v$ denoted again by $\fa$. Now 
$\Ph^\k:H_n^\k@>>>\fU\ot_\ZZ H_n^\iy$ for $\k\in\fU^*$ and 
$\Ph^v:H_n^v@>>>\fU(v)\ot_\ZZ H_n^\iy$ (see 34.11) are compatible with $\fa$. 

Let $\{E_u;u\in\cu\}$ be a set of representatives for the isomorphism classes of simple
modules for $H_n^1$ (a split semisimple $\fU$-algebra, by 31.12(b).) Define 
$\cu@>>>\cu,u\m\bu$ as in 34.14, replacing $(\fC,\fa$) by $(H_n^1,\fa)$. Let 
$\cu^\fa=\{u\in\cu;u=\bu\}$. 

Let $u\in\cu$. Clearly, if $E_u$ extends to an $H_n^{D,1}$-module then $u\in\cu^\fa$.
Conversely, we show that, if $u\in\cu^\fa$, then $E_u$ extends to an 
$H_n^{D,1}$-module. Since $H_n^1,H_n^{D,1}$ are split over $\fU$ it is enough to prove 
the analogous statement in which $E_u,H_n^1,H_n^{D,1}$ are replaced by 
$\fU'\ot_\fU E_u,\fU'\ot_\fU H_n^1,\fU'\ot_\fU H_n^{D,1}$ and $\fU'$ is an algebraic 
closure of $\fU$. Since $u\in\cu^\fa$ we can find a $\fU'$-linear isomorphism 
$X:\fU'\ot_\fU E_u@>>>\fU'\ot_\fU E_u$ such that 
$X(ce)=\fa(c)X(e)$ for all $c\in\fU'\ot_\fU H_n^1,e\in\fU'\ot_\fU E_u$. Let $k$ be the 
order of $\uD:\TT@>>>\TT$. By Schur's lemma, $X^k$ is a scalar times identity. Now 
$\fU'$ contains a $k$-th root of this scalar. Dividing $X$ by this root we see that we 
may assume that $X^k=1$. We can now define a $\fU'\ot_\fU H_n^{D,1}$-module structure
on the vector space $\fU'\ot_\fU E_u$ which extends the $\fU'\ot_\fU H_n^1$-module 
structure and in which $\tT_{\uD}$ acts as $X$.

For any $u\in\cu^\fa$ we choose an $H_n^{D,1}$-module structure on $E_u$ extending the
$H_n^1$-module structure.

Let $u\in\cu$. We regard $E_u$ as a (simple) $\fU\ot_\ZZ H_n^\iy$-module $E_u^\iy$ via 
$\Ph^1:H_n^1@>\si>>\fU\ot_\ZZ H_n^\iy$. (If $u\in\cu^\fa$ we also regard $E_u$ as a 
$\fU\ot_\ZZ H^{D,\iy}_n$-module $E_u^\iy$ via 
$\Ph^1:H_n^{D,1}@>\si>>\fU\ot_\ZZ H^{D,\iy}_n$. This extends the 
$\fU\ot_\ZZ H_n^\iy$-module structure.)

Now $\fU[v,v\i]\ot_{\fU}E_u^\iy$ is naturally a $\fU[v,v\i]\ot_\ZZ H_n^\iy$-module and 
also an $H_n$-module via the homomorphism 
$H_n@>\Ph>>\ca\ot_\ZZ H_n^\iy\sub\fU[v,v\i]\ot_\ZZ H_n^\iy$. This $H_n$-module is 
denoted by $E_u(v)$. (If $u\in\cu^\fa$ then $\fU[v,v\i]\ot_{\fU}E_u^\iy$ is naturally a
$\fU[v,v\i]\ot_\ZZ H^{D,\iy}_n$-module and also an $H^D_n$-module via the homomorphism 
$H^D_n@>\Ph>>\ca\ot_\ZZ H^{D,\iy}_n\sub\fU[v,v\i]\ot_\ZZ H^{D,\iy}_n$. This extends the
$H_n$-module structure on $E_u(v)$.) 

Let $E_u^v=\fU(v)\ot_{\fU}E_u(v)$. From 34.12(c) we see that $\{E_u^v;u\in\cu\}$ is a 
set of representatives for the isomorphism classes of simple $H_n^v$-modules. (If 
$u\in\cu^\fa$ then the $H_n^{D,v}$-module structure on $E_u^v$ coming from the 
$H^D_n$-module structure on $E_u(v)$ extends the $H_n^v$-module structure.)

For $\k$ as in 34.12(b) let $E_u^\k$ be the vector space $E_u^\iy$ regarded as an 
$H_u^\k$-module via $\Ph^\k:H_n^\k@>\si>>\fU\ot_\ZZ H_n^\iy$. From 34.12(b) we see that
$\{E_u^\k;u\in\cu\}$ is a set of representatives for the isomorphism classes of simple 
$H_n^\k$-modules. Now $E_u^\k$ can also be obtained from the $H_n$-module $E_u(v)$ 
under the specialization $\fU[v,v\i]@>>>\fU,v\m\k$. Moreover, we have $E_u^1=E_u$ as 
$H_n^1$-modules. (If $u\in\cu^\fa$ we also regard $E_u^\iy$ as an $H_u^{D,\k}$-module 
via $\Ph^\k:H_n^{D,\k}@>\si>>\fU\ot_\ZZ H^{D,\iy}_n$. This extends the $H_u^\k$-module 
structure. This can be also obtained from the $H^D_n$-module $E_u(v)$ under the 
specialization $\fU[v,v\i]@>>>\fU,v\m\k$. Moreover, $E_u^1=E_u$ as 
$H_u^{D,1}$-modules.)

From the definitions we see that the map $\cu@>>>\cu,u\m\bu$ defined as in 34.14, 
replacing $(\fC,\fa$) by $(H_n^\k,\fa)$ ($\k$ as in 34.12(b)) or by $(H_n^v,\fa)$ is 
the same as the map $u\m\bu$ defined in terms of $(H_n^1,\fa)$. We show:

(a) {\it Let $w\in\WW,\l\in\ufs_n,u\in\cu$. Then 
$\tr(\tT_w1_\l,E_u^v)\in\fU[v,v\i]$ and, for $\k$ as in 34.12(b), 
$\tr(\tT_w1_\l,E_u^\k)\in\fU$ is obtained from this element of $\fU[v,v\i]$ by the 
specialization $\fU[v,v\i]@>>>\fU,v\m\k$. If, in addition, $u\in\cu^\fa$ and $j\in\ZZ$
then $\tr(\tT_w1_\l\tT_{\uD}^j,E_u^v)\in\fU[v,v\i]$ and, for $\k$ as in 34.12(b), 
$\tr(\tT_w1_\l\tT_{\uD}^j,E_u^\k)\in\fU$ is obtained from this element of $\fU[v,v\i]$ 
by the specialization $\fU[v,v\i]@>>>\fU,v\m\k$.}
\nl
This follows immediately from the fact that $\Ph(\tT_w1_\l)$ (resp. 
$\Ph(\tT_w1_\l\tT_{\uD}^j)$) is an $\ca$-linear combination of the standard basis 
elements of $H_n^\iy$ (resp.$H^{D,\iy}_n$.)

Combining 34.13(a), 34.14(d) we see that, for $u,u'\in\cu$ and for $\k$ as in 34.12(b)
we have
$$\align&\sum_{w\in\WW,\l\in\ufs_n}\tr(\tT_w1_\l,E_u^v)\tr(1_\l\tT_{w\i},E_{u'}^v)=
\d_{u,u'}f_u^v\dim E_u,\\&\sum_{w\in\WW,\l\in\ufs_n}\tr(\tT_w1_\l,E_u^\k)
\tr(1_\l\tT_{w\i},E_{u'}^\k)=\d_{u,u'}f_u^\k\dim E_u,\tag b\endalign$$
where $f_u^v\in\fU(v)-\{0\},f_u^\k\in\fU-\{0\}$. Using (a) we see that 
$f_u^v\in\fU[v,v\i]$ and $f_u^\k$ is obtained from $f_u^v$ by the specialization 
$\fU[v,v\i]@>>>\fU,v\m\k$. Combining 34.13(a), 34.14(e) we see that, for 
$u,u'\in\cu^\fa$, we have
$$\sum_{w\in\WW,\l\in\ufs_n}\tr(\tT_w1_\l\tT_{\uD},E_u^v)
\tr(\tT_{\uD}\i1_\l\tT_{w\i},E_{u'}^v)=\d_{u,u'}f_u^v\dim E_u.\tag c$$

\subhead 34.16\endsubhead
Let $x\m x^\sp$ be the automorphism of the field $\fU$ which sends any root of $1$ to 
its inverse. We extend this to an automorphism of the field $\fU(v)$ (denoted by
$\x\m\x^\sp$) which carries $v$ to itself. For $x'\in\fU$ we say that $x'>0$ if the 
image of $x'$ under any imbedding of $\fU$ into the complex numbers is a real number 
$>0$. For example, for $x\in\fU-\{0\}$ we have $xx^\sp>0$. For $\x'\in\fU(v)$ we say 
that $\x'>0$ if $\x'$ can be expanded in a power series $\x'=a_0v^n+a_1v^{n+1}+\do$ 
where $a_0,a_1,\a_2,\do\in\fU$ and $a_0>0$.

\proclaim{Lemma 34.17} Let $u\in\cu^\fa,w\in\WW^D,\l\in\ufs_n$. We have
$$\tr(1_\l\tT_{w\i},E_u^v)=\tr(\tT_w1_\l,E_u^v)^\sp.\tag a$$
\endproclaim
The antiautomorphism $h@>>>h^\flat$ of $H_n$ (see 32.19) extends to an antiautomorphism
$h@>>>h^\flat$ of $H^D_n$ given by $\tT_{w'}\m\tT_{w'{}\i}$ for $w'\in\WW^D$, 
$1_{\l'}\m1_{\l'}$ for $\l'\in\ufs_n$. Define a ring involution $h\m h^{\di}$ of 
$H_n^{D,v}$ by $\sum_{w,\l}a_{w,\l}\tT_w1_\l\m\sum_{w,\l}a_{w,\l}^\sp(\tT_w1_\l)^\flat$
where $a_{w,\l}\in\fU(v)$. Assume that there exists a pairing 
$\la,\ra:E_u^v\T E_u^v@>>>\fU(v)$ such that $\la,\ra$ is linear in the second variable,
semi-linear (with respect to $\x\m\x^\sp$) in the first variable, is non-degenerate, 
satisfies $\la e,e'\ra=\la e',e\ra^\sp$ for $e,e'\in E_u^v$ and 
$$\la he,e'\ra=\la e,h^{\di}e'\ra\text{ for }e,e'\in E_u^v,h\in H_n^{D,v}.\tag b$$
If $(e_j),(e'_j)$ are bases of $E_u^v$ such that $\la e_i,e'_j\ra=\d_{ij}$ then

$\tr(h,E_u^v)=\sum_j\la e_j,he'_j\ra=\sum_j\la h^{\di}e_j,e'_j\ra=
\tr(h^{\di},E_u^v)^\sp$.
\nl
Taking here $h=1_\l\tT_{w\i}$ we see that (a) would follow. It remains to prove the 
existence of $\la,\ra$ as above.

We can find a pairing $\la,\ra':E_u^v\T E_u^v@>>>\fU(v)$ which is linear in the second 
variable, semi-linear (with respect to $\x\m\x^\sp$) in the first variable, satisfies 
$\la e,e'\ra'=\la e',e\ra'{}^\sp$ for $e,e'\in E_u^v$ and $\la e,e\ra>0$ for 
$e\in E_u^v-\{0\}$. (For example, we choose a basis $(e_j)$ of $E_u^v$ and we set
$\la\sum_ja_je_j,\sum_ja'_je_j\ra'=\sum_ja_j^\sp a'_j$ where $a_j,a'_j\in\fU(v)$.) We 
define a new pairing $\la,\ra:E_u^v\T E_u^v@>>>\fU(v)$ by
$$\la e,e'\ra=\sum_{w'\in\WW^D,\l'\in\ufs_n}\la\tT_{w'}1_{\l'}e,\tT_{w'}1_{\l'}e'\ra'.
$$
Note that $\la,\ra'$ is linear in the second variable, semi-linear (with respect to 
$\x\m\x^\sp$) in the first variable, satisfies $\la e,e'\ra=\la e',e\ra'{}^\sp$ for 
$e,e'\in E_u^v$ and $\la e,e\ra>0$ for $e\in E_u^v-\{0\}$. In particular, $\la,\ra$ is 
non-degenerate. We show that (b) holds. It is enough to show this when $h$ runs through
a set of generators of the algebra $H_n^{D,v}$ that is, for $h=1_\l$ or $h=\tT_{\uD}$ 
or $h=\tT_s(s\in\II)$. Assume first that $h=1_\l,\l\in\ufs_n$. We must show that

$\sum_{w'\in\WW^D,\l'\in\ufs_n}\la\tT_{w'}1_{\l'}1_\l e,\tT_{w'}1_{\l'}e'\ra'=
\sum_{w'\in\WW^D,\l'\in\ufs_n}\la\tT_{w'}1_{\l'}e,\tT_{w'}1_{\l'}1_\l e'\ra'$.
\nl
Both sides are equal to $\sum_{w'\in\WW^D}\la\tT_{w'}1_\l e,\tT_{w'}1_\l e'\ra'$.
Assume next that $h=\tT_{\uD}$. We must show that
$$\sum\Sb w'\in\WW^D\\ \l'\in\ufs_n\eSb
\la\tT_{w'}1_{\l'}\tT_{\uD}e,\tT_{w'}1_{\l'}e'\ra'=
\sum\Sb w'\in\WW^D\\ \l'\in\ufs_n\eSb
\la\tT_{w'}1_{\l'}e,\tT_{w'}1_{\l'}\tT_{\uD\i}e'\ra'$$
that is, $\sum\Sb w'\in\WW^D\\ \l'\in\ufs_n\eSb
\la\tT_{w'\uD}1_{\uD\i\l'}e,\tT_{w'}1_{\l'}e'\ra'=
\sum\Sb y\in\WW^D\\ \l''\in\ufs_n\eSb
\la\tT_y1_{\l''}e,\tT_{y\uD\i}1_{\uD\l''}e'\ra'$, \lb which is clear.
Finally, assume that $h=\tT_s,s\in\II$. We must show that 

$\sum_{w'\in\WW^D,\l'\in\ufs_n}\la\tT_{w'}1_{\l'}\tT_s e,\tT_{w'}1_{\l'}e'\ra'=
\sum_{w'\in\WW^D,\l'\in\ufs_n}\la\tT_{w'}1_{\l'}e,\tT_{w'}1_{\l'}\tT_s e'\ra'$.
\nl
Both sides are equal to
$$\sum\Sb w'\in\WW^D\\ \l'\in\ufs_n\eSb\la\tT_{w's}1_{s\l'}e,\tT_{w'}1_{\l'}e'\ra'
+(v-v\i)\sum\Sb w'\in\WW^D\\ l(w's)=l(w')-1\\ \l'\in\ufs_n,s\in\WW_{\l'}\eSb
\la\tT_{w'}1_{\l'}e,\tT_{w'}1_{\l'}e'\ra'.$$
This proves (b). The lemma is proved.

\subhead 34.18\endsubhead
Using 34.17(a) we can rewrite 34.15(c) for $u,u'\in\cu^\fa$ as follows:
$$\sum_{w\in\WW,\l\in\ufs_n}\tr(\tT_w1_{\uD\l}\tT_{\uD},E_u^v)
\tr(\tT_w1_{\uD\l}\tT_{\uD},E_{u'}^v)^\sp=\d_{u,u'}f_u^v\dim E_u.\tag a$$
Specializing this under $\fU[v,v\i]@>>>\fU,v\m\k$, we obtain
$$\sum_{w\in\WW,\l\in\ufs_n}\tr(\tT_w1_{\uD\l}\tT_{\uD},E_u^\k)
\tr(\tT_w1_{\uD\l}\tT_{\uD},E_{u'}^{\k^\sp})^\sp=\d_{u,u'}f_u^\k\dim E_u.\tag b$$

\subhead 34.19\endsubhead
Let $A$ be a character sheaf on $D$. Define an $\ca$-linear map $\hat\z:H_n@>>>\ca$ by 
$h\m\z^A(h[D])$ where $\z^A:H_n[D]@>>>\ca$ is as in 31.7. From 31.8 we see that 
$\hat\z(hh')=\hat\z(h'\fa(h))$ for $h,h'\in H_n$. Applying 34.14(c) to the linear
map $\fU(v)\ot_\ca H_n@>>>\fU(v)$ obtained from $\hat\z$ by extension of scalars, we 
see that there exist elements $b_{A,u}^v\in\fU(v)(u\in\cu^\fa)$ such that
$$\z^A(\tT_w1_{\uD\l}[D])
=\sum_{u'\in\cu^\fa}b_{A,u'}^v\tr(\tT_w1_{\uD\l}\tT_{\uD},E_{u'}^v),\tag a$$
for $w\in\WW,\l\in\ufs_n$. We multiply both sides of (a) by 
$\tr(\tT_w1_{\uD\l}\tT_{\uD},E_u^v)^\sp$ (with $u\in\cu^\fa$) and sum over all $w,\l$. 
Using 34.18(a), we obtain
$$b_{A,u}^v=\fra{1}{f_u^v\dim E_u}\sum_{w\in\WW,\l\in\ufs_n}\z^A(\tT_w1_{\uD\l}[D])
\tr(\tT_w1_{\uD\l}\tT_{\uD},E_u^v)^\sp.\tag b$$
Using 28.17(a),(b) and the notation there, we see that 
$\fD({}^pH^j(\bK^{\ss,\cl}_D))={}^pH^j(\bK^{\ss,\che\cl}_D)$ hence
$$\align&\sum_{u\in\cu^\fa}b_{\fD(A),u}^v\tr(C^\ss_{\uD\l}\tT_{\uD},E_u^v)=
\z^{\fD(A)}(C^\ss_{\uD\l}[D])\\&
=\sum_j(-v)^jv^{-\dim G}(\fD(A):{}^pH^j(\bK^{\ss,\cl}_D))
=\sum_j(-v)^jv^{-\dim G}(A:\fD({}^pH^j(\bK^{\ss,\cl}_D))\\&
=\sum_j(-v)^jv^{-\dim G}(A:{}^pH^j(\bK^{\ss,\che\cl}_D))=\z^A(C^\ss_{\uD\l\i}[D])\\&
=\sum_{u\in\cu^\fa}b_{A,u}^v\tr(C^\ss_{\uD\l\i}\tT_{\uD},E_u^v).\tag c\endalign$$

\proclaim{Lemma 34.20} Let $u\in\cu^\fa$. Assume that $E_u$ is quasi-rational in the
following sense: there exists a function 
$\et:\WW^D\T\ufs_n@>>>\{\text{roots of $1$ in }\fU\},(w,\l)\m\et_{w,\l}$ such that 
$\et$ is constant on any equivalence class for $\asi$ in $\WW^D\T\ufs_n$ (see 32.26) 
and $\tr(\tT_w1_\l,E_u)\in\et_{w,\l}\ZZ$ for all $(w,\l)\in\WW^D\T\ufs_n$. Then
$\tr(\tT_w1_\l,E_u^v)\in\et_{w,\l}\ca$ for all $(w,\l)\in\WW^D\T\ufs_n$. 
\endproclaim
From the definitions and 34.7-34.10 we see that, for $w\in\WW^D,\l\in\ufs_n$, the basis
elements $c_{w,\l}$ of $H^D_n$ (see 34.5) satisfy

(a) $c_{w,\l}\in\tT_{wy}1_\l+\sum_{y\in\WW_\l;wy<w}\ca\tT_{wy}1_\l$.
\nl
For $w\in\WW^D,\l\in\ufs_n,x\in\WW_\l$ we have 
$\tT_w1_\l\tT_x\in\sum_{x'\in\WW_\l}\ca\tT_{wx'}1_\l$ (in $H^D_n$). This follows by 
writing $x$ as product $x=\s_1\s_2\do\s_m,\s_m\in\II_\l,m=l_\l(x)$ and using 
repeatedly 34.7(a). Using this and (a) we see that

(b) $c_{w,\l}\tT_x\in\sum_{x'\in\WW_\l}\ca\tT_{wx'}1_\l$.
\nl
Now let $c_{w',\l}$ be a distinguished basis element of $H^D_n$ (see 34.1). By 34.10(e)
we have $w'\in\WW_\l$ hence $c_{w',\l}\in\sum_{x\in\WW_\l}1_\l\tT_x$. Hence from (b) we
deduce

(c) $c_{w,\l}c_{w',\l}\in\sum_{x'\in\WW_\l}\ca\tT_{wx'}1_\l$.
\nl
From (a) we deduce by inversion:

(d) $\tT_w1_\l\in\sum_{y\in\WW_\l;wy\le w}\ca c_{wy,\l}$.
\nl
Using this for $wx'$ instead of $x$ and using also (c) we deduce

$c_{w,\l}c_{w',\l}\in\sum_{x''\in\WW_\l}\ca c_{wx'',\l}$.
\nl
Hence, if $\Ph:H^D_n@>>>\ca\ot_\ZZ H^{D,\iy}_n$ is as in 34.12, then

(e) $\Ph(c_{w,\l})\in\sum_{x''\in\WW_\l}\ca t_{wx'',\l}$
\nl
where $t_{y,\l}$ be the basis element of $H^{D,\iy}_n$ corresponing to 
$c_{y,\l}\in H^D_n$ (see 34.1). 

Let $c_{w,\l;1}$ be the image of $c_{w,\l}$ in $H_n^{D,1}$. Let 
$\Ph':\QQ\ot_\ca H^D_n@>>>\QQ\ot_\ZZ H^{D,\iy}_n$ be the homomorphism obtained from 
$\Ph$ under the specialization $\ca@>>>\QQ,v\m 1$. This is an isomorphism of algebras 
since $\Ph^1$ (see 34.12) is an isomorphism. From (e) we see that for any 
$(y,\l)\in\WW^D\T\ufs_n$, $\Ph'$ restricts to a $\QQ$-linear map 

(f) $\sum_{x\in\WW_\l}\QQ c_{yx,\l;1}@>>>\sum_{x\in\WW_\l}\QQ t_{yx,\l}$.
\nl
The vector spaces in (f) form direct sum decompositions of
$\QQ\ot_\ca H^D_n,\QQ\ot_\ZZ H^{D,\iy}_n$ hence (f) must be an isomorphism. We deduce 
that $\Ph'{}\i$ carries $\sum_{x\in\WW_\l}\QQ t_{yx,\l}$ onto
$\sum_{x\in\WW_\l}\QQ c_{yx,\l;1}$. We see that $t_{y,\l}$ acts on the 
$\fU\ot_\ZZ H^{D,\iy}_n$-module $E_u^\iy$ as a $\QQ$-linear combination of the 
operators
$c_{yx,\l;1}:E_u@>>>E_u$ where $x\in\WW_\l$ hence also as a $\QQ$-linear combination of
the operators $\tT_{yx}1_\l:E_u@>>>E_u$. (From (a) specialized for $v=1$ we see that 
$c_{y,\l}\in\sum_{x\in\WW_\l}\ZZ\tT_{yx}1_\l$.) It follows that 
$\tr(t_{y,\l},E_u^\iy)\in\et_{y,\l}\QQ$. Using (e) we see that $c_{y,\l}$ acts on 
$E_u^v$ as an $\ca$-linear combination of the operators 
$1\ot t_{yx,\l}$ on $\fU(v)\ot_\fU E_u^\iy$ where $x\in\WW_\l$. The same holds for 
$\tT_{y,\l}$ (instead of $c_{y,\l}$), by (d). It follows that
$\tr(\tT_y1_\l,E_u^v)\in\et_{y,\l}\QQ[v,v\i]$.

Since the algebra $\fU\ot_\ZZ H^{D,\iy}_n$ is of finite dimension (say $m$), with $1$,
and its structure constants with respect to the basis $(t_{y,\l})$ are integers, we see
that any basis element $t_{y,\l}$ satisfies an equation of the form 
$t_{y,\l}^m+c_1t_{y,\l}^{m-1}+\do+c_m=0$ where $c_1,c_2,\do,c_m$ are integers. It 
follows that $\tr(t_{y,\l},E_u^\iy)$ is an algebraic integer (necessarily in $\fU$).
Since the definition of $\Ph$ involves only coefficients in $\ca$ it follows that the 
coefficients of $\tr(c_{y,w},E_u^v)\in\fU[v,v\i]$ are algebraic integers in $\fU$. The
same holds then for the coefficients of $\tr(\tT_y1_\l,E_u^v)$. An element of 
$\et_{y,\l}\QQ[v,v\i]$ whose coefficients are algebraic integers in $\fU$ is 
necessarily in $\et_{y,\l}\ca$. We see that $\tr(\tT_y1_\l,E_u^v)\in\et_{y,\l}\ca$, as 
required.

\proclaim{Lemma 34.21} Let $u,E_u,\et$ be as in 34.20. Let $A$ be a character sheaf on 
$D$. Let $\fE_A\sub\WW^D\T\ufs_n$ be the equivalence class under $\asi$ attached to $A$
in 32.25(a) (with $J=\II$). Let $\et_0$, a root of $1$, be the (constant) value of $\et$ 
on $\fE_A$. We have $b_{A,u}^v\in\et_0\i\QQ(v)$.
\endproclaim
From 34.18(a) we have
$$f_u^v=(\dim E_u)\i\sum_{w\in\WW,\l\in\ufs_n}\tr(\tT_{w\uD}1_\l,E_u^v)
\tr(\tT_{w\uD}1_\l,E_u^v)^\sp.$$
By 34.20, for any $w\in\WW,\l\in\ufs_n$ we have
$\tr(\tT_{w\uD}1_\l,E_u^v)=\et_{w\uD,\l}Q_{w,\l}$ where $Q_{w,\l}\in\ca$; hence
$\tr(\tT_{w\uD}1_\l,E_u^v)\tr(\tT_{w\uD}1_\l,E_u^v)^\sp=
\et_{w\uD,\l}Q_{w,\l}\et_{w\uD,\l}\i Q_{w,\l}$. Thus,
$$f_u^v=(\dim E_u)\i\sum_{w\in\WW,\l\in\ufs_n}Q_{w,\l}^2\in\QQ[v,v\i].\tag a$$
Using 32.26(b) we rewrite 34.19 as follows
$$b_{A,u}^v=\fra{1}{f_u^v\dim E_u}\sum\Sb x\in\WW,\l\in\ufs_n\\(x\uD,\l)\in\fE_A\eSb
\z^A(\tT_x1_{\uD\l}[D])\tr(\tT_{x\uD}1_\l,E_u^v)^\sp.$$  
For each $x,\l$ in the sum, we have $\z^A(\tT_x1_{\uD\l}[D])\in\ca$ (by definition) and

$\tr(\tT_{x\uD}1_\l,E_u^v)^\sp\in\et_0\i\ca$
\nl
(by 34.20). Using this and (a), we see that $b_{A,u}^v\in\et_0\i\QQ(v)$.

\head 35. Functions on $G^{0F}/U$\endhead
\subhead 35.1\endsubhead
In this section we assume that $\kk$ is an algebraic closure of a finite field $\FF_q$ 
of characteristic $p$ and that $G$ has a fixed $\FF_q$-rational structure whose 
Frobenius map induces the identity map on $\WW$ and on $G/G^0$ and the map $t\m t^q$ on
$\TT$. For any algebraic variety $V$ over $\kk$ with a given $\FF_q$-structure we 
denote by $F:V@>>>V$ the corresponding Frobenius map. We fix an integer $n\ge 1$ that 
divides $q-1$.

\subhead 35.2\endsubhead
In this section we fix an {\it \'epinglage} of $G^0$ compatible with the 
$\FF_q$-structure. Thus, we fix $B^*,T,U^*$ as in 28.5 such that $F(B^*)=B^*,F(T)=T$ 
and we fix for each $s\in\II$ an isomorphism $a\m\x_s(a)$ of $\kk$ onto a subgroup of 
$U^*$ such that $t\x_s(a)t\i=\x_s(\a_s(t)a)$ for all $t\in T,a\in\kk$ and 
$F(\x_s(a))=\x_s(a^q)$ for all $a\in\kk$; here $\a_s\in R^+$ and the corresponding 
coroot $\cha_s$ satisfy $t=s(t)\cha_s(\a_s(t))$ for all $t\in T$. (Clearly, such an
\'epinglage exists and any two such \'epinglages are conjugate under the action of 
$(G^0/\cz_{G^0})^F$.) We identify $T=\TT$ as in 28.5. 
For $s\in\II$ let ${}'U^*_s$ be the root subgroup of $G^0$ 
corresponding to the root $\a_s\i$. Define $y\in {}'U^*_s-\{1\}$, $\ds\in N_{G^0}T$ by
$\ds=\x_s(1)y\x_s(1)$; then $\ds$ is a representative of $s$ in $(N_{G^0}T)^F$. For 
$s\in\II\cup\{1\}$ we define $\cha_s:\kk^*@>>>T$ as above if 
$s\in\II$ and to be $1$, if $s=1$.

Following Tits, we can define uniquely for each $w\in\WW$ a representative $\dw$ in 
$N_{G^0}T$ by the requirements:

(i) if $s\in\II$ then $\ds$ is as above;

(ii) if $w,w'\in\WW$ and $l(ww')=l(w)+l(w')$ then $(ww')\dot{}=\dw\dw'$;

(iii) $\dot 1=1$.
\nl
We have $F(\dw)=\dw$ for any $w\in\WW$. Let $\che T^F=\Hom(T^F,\bbq^*)$. If 
$\cl\in\fs_{q-1}$ (see 31.2) then 
$\cl^{\ot(q-1)}\cong\cl$ hence $F^*\cl\cong\cl$ and there is a unique isomorphism 
$\t_0:F^*\cl@>\si>>\cl$ that induces the identity on the stalk at $1$. Hence we can 
form the characteristic function $\c_{\cl,\t_0}:T^F@>>>\bbq^*$ (a group homomorphism). 
Let $\l\in\ufs_{q-1}$ be the isomorphism class of $\cl$. We set 
$\th^\l_F=\c_{\cl,\t_0}\in\che T^F$. Now $\l\m\th^\l_F$ is a bijection 
$\ufs_{q-1}@>\si>>\che T^F$.

If $\th\in\che T^F,\a\in R$ we write $\th\cha$ for the composition
$\FF_q^*@>\cha|_{\FF_q^*}>>T^F@>\th>>\bbq^*$. For $\a\in R,b\in\Wb$ we write $b\cha$ 
for the coroot $a\m b(\cha(a))$.

\subhead 35.3\endsubhead
In this section we write $U$ instead of $U^{*F}$. Let $\fU$ be as in 34.12. Let $\fT$ 
be the vector space of all functions $G^{0F}@>>>\fU$ that are constant on $U,U$ double
cosets. Now $\fT$ has a basis $\{k_\nu;\nu\in(N_{G^0}T)^F\}$ where $k_\nu$ is $1$ on 
$U\nu U$ and is $0$ on $G^{0F}-U\nu U$. We regard $\fT$ as an associative $\fU$-algebra
in which the product of $h_1,h_2$ is given by
$$(h_1h_2)(g)=|U|\i\sum_{g_1,g_2\in G^{0F};g_1g_2=g}h_1(g_1)h_2(g_2).$$
This algebra has $1=k_1$. As in \cite{\YO}, we have 
$$k_{\ds}k_{\ds}=qk_{\cha_s(-1)}+\sum_{a\in\FF_q^*}k_{\ds}k_{\cha_s(a)},$$
$$k_\nu k_{\nu'}=k_{\nu\nu'},$$
where $s\in\II$ and $\nu,\nu'$ represent $w,w'$ in $\WW$ such that $l(ww')=l(w)+l(w')$.
For any $\l\in\ufs_{q-1}$ we set
$$1_\l=|T^F|\i\sum_{t\in T^F}\th^\l_F(t)k_t\in\fT.$$
Then $1_\l1_{\l'}=\d_{\l,\l'}1_\l$ for $\l,\l'\in\ufs_{q-1}$. Let 
$\sqrt{-1}$ be a fixed element of $\fU^*$ whose square is $-1$; we set $\sqrt{1}=1$. 
For $s\in\II$ we set 
$$T_s=k_{\ds}\sum_{\l\in\ufs_{q-1}}\sqrt{\th^\l_F(\cha_s(-1))}1_\l\in\fT.$$
More generally, for $w\in\WW$ we set
$$T_w=\sum_{\l\in\ufs_{q-1}}\prod_{\a\in R^+,w\i\a\in R^-}\sqrt{\th^\l_F(\cha(-1))}
k_{\dw}1_\l\in\fT.$$
From the definitions we see that $T_wT_{w'}=T_{ww'}$ for $w,w'\in\WW$ with 
$l(ww')=l(w)+l(w')$ and $T_w1_\l=1_{w\l}T_w$. We have
$$\align&T_sT_s=\sum_\l\th^\l_F(\cha_s(-1))k_{\ds}k_{\ds}1_\l\\&
=q\sum_\l\th^\l_F(\cha_s(-1))k_{\cha_s(-1)}1_\l+\sum_\l\th^\l_F(\cha_s(-1))
\sum_{a\in\FF_q^*}k_{\ds}k_{\cha_s(a)}1_\l\\&=q+
\sum_\l\th^\l_F(\cha_s(-1))\sum_{a\in\FF_q^*}\th^\l_F(\cha_s(a))k_{\ds}1_\l\\&
=q+(q-1)\sum\Sb\l\\ \th^\l_F\cha_s=1\eSb\th^\l_F(\cha_s(-1))k_{\ds}1_\l=
q+(q-1)\sum\Sb\l\\ \th^\l_F\cha_s=1\eSb\sqrt{\th^\l_F(\cha_s(-1))}k_{\ds}1_\l.\endalign
$$
Thus,
$$T_sT_s=q+(q-1)\sum_{\l\in\ufs_{q-1};\th^\l_F\cha_s=1}T_s1_\l.$$

\subhead 35.4\endsubhead
We fix a square root $\sqrt{p}$ of $p$ in $\fU$. For any $e\in\ZZ$ we set
$\sqrt{p^e}=(\sqrt{p})^e$. In particular, $\sqq\in\fU$ is defined. Then 
$H_{q-1}^{\sqq}$ is defined as in 34.12. From 35.3 we see that the elements $T_w,1_\l$ 
of $\fT$ define a $\fU$-algebra homomorphism $H_{q-1}^{\sqq}@>>>\fT$. (We use the fact 
that if $\l\in\ufs_{q-1}$ and $s\in\II$ then $\th^\l_F\cha_s=1$ if and only if 
$s\in\WW_\l$.) This is an isomorphism: from the definitions we see that 
$\{T_w1_\l;w\in\WW,\l\in\ufs_{q-1}\}$ is a $\fU$-basis of $\fT$; we then use 31.2(a).

\subhead 35.5\endsubhead
Let $\fP$ be the vector space of all functions $f:G^{0F}@>>>\fU$ such that $f(xu)=f(x)$
for all $x\in G^{0F},u\in U$. For $g\in G^F$, $g'\in(N_GB^*\cap N_GT)^F$ such that 
$gG^0=g'G^0$ we define a linear map $\r_{g,g'}:\fP@>>>\fP$ by 
$(\r_{g,g'}f)(x)=f(g\i xg')$. Then $g_0:f\m\r_{g_0,1}f$ makes $\fP$ into a 
$G^0F$-module.

Any element $c\in\fT$ defines a linear map $\fP@>>>\fP$:
$$f\m cf,(cf)(x)=|U|\i\sum_{x'\in G^{0F}}c(x')f(xx').$$
Clearly, $c\m[f\m cf]$ is an isomorphism $\fT@>\si>>\End_{G^{0F}}(\fP)$ and a left
$\fT$-module structure on $\fP$. For $\nu\in(N_{G^0}T)^F,\l\in\ufs_{q-1},f\in\fP$ we 
have
$$(k_\nu f)(x)=|U|\i\sum_{x'\in U\nu U}f(xx'),$$
$$(1_\l f)(x)=|T^F|\i\sum_{t\in T^F}\th^\l_F(t)f(xt).$$
In the remainder of this section we fix a connected component $D$ of $G$ and an element
$d\in(N_DB^*\cap N_DT)^F$. Let $s\in\II\cup\{1\},\l\in\ufs_{q-1},f\in\fP$. For 
$x\in G^{0F}$ we have:
$$(T_s1_{\uD\l}f)(x)=\sqrt{\th^{\uD\l}_F(\cha_s(-1))}|B^{*F}|\i\sum_{x'\in U\ds U,
t\in T^F}\th^{\uD\l}_F(t)f(xx't).\tag a$$
Now let $\ss=(s_1,s_2,\do,s_r)$ be a sequence in $\II\cup\{1\},\l\in\ufs_{q-1}$,
$$T_\ss1_{\uD\l}=T_{s_1}T_{s_2}\do T_{s_r}1_{\uD\l}=
T_{s_1}1_{s_2\do s_r\uD\l}T_{s_2}1_{s_3\do s_r\uD\l}\do T_{s_r}1_{\uD\l}.$$
Applying (a) $r$ times gives (for $f\in\fP,x\in G^{0F},g\in D^F$):
$$\align&(T_\ss1_{\uD\l}\r_{g,d}f)(x)=a_{\uD\l,F,\ss}|B^{*F}|^{-r}\\&
\T\sum\Sb g_1,g_2,\do,g_r\\t_1,t_2,\do t_r\\g_i\in U\ds_iU\\t_i\in T^F\eSb
\th^{s_2\do s_r\uD\l}_F(t_1)
\th^{s_3\do s_r\uD\l}_F(t_2)\do\th^{\uD\l}_F(t_r)f(g\i xg_1t_1g_2t_2\do g_rt_rd)\\&
=a_{\uD\l,F,\ss}|B^{*F}|^{-r}\\&
\T\sum\Sb g_1,g_2,\do,g_r\\t_1,t_2,\do t_r\\g_i\in U\ds_iU\\
t_i\in T^F\eSb\th^{\uD\l}_F((s_r\do s_2t_1)(s_r\do s_3t_2)\do(t_r))
f(g\i xg_1t_1g_2t_2\do g_rt_rd)\tag a${}'$\endalign$$
where 
$$\align&a_{\uD\l,F,\ss}\\&=\sqrt{\th^{\uD\l}_F(s_r\do s_2(\cha_{s_1}(-1)))}
\sqrt{\th^{\uD\l}_F(s_r\do s_3(\cha_{s_2}(-1)))}\do
\sqrt{\th^{\uD\l}_F(\cha_{s_r}(-1))}.\endalign$$
Let $J\sub\II$ and let $Q\in\cp_J$ be such that $F(Q)=Q$. Define a linear map 
$\pr_Q:\fP@>>>\fP$ by $(\pr_Qf)(x)=f(x)$ if $x\in G^{0F},xQ_{J,B^*}x\i=Q$ and 
$(\pr_Qf)(x)=f(x)$ if $x\in G^{0F},xQ_{J,B^*}x\i\ne Q$. We compute the trace of the 
linear map 
$$T_\ss1_{\uD\l}\r_{g,d}\pr_Q:\fP@>>>\fP\tag b$$
using the $\fU$-basis of $\fP$ consisting of the characteristic functions of the 
various right $U$-cosets in $G^{0F}$. Using (a${}'$) and the definitions we see that 
this trace equals
$$\fra{a_{\uD\l,F,\ss}}{|B^{*F}|^r|U|}\sum\Sb g_1,g_2,\do,g_r,\\t_1,t_2,\do 
t_r;\\g_i\in U\ds_iU,t_i\in T^F,x\in D^F;\\g\i xg_1t_1g_2t_2\do g_rt_rd\in xU;\\xQ_{J,
B^*}x\i=Q\eSb\th^\l_F(d\i(s_r\do s_2t_1)(s_r\do s_3t_2)\do(t_r)d).\tag c$$

\subhead 35.6\endsubhead
Let $\cl,\l,\t_0:F^*\cl@>\si>>\cl$ be as in 35.2. Let $\ss=(s_1,s_2,\do,s_r)$ be a 
sequence in $\II\cup\{1\}$ such that $s_1s_2\do s_r\uD\l=\l$. In 28.7 we have defined a
local system $\tcl$ on $Z^\ss_{\em,J,D}$ in terms of $d$ and a representative for 
$s_1s_2\do s_r$ in $N_{G^0}T$. We now take as a representative for $s_1s_2\do s_r$ the 
element $\ds_1\ds_2\do\ds_r$ with $\ds_i$ as in 35.2. We reformulate the definition of 
$\tcl$ as follows (see also the proof of 28.10). Define 
$\g:Z'{}^\ss@>>>Z^\ss_{\em,J,D}$ by the formula 28.10(a). Define $\ps:Z'{}^\ss@>>>T$ by
$$(h_0,h_1,\do,h_r,g)\m d\i(\ds_1\ds_2\do\ds_r)\i n_1n_2\do n_rn_0$$
with $n_i\in N_{G^0}T$ given by $h_{i-1}\i h_i\in U^*n_iU^*$ and 
$n_0\in N_GB^*\cap N_GT$ given by 
$h_r\i gh_0\in U^*n_0$. Then $\tcl$ is the local system on $Z^\ss_{\em,J,D}$ such that 
$\g^*\tcl=\ps^*\cl$. Note that $\g,\ps$ are naturally defined over $\FF_q$. Let 
$\ti\t:F^*\ps^*\cl@>\si>>\ps^*\cl$ be the isomorphism induced by $\t_0$. There is a 
well defined isomorphism $\t:F^*\tcl@>\si>>\tcl$ such that $\t$ induces via $\g$ the 
isomorphism $\ti\t$.

Let $\p_\ss:Z^\ss_{\em,J,D}@>>>Z_{J,D}$ be as in 28.12(a) and let 
$K=K^{\ss,\cl}_{J,D}=\p_{\ss!}\tcl\in\cd(Z_{J,D})$. Now $\p_\ss$ is naturally defined 
over $\FF_q$. Hence $\t:F^*\tcl@>\si>>\tcl$ induces an isomorphism $\o:F^*K@>\si>>K$ 
and $\c_{K,\o}:Z_{J,D}^F@>>>\bbq$ is defined. Let $\x=(Q,Q',gU_Q)\in Z_{J,D}^F$ (we can
take $g\in D^F$). Using the definitions and the Grothendieck trace formula we have
$$\c_{K,\o}(\x)=\sum_i(-1)^i\tr(\t^*,H^i_c(\p_\ss\i(\x),\tcl))=
\sum_{\et\in\p_\ss\i(\x)^F}\tr(\t,\tcl_\et)$$
where $\tcl_\et$ is the stalk of $\tcl$ at $\et$. Now $\g$ induces a map 
$Z'{}^{\ss F}@>>>Z^\ss_{\em,J,D}{}^F$ all of whose fibres have cardinal 
$|B^{*F}|^{r+1}|U_{J,B^*}^F|$. It follows that
$$|B^{*F}|^{r+1}|U_{J,B^*}^F|\c_{K,\o}(\x)=
\sum\Sb\ti\et\in Z'{}^{\ss F}\\\p_\ss\g(\ti\et)=\x\eSb\tr(\ti\t,(\ps^*\cl)_{\ti\et})=
\sum\Sb\ti\et\in Z'{}^{\ss F}\\\p_\ss\g(\ti\et)=\x\eSb\tr(\t_0,\cl_{\ps(\ti\et)}).$$
Let 
$$\align\Xi=\{(h_0,h_1&,\do,h_r)\in(G^{0F})^{r+1};h_{i-1}\i h_i\in B^*\ds_iB^*
(i\in[1,r]),\\&h_r\i gh_0\in N_GB^*,h_0Q_{J,B^*}h_0\i=Q\}.\endalign$$
Then $\{\ti\et\in Z'{}^{\ss F};\p_\ss\g(\ti\et)=\x\}$ may be identified with 
$\Xi\T(gU_Q^F)$. Hence
$$\c_{K,\o}(\x)=|B^{*F}|^{-r-1}\sum_{(h_0,h_1,\do,h_r)\in\Xi}\c_{\cl,\t_0}
(\d(h_0,h_1,\do,h_r))$$ 
with $\d:\Xi@>>>T$ given by 
$(h_0,h_1,\do,h_r)\m d\i(\ds_1\ds_2\do\ds_r)\i n_1n_2\do n_rn_0$ ($n_i,n_0$ as above). 
For any $(h_0,h_1,\do,h_r)\in\Xi$ we define $g_i\in U\ds_iU$,
$t_i\in T^F(i\in[1,r])$ by $h_{i-1}\i h_i=g_it_i$. We also set $h_0=x$. Then $\Xi$ 
becomes
$$\align\{(x,g_1,&g_2,\do,g_r,t_1,t_2,\do,t_r);x\in G^{0F},g_i\in U\ds_i U,t_i\in T^F;
\\&gx=xg_1t_1g_2t_2\do g_rt_rtdu\text{ for some }t\in T^F,u\in U;xQ_{J,B^*}x\i=Q\}
\endalign$$ 
and $\d:\Xi@>>>T$ becomes 
$$\align&(x,g_1,g_2,\do,g_r,t_1,t_2,\do,t_r)\m d\i(\ds_1\ds_2\do\ds_r)\i ds_1t_1\ds_2
t_2\do\ds_rt_rtd\\&=d\i(s_r\do s_2t_1)(s_r\do s_3t_2)\do(s_rt_{r-1})(t_rt)d.\endalign$$
We make the change of variable $t_rt\m t_r,t\m t$. Then $t$ no longer appears 
explicitly; it only introduces a factor $|T^F|$. We see that
$$\align&\c_{K,\o}(\x)\\&=
\fra{|T^F|}{|B^{*F}|^{r+1}}\sum\Sb g_1,g_2,\do,g_r,\\t_1,t_2,\do 
t_r,\\g_i\in U\ds_iU,t_i\in T^F,x\in D^F,\\g\i xg_1t_1g_2t_2\do g_rt_rd\in xU\\
xQ_{J,B^*}x\i=Q\eSb\th^\l_F(d\i(s_r\do s_2t_1)(s_r\do s_3t_2)\do(t_r)d).\endalign$$
This is the same (up to the factor $a_{\uD\l,F,\ss}$) as the expression 35.5(c). Using 
the equality between 35.5(c) and the trace of 35.5(b), we see that
$$\c_{K,\o}(Q,gQg\i,gU_Q)=a_{\uD\l,F,\ss}\i\tr(T_\ss1_{\uD\l}\r_{g,d}\pr_Q,\fP).\tag a
$$

\subhead 35.7\endsubhead
In the setup of 35.6, let $\cj_\ss\sub\cj^0\sub[1,r]$ be as in 28.9. For any 
$\cj\sub\cj_\ss$ let $\ss_\cj$ be as in 28.9. We have $\ss_\em=\ss$. Define 
$\g_\cj:Z'{}^{\ss_\cj}@>>>Z^{\ss_j}_{\em,J,D}$ by the formula 28.10(a). Define 
$\ps_\cj:Z'{}^{\ss_\cj}@>>>T$ as in 28.10 (with $\ds$ as in 35.2). Let $\tcl_\cj$ be 
the local system on $Z^{\ss_\cj}_{\em,J,D}$ such that $\g_\cj^*\tcl_\cj=\ps_\cj^*\cl$. 
Let $\ti\t^\cj:F^*\ps_\cj^*\cl@>\si>>\ps_\cj^*\cl$ be the isomorphism induced by 
$\t_0$. There is a well defined isomorphism $\t^\cj:F^*\tcl_\cj@>\si>>\tcl_\cj$ such 
that $\t^\cj$ induces via $\g_\cj$ the isomorphism $\ti\t^\cj$. Note that for 
$\cj=\em$, $\tcl_\cj,\ti\t^\cj,\t^\cj$ reduce to $\tcl,\ti\t,\t$ of 35.6.

Let $\bZ^\ss_{\em,J,D}$ be as in 28.9. This is a smooth irreducible variety, naturally 
defined over $\FF_q$ and $Z^{\ss_\em}_{\em,J,D}$ is open dense in $\bZ^\ss_{\em,J,D}$. 
Hence $\bcl=IC(\bZ^\ss_{\em,J,D},\tcl_\em)$ is defined and the isomorphism 
$\t^\em:F^*\tcl_\em@>\si>>\tcl_\em$ of local systems on $Z^{\ss_\em}_{\em,J,D}$ extends
canonically to an isomorphism $\t^\sh:F^*\bcl@>\si>>\bcl$ (of constructible sheaves, 
see 28.10). Restricting this isomorphism to the subset $Z^{\ss_\cj}_{\em,J,D}$ of 
$\bZ^\ss_{\em,J,D}$ (with $\cj\sub\cj_\ss$) we obtain an isomorphism 
$\t^{\sh\cj}:F^*\ce^\cj@>\si>>\ce^\cj$ where $\ce^\cj=\bcl|_{Z^{\ss_\cj}_{\em,J,D}}$. 
From 28.10 we see that $\ce^\cj$ is a local system isomorphic to $\tcl_\cj$.

\proclaim{Lemma 35.8}The isomorphisms $\t^{\sh\cj},\t^\cj$ correspond to each other 
under \lb some/any isomorphism of local systems $\ce^\cj@>\si>>\tcl_\cj$ on 
$Z^{\ss_\cj}_{\em,J,D}$. 
\endproclaim
The proof is a refinement of that of Lemma 28.10. Note that $Z'{}^{\ss_\em}$ is an open
dense subvariety of the smooth irreducible variety $\tZ^\ss$ (as in 28.10). Hence
$IC(\tZ^\ss,\ps_\em^*\cl)$ is defined and the isomorphism 
$\ti\t^\em:F^*\ps_\em^*\cl@>\si>>\ps_\em^*\cl$ of local systems on $Z'{}^{\ss_\em}$ 
extends canonically to an isomorphism 

$\ti\t^\sh:F^*IC(\tZ^\ss,\ps_\em^*\cl)@>\si>>IC(\tZ^\ss,\ps_\em^*\cl)$
\nl
which may be identified with the isomorphism induced by $\t^\sh$ through the fibration 
$\tZ^\ss@>>>\bZ^\ss_{\em,J,D}$ (see 28.10(a)) whose fibres are smooth and connected. 
Let $\tce^\cj=IC(\tZ^\ss,\ps_\em^*\cl)|_{Z'{}^{\ss_\cj}}$ and let 
$\ti\t^{\sh\cj}:F^*\tce^\cj@>\si>>\tce^\cj$ be the isomorphism induced by $\ti\t^\sh$ 
by restriction. It suffices to prove the following statement: 

{\it The isomorphisms $\ti\t^{\sh\cj},\ti\t^\cj$ correspond to each other under 
some/any isomorphism of local systems $\tce^\cj@>\si>>\ps_{\cj}^*\cl$.}
\nl
Let $\cl'=(\uD\i)^*\cl\in\fs(\TT)=\fs(T)$. Let $\t'_0:F^*\cl'@>\si>>\cl'$ be the unique
isomorphism which induces the identity map on the stalk of $\cl'$ at $1$. Let 
${}'\bZ^\ss,{}'Z^{\ss_\cj}$ be as in 28.10. Define ${}'\ps_\cj:{}'Z^{\ss_\cj}@>>>T$ as 
in 28.10 (with $\ds$ as in 35.2). Then ${}'\ps_\cj$ is compatible with the natural 
$\FF_q$-structures on ${}'Z^{\ss_\cj},T$; hence 
$\t'_0:F^*\cl'@>\si>>\cl'$ induces an isomorphism of local systems 
${}'\t^\cj:F^*{}'\ps_\cj^*\cl'@>\si>>{}'\ps_\cj^*\cl'$. Let 
$${}'\t^\sh:F^*IC({}'\bZ^\ss,{}'\ps_\em^*\cl')@>\si>>IC({}'\bZ^\ss,{}'\ps_\em^*\cl')$$
be the isomorphism induced by ${}'\t^\em:F^*{}'\ps_\em^*\cl'@>\si>>{}'\ps_\em^*\cl'$. 
Let 
$${}'\ce^\cj=IC({}'\bZ^\ss,{}'\ps_\em^*\cl')|_{{}'Z^{\ss_\cj}}$$ 
and let ${}'\t^{\sh\cj}:F^*{}'\ce^\cj@>\si>>{}'\ce^\cj$ be the isomorphism induced by 
${}'\t^\sh$ by restriction. As in the proof of 28.10, it suffices to prove the 
following statement: 

(a) {\it The isomorphisms ${}'\t^{\sh\cj},{}'\t^\cj$ correspond to each other under 
some/any isomorphism of local systems ${}'\ce^\cj@>\si>>{}'\ps_\cj^*\cl'$.}
\nl
Assume that (a) is known in the case where $|\cj|=1$. We now consider a general 
$\cj\sub\cj_\ss$. We prove (a) by induction on $|\cj|$. If $\cj=\em$, (a) is obvious. 
Assume that $\cj\ne\em$. Let $j\in\cj$ be the largest number in $\cj$. Let 
$\cj'=\cj-\{j\}$. Let ${}'\bZ^{\ss_{\cj'}}$ be the variety analogous to ${}'\bZ^\ss$ 
when $\ss$ is replaced by $\ss_{\cj'}$ (this is the same as the closure of 
${}'Z^{\ss_{\cj'}}$ in ${}'\bZ^\ss$). Let 
$${}''\t^\sh:F^*IC({}'\bZ^{\ss_{\cj'}},{}'\ps_{\cj'}^*\cl')@>\si>>
                                          IC({}'\bZ^{\ss_{\cj'}},{}'\ps_{\cj'}^*\cl')$$
be the isomorphism induced by 
${}'\t^{\cj'}:F^*{}'\ps_{\cj'}^*\cl'@>\si>>{}'\ps_{\cj'}^*\cl'$ and let 
${}''\t^{\sh\cj}$ be its restriction to ${}'Z^{\ss_\cj}$. By the induction hypothesis, 
the isomorphisms ${}'\t^{\sh\cj'},{}'\t^{\cj'}$ correspond to each other under some/any
isomorphism of local systems ${}'\ce^{\cj'}@>\si>>{}'\ps_{\cj'}^*\cl'$. It follows that
the isomorphisms ${}''\t^{\sh\cj},{}'\t^{\sh\cj}$ correspond to each other under 
some/any isomorphism of local systems
$$IC({}'\bZ^{\ss_{\cj'}},{}'\ps_{\cj'}^*\cl')|_{{}'Z^{\ss_\cj}}@>\si>>{}'\ce^\cj.$$
By our assumption (applied to $\ss_{\cj'}$ instead of $\ss$), the isomophisms
${}''\t^{\sh\cj},{}'\t^\cj$ correspond to each other under some/any isomorphism of 
local systems
$$IC({}'\bZ^{\ss_{\cj'}},{}'\ps_{\cj'}^*\cl')|_{{}'Z^{\ss_\cj}}@>\si>>{}'\ps_\cj^*\cl'.
$$
It follows that the isomorphisms ${}'\t^{\sh\cj},{}'\t^\cj$ correspond to each other 
under \lb some/any isomorphism of local systems ${}'\ce^\cj@>\si>>{}'\ps_\cj^*\cl'$. 
Thus, (a) holds for $\cj$.

We now consider the remaining case, where $\cj$ consists of a single element $j$. Note 
that $j\in\cj_\ss$ where $\cj_\ss$ is defined in terms of $D,\cl$ or equivalently, in 
terms of $G^0,\cl'$. The statement (a) involves only $G^0$. Hence to prove it, we may 
assume that $G=G^0$. We write $Z',Z''$ instead of ${}'Z^{\ss_\em},{}'Z^{\ss_{\{j\}}}$ 
and we set $Z=Z'\cup Z''$, a subvariety of ${}'\bZ^\ss$. We write $f',f''$ instead of 
${}'\ps_\em,{}'\ps_{\{j\}}$. Let $b=\ds_{j+1}\ds_{j+2}\do\ds_r$.

We have $\cl'=\k^*\cl_1$ where $\k\in\Hom(T,\kk^*)$ and $\cl_1\in\fs(\kk^*)$. Since 
$\cl'{}^{\ot(q-1)}\cong\bbq$, we may assume that $\cl_1^{\ot(q-1)}\cong\bbq$. Hence 
there is a unique isomorphism $\t_1:F^*\cl_1@>\si>>\cl_1$ which induces the identity on
the stalk of $\cl_1$ at $1$. Then $\t'_0:F^*\cl'@>\si>>\cl'$ is induced by $\t_1$, via 
$\k^*$. 

We continue the proof assuming that $G$ has simply connected derived subgroup. Let 
$\che\b$ be as in the proof of 28.10. As in that proof, we may assume that 
$\la\che\b,\k\ra=0$. Hence there exists a homomorphism of algebraic groups 
$\c:B^*\ds_jB^*\cup B^*@>>>\kk^*$ such that $\c(t)=\k(b\i tb)$ for all $t\in T$. Let 
$\tf:Z@>>>\kk^*$ be as in 28.10. If $y_j\in B^*\ds_jB^*$, we have
$\tf(y_1,\do,y_r)=\k(f'(y_1,\do,y_r))$; if $y_j\in B^*$ we have 
$\tf(y_1,\do,y_r)=\k(f''(y_1,\do,y_r))$. (See 28.10.) We show that 
$$\c(F(g))=\c(g)^q\text{ for all }g\in B^*\ds_jB^*\cup B^*.$$
We may assume that $g\in T$. Then 
$$\c(F(g))=\k(b\i F(g)b)=\k(F(b\i gb))=\k((b\i gb)^q)=\k(b\i gb)^q=\c(g)^q,$$
as required. It follows that for any $(y_1,\do,y_r)\in Z$ we have
$$\tf(F(y_1),\do,F(y_r))=(\tf(y_1,\do,y_r))^q.\tag b$$
Let $\cf=\tf^*(\cl_1)$. We have canonically 
$\cf_{Z'}=f'{}^*\cl',\cf|_{Z''}=f''{}^*\cl'$. From (b) we see that the isomorphism 
$\t_1:F^*\cl_1@>\si>>\cl_1$ induces via $\tf^*$ an isomorphism $F^*\cf@>\si>>\cf$ and 
this gives rise upon restriction to $Z',Z''$ to the isomorphisms 
${}'\t^\em,{}'\t^{\{j\}}$. Since $\cf$ is a local system on $Z$ such that 
$\cf|_{Z'}=f'{}^*\cl'$, we have canonically $\cf=IC(Z,f'{}^*\cl')$. Hence (a) holds in 
our case.

We now drop the assumption that $G$ has simply connected derived subgroup. Let 
$\p:\hG@>>>G$ be a surjective homomorphism of connected reductive groups whose kernel 
is a central torus in $\hG$ and such that $\hG$ has simply connected derived subgroup. 
We may assume that $\hG$ and $\p$ are defined over $\FF_q$. Then $\p$ restricts to a 
surjective homomorphism $\hG^F@>>>G^F$. Since the set of \'epinglages of $G,\hG$ are in
natural bijection, the \'epinglage of $G$ fixed in 35.2 gives rise to an \'epinglage of
$\hG$ (the associated Borel subgroup and maximal torus are $\hB^*=\p\i(B^*)$,
$\hT=\p\i(T)$ and the analogue of $\x_s:\kk@>>>U^*$ is the obvious one). For $s\in\II$ 
let $\ts\in(N_{\hG}\hT)^F$ be associate to this \'epinglage of $\hG$ in the same way as
$\ds$ was associated to the \'epinglage of $G$. We set $\ti 1=1\in\tG$. Define 
$\hZ,\hZ',\hZ'',\hf':\hZ'@>>>\hT,\hf'':\hZ''@>>>\hT$ in terms of $\hG,\hB^*,\hT,\ts_i$ 
in the same way as $Z,Z',Z'',f':Z'@>>>T,f'':Z''@>>>T$ are defined in terms of 
$G,B^*,T,\ds_i$. Let $\hcl'\in\fs(\hT)$ be the inverse image of $\cl'$ under 
$\p:\hT@>>>T$. There is a unique isomorphism $F^*\hcl'\cong\hcl'$ which induces the 
identity map on the stalk of $\hcl'$ at $1$. This induces isomorphisms
$${}'\htt^\em:F^*\hf'{}^*\hcl'@>>>\hf'{}^*\hcl',
{}'\htt^{\{j\}}:F^*\hf''{}^*\hcl'@>>>\hf''{}^*\hcl'.$$
Now ${}'\htt^\em$ induces an isomorphism 
$F^*IC(\hZ,\hf'{}^*\hcl')@>\si>>IC(\hZ,\hf'{}^*\hcl')$ and this restricts to an 
isomorphism of local systems 
$${}'\htt^{\sh\{j\}}:F^*IC(\hZ,\hf'{}^*\hcl')|_{\hZ''}@>\si>>
IC(\hZ,\hf'{}^*\hcl')|_{\hZ''}.$$
By an earlier part of the proof, the isomorphisms ${}'\htt^{\sh\{j\}},{}'\htt^{\{j\}}$ 
correspond to each other under some/any isomorphism of local systems 
$IC(\hZ,\hf'{}^*\hcl')|_{\hZ''}@>\si>>\hf''{}^*\hcl'$. Now the map $\hZ@>>>Z$ induced 
by $\p$ is a fibration with smooth connected fibres and $IC(\hZ,\hf'{}^*\hcl')$ is 
canonically the inverse image under this map of $IC(Z,f'{}^*\cl')$. Hence 
$IC(\hZ,\hf'{}^*\hcl')|_{\hZ''}$ is the inverse image under $\hZ''@>>>Z''$ of 
$IC(Z,f'{}^*\cl')|_{Z''}$. Similarly, $\hf''{}^*\hcl'$ is the inverse image under 
$\hZ''@>>>Z''$ of $f''{}^*\cl'$. Also, ${}'\htt^{\sh\{j\}},{}'\htt^{\{j\}}$ are 
obtained from ${}'\t^{\sh\{j\}},{}'\t^{\{j\}}$ by inverse image under $\hZ''@>>>Z''$. 
Therefore the required statement about ${}'\t^{\sh\{j\}},{}'\t^{\{j\}}$ is a 
consequence of the analogous, already known statement about 
${}'\htt^{\sh\{j\}},{}'\htt^{\{j\}}$ (we use the faithfulness of the inverse image
functor under the fibration $\hZ''@>>>Z''$ with smooth connected fibres). The lemma is 
proved.

\subhead 35.9\endsubhead
We preserve the setup of 35.7. Let $\bpi_\ss:\bZ^\ss_{\em,J,D}@>>>Z_{J,D}$ be as in 
28.12 and let $\bK=\bK^{\ss,\cl}_{J,D}=\bpi_{\ss!}\bcl\in\cd(Z_{J,D})$ (see 28.12, 35.7).
Now $\bpi_\ss$ is naturally defined over $\FF_q$. Hence the isomorphism 
$\t^\sh:F^*\bcl@>\si>>\bcl$ in 35.7 induces an isomorphism $\bao:F^*\bK@>>>\bK$ and 
$\c_{\bK,\bao}:Z_{J,D}^F@>>>\bbq$ is defined. 

\proclaim{Proposition 35.10} Let $(Q,Q',gU_Q)\in Z_{J,D}^F$ (we take $g\in D^F$). Let 
$C^\ss_{\uD\l}\in H_{q-1}$ be as in 31.5. We have
$$\c_{\bK,\bao}(Q,Q',gU_Q)=a_{\uD\l,F,\ss}\i\tr(C^\ss_{\uD\l}\r_{g,d}\pr_Q,\fP).$$
\endproclaim
Consider the partition $\bZ^\ss_{\em,J,D}=\sqc_{\cj\sub\cj^0}Z^{\ss_\cj}_{\em,J,D}$ 
(see 28.9). For each $\cj\sub\cj^0$ let $\p_\cj:Z^{\ss_\cj}_{\em,J,D}@>>>Z_{J,D}$ be 
the restriction of $\bp_\ss$, let $K_\cj=\p_{\cj!}(\bcl|_{Z^{\ss_\cj}_{\em,J,D}})$ and 
let $\o_\cj:F^*K_\cj@>>>K_\cj$ be the isomorphism induced by $\t^\sh$. Using the
additivity property of characteristic functions, we see that
$\c_{\bK,\bao}=\sum_{\cj\sub\cj^0}\c_{K_\cj,\o_\cj}$. By 28.10, we have $K_\cj=0$ 
unless $\cj\sub\cj_\ss$. By Lemma 35.8, if $\cj\sub\cj_\ss$, $\c_{K_\cj,\o_\cj}$ is 
just like $\c_{K,\o}$ in 35.6, with $\ss$ replaced by $\ss_\cj$. Hence 35.6(a) can be 
applied and it yields
$$\c_{K_\cj,\o_\cj}(Q,Q',gU_Q)
=a_{\uD\l,F,\ss_\cj}\i\tr(T_{\ss_\cj}1_{\uD\l}\r_{g,d}\pr_Q,\fP).$$
We will verify below that
$$a_{\uD\l,F,\ss}=a_{\uD\l,F,\ss_\cj}\text{ for any }\cj\sub\cj_\ss.\tag a$$
We see that 
$$\c_{\bK,\bao}(Q,gQg\i,gU_Q)=a_{\uD\l,F,\ss}\i
                \sum_{\cj\sub\cj_\ss}\tr(T_{\ss_\cj}1_{\uD\l}\r_{g,d}\pr_Q,\fP)$$
is as desired. (We have used the identity 
$C^\ss_{\uD\l}=\sum_{\cj\sub\cj_\ss}T_{\ss_\cj}1_{\uD\l}$ which follows easily from the
definitions.) 

We now verify (a). We may assume that $\cj$ has a single element $j\in\cj_\ss$ (the 
general case can then be obtained by iteration). We have $\ss_\cj=(s'_1,s'_2,\do,s'_r)$
where $s'_i=s_i$ for $i\ne j$ and $s'_j=1$. It suffices to show that, for any 
$k\in[1,r]$ we have 
$$\sqrt{\th^{\uD\l}_F(s_r\do s_{k+1}\cha_{s_k}(-1))}=
\sqrt{\th^{\uD\l}_F(s'_r\do s'_{k+1}\cha_{s'_k}(-1))}.\tag b$$
If $s_k=1$, then $s'_k=1$ and both sides are $1$. If $s_k\in\II$ and $k>j$ then 
$s_{k'}=s'_{k'}$ for $k'\ge k$ and (b) is obvious. If $k=j$ then (b) states that

$\sqrt{\th^{\uD\l}_F(s_r\do s_{j+1}\cha_{s_j}(-1))}=1$;
\nl
this follows from
$$\th^{\uD\l}_F(s_r\do s_{j+1}\cha_{s_j}(a))=1\text{ for all }a\in\FF_q^*\tag c$$
which comes from $j\in\cj_\ss$. Assume now that $s_k\in\II$ and $k<j$. Then for some 
$m\in\ZZ$ we have 
$$s_js_{j-1}\do s_{k+1}\cha_{s_k}=(s_{j-1}\do s_{k+1}\cha_{s_k})\cha_{s_j}^m.$$
Applying $s_rs_{r-1}\do s_{j+1}$ to both sides gives
$$s_rs_{r-1}\do s_{k+1}\cha_{s_k}=(s_rs_{r-1}\do s_{j+1}s_{j-1}\do s_{k+1}\cha_{s_k})
(s_rs_{r-1}\do s_{j+1}\cha_{s_j})^m.$$
Applying $\th^{\uD\l}_F$ to both sides and using (c) gives
$$\th^{\uD\l}_F(s_rs_{r-1}\do s_{k+1}\cha_{s_k}(a))
=\th^{\uD\l}_F(s_rs_{r-1}\do s_{j+1}s_{j-1}\do s_{k+1}\cha_{s_k}(a))$$
for all $a\in\FF_q^*$. We set $a=-1$ and we see that (b) holds for this $k$. This 
proves (a). The proposition is proved.

\subhead 35.11\endsubhead
Let $\cl,\cl'\in\fs_n$. Then $\cl,\cl'\in\fs_{q-1}$. Let $\t_0:F^*\cl@>\si>>\cl$ be as 
in 35.2; let $\t'_0:F^*\cl'@>\si>>\cl'$ be the analogous isomorphism. Let $\l\in\ufs_n$
(resp. $\l'\in\ufs_n$) be the isomorphism class of $\cl$ (resp. $\cl'$). We have 
$\l\in\ufs_{q-1}$, $\l'\in\ufs_{q-1}$. Let $\ss=(s_1,s_2,\do,s_r)$, 
$\ss'=(s'_1,s'_2,\do,s'_{r'})$ be sequences in $\II$ such that 
$s_1s_2\do s_r\uD\l=\l$, $s'_1s'_2\do s'_{r'}\uD\l'=\l'$. Let 
$\bK=\bK^{\ss,\cl}_{J,D},\bao:F^*\bK@>\si>>\bK$ be as in 35.9. Let
$\bK'=\bK^{\ss',\cl'}_{J,D}$ and let $\bao':F^*\bK'@>\si>>\bK'$ be the analogue of
$\bao$ (defined in terms of $\t'_0$). Then $\c_{\bK,\bao}:Z_{J,D}^F@>>>\bbq$,
$\c_{\bK',\bao'}:Z_{J,D}^F@>>>\bbq$ are defined. Let
$$E=\sum_{(Q,Q',gU_Q)\in Z_{J,D}^F}\c_{\bK,\bao}(Q,Q',gU_Q)\c_{\bK',\bao'}(Q,Q',gU_Q).
$$
Using 35.10 for $\bK$ and for $\bK'$, we see that 
$$E=|U_{J,B^*}^F|\i a_{\uD\l,F,\ss}\i a_{\uD\l',F,\ss'}\i\sum_{Q\in\cp_J^F}
\sum_{g\in D^F}\tr(C^\ss_{\uD\l}\r_{g,d}\pr_Q,\fP)\tr(C^{\ss'}_{\uD\l'}\r_{g,d}
\pr_Q,\fP).$$
Setting $g=dg_0$ where $g_0\in G^{0F}$ we have $\r_{g,d}=\r_{d,d}\r_{g_0,1}$. Hence
$$\align&E=|U_{J,B^*}^F|\i a_{\uD\l,F,\ss}\i a_{\uD\l',F,\ss'}\i
\sum\Sb Q\in\cp_J^F\\ g_0\in G^{0F}\eSb\tr(\pr_QC^\ss_{\uD\l}\r_{d,d}\r_{g_0,1},\fP)\\&
\T\tr(\pr_QC^{\ss'}_{\uD\l'}\r_{d,d}\r_{g_0,1},\fP)=|G^{0F}||U_{J,B^*}^F|\i 
a_{\uD\l,F,\ss}\i a_{\uD\l',F,\ss'}\i\tr(XY,\fP\ot\fP)\endalign$$
where
$$X=\sum_{Q\in\cp_J^F}(\pr_Q\ot\pr_Q)((C^\ss_{\uD\l}\r_{d,d})\ot(C^{\ss'}_{\uD\l'}
\r_{d,d})):\fP\ot\fP@>>>\fP\ot\fP,$$
$$Y=|G^{0F}|\i\sum_{g_0\in G^{0F}}(\r_{g_0,1}\ot\r_{g_0,1}):\fP\ot\fP@>>>\fP\ot\fP.$$
We have $XY=YX$ and $Y$ is a projection of $\fP\ot\fP$ onto the subspace 
$(\fP\ot\fP)^{G^{0F}}$ of $G^{0F}$-invariants for the $G^{0F}$-action in which $g_0$ 
acts as $\r_{g_0,1}\ot\r_{g_0,1}$. Hence $X$ maps $(\fP\ot\fP)^{G^{0F}}$ into itself 
and
$$E=|G^{0F}||U_{J,B^*}^F|\i a_{\uD\l,F,\ss}\i a_{\uD\l',F,\ss'}\i
\tr(X,(\fP\ot\fP)^{G^{0F}}).$$
The non-singular symmetric bilinear form 
$$(,):\fP\T\fP@>>>\fU,(f,f')=\sum_{x\in G^{0F}}f(x)f'(x)$$
gives rise to an isomorphism $\fP\ot\fP@>\si>>\End(\fP)$, $f'\ot f''\m[f\m(f,f')f'']$. 
Under this isomorphism, $X$ corresponds to a linear map $X':\End(\fP)@>>>\End(\fP)$,
$$\ph\m\sum_{Q\in\cp_J^F}(\pr_QC^{\ss'}_{\uD\l'}\r_{d,d})\ph({}^t(\pr_QC^\ss_{\uD\l}
\r_{d,d}))$$
where ${}^t$ denotes taking transpose with respect to $(,)$. We have 
$(\r_{g_0,1}f,\r_{g_0,1}f')=(f,f')$ for all $f,f'\in\fP,g_0\in G^{0F}$. Hence 
$\fP\ot\fP@>\si>>\End(\fP)$ restricts to an isomorphism
$(\fP\ot\fP)^{G^{0F}}@>\si>>\End_{G^{0F}}(\fP)$ under which

$X:(\fP\ot\fP)^{G^{0F}}@>>>(\fP\ot\fP)^{G^{0F}}$
\nl
corresponds to the restriction of $X'$ to $\End_{G^{0F}}(\fP)$. It follows that
$$E=|G^{0F}||U_{J,B^*}^F|\i a_{\uD\l,F,\ss}\i a_{\uD\l',F,\ss'}\i
\tr(X',\End_{G^{0F}}(\fP)).$$
From the definitions we have 
$${}^t\r_{d,d}=\r_{d\i,d\i}:\fP@>>>\fP,$$
$${}^t(\pr_Q)=\pr_Q\text{ for all }Q\in\cp_J,$$
$${}^t(k_\nu)=k_{\nu\i}:\fP@>>>\fP\text{ for all }\nu\in N_{G^0}T.$$
In particular,
$${}^tk_{\ds}=k_{\ds\i}=k_{\ds}k_{\cha_s(-1)}:\fP@>>>\fP\text{ for all }s\in\II.$$
We also see that 
$${}^t(1_{\l_1})=1_{\l_1\i}:\fP@>>>\fP\text{ for all }\l_1\in\ufs_{q-1}.$$
For $s\in\II,\l_1\in\ufs_{q-1}$ we have
$$\align{}^t(T_s1_{\l_1})&=\sqrt{\th^{\l_1}_F(\cha_s(-1))}1_{\l_1\i}k_{\ds}
k_{\cha_s(-1)}=\sqrt{\th^{\l_1}_F(\cha_s(-1))}\th^{\l_1}_F(\cha_s(-1))1_{\l_1\i}k_{\ds}
\\&=\th^{\l_1}_F(\cha_s(-1))T_s1_{s\l_1\i}:\fP@>>>\fP,\endalign$$
hence ${}^t(C^s_{\l_1})=\th^{\l_1}_F(\cha_s(-1))C^s_{s\l_1\i}:\fP@>>>\fP$. It follows 
that
$${}^t(C^\ss_{\uD\l})=\d_0C^{\tss}_{s_1s_2\do s_r(\uD\l)\i}=\d_0C^{\tss}_{\l\i}$$
where $\tss=(s_r,\do,s_2,s_1)$ and
$$\d_0=(\th^{\uD\l}_F(\cha_{s_r}(-1)))(\th^{\uD\l}_F(s_r\cha_{s_{r-1}}(-1)))\do
(\th^{\uD\l}_F(s_r\do s_2\cha_{s_1}(-1)))=a_{\uD\l,F,\ss}^2.$$
We see that 
$X'(\ph)=\d_0\sum_{Q\in\cp_J^F}\pr_QC^{\ss'}_{\uD\l'}\r_{d,d}\ph\r_{d\i,d\i}
C^{\tss}_{\l\i}\pr_Q$ for $\ph\in\End(\fP)$.

\subhead 35.12\endsubhead
For $w\in\WW$ we set
$$t_{d,w}=((\e_D(w))\dot{})\i d\dw d\i\in T^F.$$
We show that, for $\l_1\in\ufs_{q-1}$, we have
$$\r_{d,d}T_w1_{\l_1}\r_{d\i,d\i}=\th^{\l_1}_F(d\i t_{d,w}d)\i T_{\e_D(w)}1_{\uD\l_1}:
\fP@>>>\fP.\tag a$$
(We regard $T_w1_{\l_1}$ as an element of $\End_{G^{0F}}\fP=\fT=H_{q-1}^{\sqq}$, see
35.4, 35.5.) We first show that for $\nu\in(N_{G^0}T)^F$ we have
$$\r_{d,d}k_\nu\r_{d\i,d\i}=k_{d\nu d\i}:\fP@>>>\fP.$$
Indeed, for $f\in\fP,x\in G^{0F}$, we have
$$\align&(\r_{d,d}k_\nu\r_{d\i,d\i}f)(x)=(k_\nu\r_{d\i,d\i}f)(d\i xd)\\&
=|U|\i\sum_{x'\in U\nu U}(\r_{d\i,d\i}f)(d\i xdx')=|U|\i\sum_{x'\in U\nu U}
f(dd\i xdx'd\i)\\&=|U|\i\sum_{x''\in Ud\nu d\i U}f(xx'')=(k_{d\nu d\i}f)(x),\endalign$$
as required. Using the equality 
$$\align&T_w1_{\l_1}=\prod_{\a\in R^+,w\i\a\in R^-}\sqrt{\th^{\l_1}_F(\cha(-1))}k_{\dw}
1_{\l_1}\\&=|T^F|\i\prod_{\a\in R^+,w\i\a\in R^-}
\sqrt{\th^{\l_1}_F(\cha(-1))}\sum_{t\in T^F}\th^{\l_1}_F(t)k_{\dw t},\endalign$$
we have
$$\align&\r_{d,d}T_w1_{\l_1}\r_{d\i,d\i}\\&
=|T^F|\i\prod_{\a\in R^+,w\i\a\in R^-}\sqrt{\th^{\l_1}_F(\cha(-1))}
\sum_{t\in T^F}\th^{\l_1}_F(t)\r_{d,d}k_{\dw t}\r_{d\i,d\i}\\&
=|T^F|\i\prod_{\a\in R^+,w\i\a\in R^-}\sqrt{\th^{\l_1}_F(\cha(-1))}
\sum_{t\in T^F}\th^{\l_1}_F(t)k_{d\dw td\i}\\&=|T^F|\i\prod_{\a\in R^+,w\i\a\in R^-}
\sqrt{\th^{\l_1}_F(\cha(-1))}\sum_{t\in T^F}\th^{\l_1}_F(t)k_{d\dw d\i}k_{dtd\i}\\&=
\prod_{\a\in R^+,w\i\a\in R^-}\sqrt{\th^{\uD\l_1}_F(\uD\cha(-1))}k_{(\e_D(w))\dot{}}
k_{t_{d,w}}1_{\uD\l_1}\\&=\th^{\uD\l_1}_F(t_{d,w}\i)\prod_{\a'\in R^+,\e_D(w)\i\a'\in 
R^-}\sqrt{\th^{\uD\l_1}_F(\cha'(-1))}k_{(\e_D(w))\dot{}}
1_{\uD\l_1}\\&=\th^{\l_1}_F(d\i t_{d,w}d)\i T_{\e_D(w)}1_{\uD\l_1}\endalign$$
and (a) is proved.

\subhead 35.13\endsubhead
We write $\Th^J_n:H_n@>>>H_n$ instead of $\Th^J:H_n@>>>H_n$, see 32.22. Let 
$\Th^J_{n,q}:H_n^{\sqq}@>>>H_n^{\sqq}$ be the linear map defined by $\Th^J_n$ by 
extension of scalars. Replacing $n$ by $q-1$ we obtain a linear map 
$\Th^J_{q-1,q}:H_{q-1}^{\sqq}@>>>H_{q-1}^{\sqq}$. We identify
$\End_{G^{0F}}\fP=\fT=H_{q-1}^{\sqq}$, see 35.4, 35.5. We show:
$$\Th^J_{q-1,q}(\ph)=\sum_{Q\in\cp_J^F}\pr_Q\ph\pr_Q:\fP@>>>\fP\tag a$$
for any $\ph\in\End_{G^{0F}}\fP$. First we show that for $\nu\in(N_{G^0}T)^F$ we have
$$\sum_{Q\in\cp_J^F}\pr_Qk_\nu\pr_Q=k_\nu\text{ if }\nu\in Q_{J,B^*};
\sum_{Q\in\cp_J^F}\pr_Qk_\nu\pr_Q=0\text{ if }\nu\n Q_{J,B^*}.\tag b$$
Let $f\in\fP,x\in G^{0F}$. We have
$$\align&(\sum_{Q\in\cp_J^F}\pr_Qk_\nu\pr_Qf)(x)=\sum\Sb Q\in\cp^J;\\
xQ_{J,B^*}x\i=Q\eSb(k_\nu\pr_Qf)(x)\\&=|U|\i\sum\Sb Q\in\cp^J;\\xQ_{J,B^*}x\i=Q;
\\x'\in U\nu U\eSb(\pr_Qf)(xx')=|U|\i\sum\Sb Q\in\cp^J;\\xQ_{J,B^*}x\i=Q;
\\x'\in U\nu U;xx'Q_{J,B^*}x'{}\i x\i=Q\eSb f(xx')
\\&=|U|\i\sum\Sb x'\in U\nu U;\\x'Q_{J,B^*}x'{}\i=Q_{J,B^*}\eSb f(xx')
=|U|\i\sum_{x'\in U\nu U\cap Q_{J,B^*}}f(xx'),\endalign$$
and (b) follows.

It suffices to prove (a) for $\ph=T_w1_\l$ where $w\in\WW,\l\in\ufs_{q-1}$. We have
$$\sum_{Q\in\cp_J^F}\pr_QT_w1_\l\pr_Q=|T^F|\i\prod\Sb\a\in R^+\\ w\i\a\in R^-\endSb
\sqrt{\th^\l_F(\cha(-1))}\sum_{t\in T^F}\th^\l_F(t)\sum_{Q\in\cp_J^F}\pr_Qk_{\dw t}
\pr_Q.$$
By (b) this is zero, unless $\dw\in Q_{J,B^*}$ that is, $w\in W_J$; assuming that 
$w\in W_J$ this equals
$$|T^F|\i\prod_{\a\in R^+,w\i\a\in R^-}\sqrt{\th^\l_F(\cha(-1))}\sum_{t\in T^F}
\th^\l_F(t)k_{\dw t}=T_w1_\l;$$
(a) is proved.

\subhead 35.14\endsubhead
Let $m$ be an integer $\ge 1$. Let $\FF_{q^m}$ be the subfield of $\kk$ consisting of 
$q^m$ elements. Then $F^m:G@>>>G$ is the Frobenius map for an $\FF_{q^m}$-rational 
structure on $G$. The \'epinglage in 35.2 relative to $\FF_q$ can be also regarded as 
an \'epinglage relative to $\FF_{q^m}$. In the setup of 35.11, define 
$\t_0^{(m)}:F^{m*}\cl@>\si>>\cl$, $\t'_0{}^{(m)}:F^{m*}\cl'@>\si>>\cl'$ in terms of
$\t_0,\t'_0$ as in 33.5. We have $\c_{\cl,\t_0^{(m)}}=\th^\l_{F^m}$, 
$\c_{\cl',\t'_0{}^{(m)}}=\th^{\l'}_{F^m}$ where 

$\th^\l_{F^m}(t)=\th^\l_F(t^{1+q+\do+q^{m-1}})$,
$\th^{\l'}_{F^m}(t)=\th^{\l'}_F(t^{1+q+\do+q^{m-1}})$
\nl
for all $t\in T^{F^m}$. Now $\bao^{(m)}:F^{m*}\bK@>\si>>\bK$ (see 33.5) has the same 
relation to $\t_0^{(m)}$ as $\bao:F^*\bK@>\si>>\bK$ (see 35.9) to $\t_0$. Let 
$\bao'{}^{(m)}:F^{m*}\bK'@>\si>>\bK'$ be the analogous isomorphism defined in terms of 
$\t'_0{}^{(m)}$. Then $\c_{\bK,\bao^{(m)}}:Z_{J,D}^{F^m}@>>>\bbq$, 
$\c_{\bK',\bao'{}^{(m)}}:Z_{J,D}^{F^m}@>>>\bbq$ are well defined. 

We choose (as we may) $m_0\in\NN^*_\kk$ such that $(-1)^{m_0}=1$ (in $\kk^*$) and 
$t_{d,w}^{m_0}=1$ (in $T$) for all $w\in\WW$. 

We show that if $m\in m_0\ZZ,m\ge 1$, then $a_{\uD\l,F^m,\ss}=1$. It suffices to show 
that $\th^{\uD\l}_{F^m}\cha(-1)=1$ for any $\a\in R$. Since 
$\uD\l\in\ufs_{q-1}$ we have $\th^{\uD\l}_{F^m}(t)=\th^\l_F(t^{1+q+\do+q^{m-1}})$. 
Hence it suffices to show that $\cha((-1)^{1+q+\do+q^{m-1}})=1$ for any $\a\in R$ or 
that $(-1)^{1+q+\do+q^{m-1}}=1$ (in $\kk^*$). Since $(-1)^q=-1$ (in $\kk^*$), it 
suffices to show that $(-1)^m=1$. This follows from our assumption on $m$ and $m_0$.

Similarly, we see that if $m\in m_0\ZZ,m\ge 1$ then $a_{\uD\l',F^m,\ss'}=1$. 

We show that, if $\ti\l\in\ufs_n,w\in\WW$ and $m\in m_0\ZZ,m\ge 1$, then 
$\th^{\ti\l}_{F^m}(d\i t_{d,w}d)=1$. Since $\ti\l\in\ufs_{q-1}$, we have 
$\th^{\ti\l}_{F^m}(d\i t_{d,w}d)=\th^{\ti\l}_F(d\i t_{d,w}d)^{1+q+\do+q^{m-1}}$. Hence 
it suffices to show that $t_{d,w}^{1+q+\do+q^{m-1}}=1$. Since $t_{d,w}^q=t_{d,w}$, it 
suffices to show that $t_{d,w}^m=1$. This follows from our assumption on $m_0,m$.

We replace $\FF_q$ in 35.1 by $\FF_{q^{m_0}}$ which we rename as $\FF_q$. The results 
above can be reformulated as follows.

(a) {\it If $m\in\ZZ,m\ge 1$, then $a_{\uD\l,F^m,\ss}=1$, $a_{\uD\l',F^m,\ss'}=1$.}

(b) {\it If $\ti\l\in\ufs_n,w\in\WW$ and $m\in\ZZ,m\ge 1$ then }
$\th^{\ti\l}_{F^m}(d\i t_{d,w}d)=1$.

\proclaim{Proposition 35.15}Let $\fa=\fa_D:H_n@>>>H_n$ be as in 34.15. Define 
$\Ph'':H_n@>>>H_n$ by $h\m\Th^J_n(C^{\ss'}_{\uD\l'}\fa(h)C^{\tss}_{\l\i})$ (an 
$\ca$-linear map) and let $\mu(G^0)$ be as in 32.22. If $m\in\ZZ,m\ge 1$, then
$$E_m=\sum_{(Q,Q',gU_Q)\in Z_{J,D}^{F^m}}
\c_{\bK,\bao^{(m)}}(Q,Q',gU_Q)\c_{\bK',\bao'{}^{(m)}}(Q,Q',gU_Q)$$
is obtained by substituting $v^2=q^m$ in $v^{2l(w^0_J)}\mu(G^0)\tr(\Ph'',H_n)$, which 
is a polynomial in $\NN[v^2]$.
\endproclaim
By the arguments in 35.11-35.13 applied with $F^m$ instead of $F$ we see that 
$$E_m=|G^{0F^m}||U_{J,B^*}^{F^m}|\i a_{\uD\l,F^m,\ss}a_{\uD\l',F^m,\ss'}\i
\tr(X'_m,H_{q^m-1}^{\sqrt{q^m}})$$
where $X'_m(h)=\Th^J_{q^m-1,q^m}(C^{\ss'}_{\uD\l'}\x'_m(h)C^{\tss}_{\l\i})$ for 
$h\in H_{q^m-1}^{\sqrt{q^m}}$ and \lb 
$\x'_m:H_{q^m-1}^{\sqrt{q^m}}@>>>H_{q^m-1}^{\sqrt{q^m}}$ is the linear map given by
$$T_w1_{\ti\l}\m\th^{\ti\l}_{F^m}(d\i t_{d,w}d)\i T_{\e_D(w)}1_{\uD\ti\l}$$
for $w\in\WW,\ti\l\in\ufs_{q^m-1}$. Clearly, $X'_m(T_w1_{\ti\l})=0$ unless 
$\ti\l\in\ufs_n$ and, if this condition is satisfied, then $X'_m(T_w1_{\ti\l})$ is a 
linear combination of elements of the form $T_{w'}1_{\ti\l'}$ with $\ti\l'\in\ufs_n$. 
It follows that
$$\tr(X'_m,H_{q^m-1}^{\sqrt{q^m}})=\tr(\Ph''_m,H_n^{\sqrt{q^m}})$$
where $\Ph''_m(h)=\Th^J_{n,q^m}(C^{\ss'}_{\uD\l'}\x''_m(h)C^{\tss}_{\l\i})$ and
$\x''_m:H_n^{\sqrt{q^m}}@>>>H_n^{\sqrt{q^m}}$ is the restriction of $\x'_m$. We now use
35.14(a),(b). The proposition follows, except for the assertion "which is a polynomial 
in $\NN[v^2]$". That assertion follows from 32.22(a) and the second equality in 
32.23(a). 

\subhead 35.16\endsubhead
We now assume that $J=\II$. Let $\ci_n$ be a set of representatives for the isomorphism
classes of character sheaves contained in $\hat D^\cl$ for some $\cl\in\fs_n$. Then 
$\ci_n$ is finite and we can find an integer $m_1\ge 1$ such that $F^{m_1*}A\cong A$ 
for any $A\in\ci_n$. We replace $\FF_q$ in 35.14 by $\FF_{q^{m_1}}$ which we rename as 
$\FF_q$. We now have $F^*A\cong A$ for any $A\in\ci_n$. For $A\in\ci_n$, the 
Verdier dual $\fD(A)$ is isomorphic to an object in $\ci_n$. (See 28.18(a).) For 
each $A\in\ci_n$ we choose isomorphisms $\k_A:F^*A@>\si>>A$,
$\k'_A:F^*\fD(A)@>\si>>\fD(A)$ so that the following holds: if $\co$ is an open dense 
$F$-stable subset of $\supp(A)=\supp(\fD(A))$ on which $\ch^{-e}(A),\ch^{-e}(\fD(A))$
are local systems ($e=\dim\co$) then for any $y\in\co$ and 
any $m\ge 1$ such that $F^m(y)=y$,

(i) $\sqq^{m(e-\dim D)}\k_A^{(m)}:\ch^{-e}(A)_y@>>>\ch^{-e}(A)_y$ is of finite order;

(ii) $\sqq^{m(e-\dim D)}\k'_A{}^{(m)}:\ch^{-e}(\fD(A))_y@>>>\ch^{-e}(\fD(A))_y$ is of 
finite order;

(iii)  the isomorphism $\ch^{-e}\fD(A)_{F(y)}@>\si>>\ch^{-e}\fD(A)_y$ (that is, 
$\ch^{-e}(A)\che{}_{F(y)}@>\si>>\ch^{-e}(A)\che{}_y$) induced by $\k'_A$ is 
$q^{\dim D-e}$ times the transpose  inverse of the isomorphism 
$\ch^{-e}(A)_{F(y)}@>\si>>\ch^{-e}(A)_y$ induced by $\k_A$.

Note that $\k'_A$ is uniquely determined by $\k_A$ and that (ii) follows from (i) and
(iii).

Let $\cl,\cl',\ss,\ss',\bK,\bK',\bao,\bao'$ be as in 35.11. Since $\bK,\bK'$ are 
semisimple, we have canonically (for $i,i'\in\ZZ$):
$${}^pH^i(\bK)=\op_{A\in\ci_n}(A\ot V_{A,i,\ss,\cl}),\tag a$$
$${}^pH^{i'}(\bK')=\op_{A'\in\ci_n}(\fD(A')\ot V'_{A',i',\ss',\cl'})\tag b$$
where $V_{A,i,\ss,\cl},V'_{A',i',\ss',\cl'}$ are finite dimensional $\bbq$-vector 
spaces endowed with endomorphisms 

$\ps_A:V_{A,i,\ss,\cl}@>>>V_{A,i,\ss,\cl}$, 
$\ps'_{A'}:V'_{A',i',\ss',\cl'}@>>>V'_{A',i',\ss',\cl'}$
\nl
such that under (a) (resp. (b)) the map $\op_A(\k_A\ot\ps_A)$ (resp. 
$\op_{A'}(\k'_{A'}\ot\ps'_{A'})$) corresponds to the isomorphism 
$F^*({}^pH^i(\bK))@>\si>>{}^pH^i(\bK)$ (resp. 
$F^*({}^pH^{i'}(\bK'))@>\si>>{}^pH^{i'}(\bK')$) induced by $\bao$ (resp. $\bao'$). 

{\it In the remainder of this section we assume that $D$ is clean (see 33.4(b)).}

\proclaim{Proposition 35.18}(a) With each $A\in\ci_n$ one can associate 
$\sgn_A\in\{1,-1\}$ with the following property: if $A$ is a direct summand of 
${}^pH^i(\bK^{\ss,\cl}_D)$ where $\ss,\cl$ are as in 35.11, then 
$(-1)^{i+\dim G}=\sgn_A$.

(b) With each $A\in\ci_n$ one can attach an element $\x_A\in\bbq^*$ such that the 
following hold: $\x_A$ is an algebraic number all of whose complex conjugates have 
absolute value $1$; for any $\ss,\cl$ as in 35.11 and any $i\in\ZZ$, $\ps_A$ is equal 
to $\x_A\sqq^{i-\dim G}$ times a unipotent automorphism of $V_{A,i,\ss,\cl}$; for any 
$\ss',\cl'$ as in 35.11 and any $i\in\ZZ$, $\ps'_A$ is equal to $\x_A\i\sqq^{i-\dim G}$
times a unipotent automorphism of $V'_{A,i,\ss',\cl'}$.
\endproclaim
From the definitions we we see that, in the setup of 35.16, we have
$$\tr(\bao,\ch^j_g({}^pH^i(\bK)))=\sum_{A\in\ci_n}\tr(\k_A,\ch^j_gA)
\tr(\ps_A,V_{A,i,\ss,\cl}),$$
$$\tr(\bao',\ch^{j'}_g({}^pH^{i'}(\bK')))
=\sum_{A'\in\ci_n}\tr(\k'_{A'},\ch^{j'}_g\fD(A'))\tr(\ps'_{A'},V'_{A',i',\ss',\cl'}),$$
for all $g\in D^F$ and all $i,j,i',j'$. Taking alternating sums over $i,j$ or $i',j'$ 
and using
$$\c_{\bK,\bao}(g)=\sum_{i,j}(-1)^{i+j}\tr(\bao,\ch^j_g({}^pH^i(\bK))),$$
$$\c_{\bK',\bao'}(g)=\sum_{i',j'}(-1)^{i'+j'}\tr(\bao',\ch^{j'}_g({}^pH^{i'}(\bK'))),$$
we obtain
$$\align&\c_{\bK,\bao}(g)=\sum_{A\in\ci_n}\c_{A,\k_A}(g)\sum_i(-1)^i
\tr(\ps_A,V_{A,i,\ss,\cl}),\\&\c_{\bK',\bao'}(g)=\sum_{A'\in\ci_n}\c_{\fD(A'),\k'_{A'}}
(g)\sum_{i'}(-1)^{i'}\tr(\ps'_{A'},V'_{A',i',\ss',\cl'}).\tag c\endalign$$
Since $\bK$ (resp. $\bK'$) is pure of weight $0$, we see that ${}^pH^i(\bK)$ (resp. 
${}^pH^{i'}(\bK')$) is pure of weight $i$ (resp. $i'$). By our choice of $\k_A$ (resp. 
$\k'_A$) we see that $(A,\k_A)$ and $(\fD(A),\k'_A)$ are pure of weight $\dim G$ for 
$A\in\ci_n$. Using 35.16(a),(b), we deduce that

(d) {\it $(V_{A,i,\ss,\cl},\ps_A)$ is pure of weight $i-\dim G$ and 
$(V'_{A,i',\ss',\cl'},\ps'_A)$ is pure of weight $i'-\dim G$.}
\nl
Using (c) we have
$$\align&|G^{0F}|\i\sum_{g\in D^F}\c_{\bK,\bao}(g)\c_{\bK',\bao'}(g)
=\sum_{A,A'\in\ci_n}|G^{0F}|\i\sum_{g\in D^F}\c_{A,\k_A}(g)\c_{\fD(A'),\k'_{A'}}(g)\\&
\T\sum_{i,i'}(-1)^{i+i'}\tr(\ps_A,V_{A,i,\ss,\cl})\tr(\ps'_{A'},V'_{A',i',\ss',\cl'}).
\tag e\endalign$$
Using 24.18 (which is applicable by our cleanness assumption and by Corollary 31.15) we
see that for any $A,A'\in\ci_n$ we have
$$|G^{0F}|\i\sum_{g\in D^F}\c_{A,\k_A}(g)\c_{\fD(A'),\k'_{A'}}(g)=\d_{A,A'}.$$
Hence (e) becomes
$$|G^{0F}|\i\sum_{g\in D^F}\c_{\bK,\bao}(g)\c_{\bK',\bao'}(g)=\sum\Sb A\in\ci_n\\i,i'
\eSb(-1)^{i+i'}\tr(\ps_A,V_{A,i,\ss,\cl})\tr(\ps'_A,V'_{A,i',\ss',\cl'}).$$
This remains true if $F$ is replaced by $F^m$, where $m\ge 1$. Thus we have
$$\align&|G^{0F^m}|\i\sum_{g\in D^{F^m}}\c_{\bK,\bao^{(m)}}(g)
\c_{\bK',\bao'{}^{(m)}}(g)\\&
=\sum_j(-1)^j\sum_{A\in\ci_n}\sum_{i,i';i+i'=j}\tr(\ps_A^m,V_{A,i,\ss,\cl})
\tr(\ps'_A{}^m,V'_{A,i',\ss',\cl'})\tag f\endalign$$
with $\bao^{(m)},\bao'{}^{(m)}$ as in 35.14. Using 35.15, we may rewrite the previous 
equality as follows:
$$\sum_j(-1)^j\sum_fa_{j,f}^m=\Pi(q^m)$$
where $a_{j,f}$ are the eigenvalues of $\ps_A\ot\ps'_A$ on 
$\op_{i,i';i+i'=j}V_{A,i,\ss,\cl}\ot V'_{A,i',\ss',\cl'}$ and $\Pi$ is a polynomial 
with coefficients in $\NN$. By (d), each $a_{j,f}$ is an algebraic number all of whose 
complex conjugates have absolute value squared equal to $q^{j-2\dim G}$. It follows 
that, for fixed $j$, the set $\{a_{j,f}\}$ is empty if $j$ is odd and that each 
$a_{j,f}$ is equal to $q^{j/2-\dim G}$ if $j$ is even. This implies that, for any 
$A\in\ci_n$, $V_{A,i,\ss,\cl}\ot V'_{A,i',\ss',\cl'}$ is $0$ if $i+i'$ is odd and, if 
$i+i'$ is even, any eigenvalue of $\ps_A$ on $V_{A,i,\ss,\cl}$ multiplied by any 
eigenvalue of $\ps'_A$ on $V'_{A,i',\ss',\cl'}$ gives $q^{(i+i')/2-\dim G}$. Since for 
$A\in\ci_n$ we have $V'_{A,i',\ss',\cl'}\ne 0$ for some $\ss',\cl'$ as in 
35.11, we see that the parity of $i$ such that $V_{A,i,\ss,\cl}\ne 0$ for some 
$\ss,\cl$ as in 35.11 is an invariant of $A$ and that, for any eigenvalue $\x$ of 
$\ps_A$ on $V_{A,i,\ss,\cl}$, the product $\x\sqq^{\dim G-i}$ is also an invariant of 
$A$. The proposition follows.

\subhead 35.19\endsubhead 
As in 34.15, let $\{E_u;u\in\cu\}$ be a set of representatives for the isomorphism 
classes of simple modules for $H_n^1$. Let $m\ge 1$. Using 35.18, we can rewrite 
35.18(f) as follows:
$$\align&|G^{0F^m}|\i\sum_{g\in D^{F^m}}\c_{\bK,\bao^{(m)}}(g)\c_{\bK',
\bao'{}^{(m)}}(g)\\&=\sum_{A\in\ci_n}\sum_{i,i'}(-1)^{i+i'}\dim V_{A,i,\ss,\cl}
\dim V'_{A,i',\ss',\cl'}q^{m(i-\dim G)/2}q^{m(i'-\dim G)/2}.\tag a\endalign$$
By 35.15 (with $J=\II$), the left hand side of (a) is $\tr(\Ph'',H_n)|_{v=\sqrt{q^m}}$.
We have
$$\align&\sum_{A\in\ci_n}(\sum_i(-v)^iv^{-\dim G}\dim V_{A,i,\ss,\cl})
(\sum_{i'}(-v)^{i'}v^{-\dim G}\dim V'_{A,i',\ss',\cl'})\\&=\tr(\Ph'',H_n)
=\sum_{u\in\cu^\fa}\tr(C^{\ss'}_{\uD\l'}\tT_{\uD},E_u^v)
\tr(\tT_{\uD}\i C^{\tss}_{\l\i},E_u^v)\\&
=\sum_{u\in\cu^\fa}\tr(C^{\ss'}_{\uD\l'}\tT_{\uD},E_u^v)
\tr((C^{\tss}_{\l\i})^\flat\tT_{\uD},E_u^v)^\sp\\&
=\sum_{u\in\cu^\fa}\tr(C^{\ss'}_{\uD\l'}\tT_{\uD},E_u^v)
\tr(C^\ss_{\uD\l\i}\tT_{\uD},E_u^v)^\sp.\tag b\endalign$$
(The first equality in (b) comes from the fact that (a) holds for any integer $m\ge 1$.
Here $\tr(\Ph'',H_n)$, as in 35.15 with $J=\II$, can be replaced by a sum of products 
of traces, see 34.14(a). This gives the second equality in (b) where the notation of 
34.15 is used. The third equality in (b) comes from 34.17 since 
$\tT_{\uD}\i C^{\tss}_{\l\i}$ is an $\ca$-linear combination of elements of the form
$\tT_w1_{\l\i},w\in\WW^D$. The fourth equality in (b) follows from the definitions
using $s_r\do s_2s_1\l\i=\uD\l\i$.) Using the definitions and 34.19(a),(c), we obtain
for $A\in\ci_n$:
$$\align&\sum_i(-v)^iv^{-\dim G}\dim V_{A,i,\ss,\cl}
=(\sum_i(-v)^iv^{-\dim G}\dim V_{A,i,\ss,\cl})^\sp\\&
=(\z^A(C^\ss_{\uD\l}[D]))^\sp=\sum_{u\in\cu^\fa}
(b_{A,u}^v)^\sp\tr(C^\ss_{\uD\l}\tT_{\uD},E_u^v)^\sp,\tag b\endalign$$
$$\align&\sum_{i'}(-v)^{i'}v^{-\dim G}\dim V'_{A,i',\ss',\cl'}=
\z^{\fD(A)}(C^{\ss'}_{\uD\l'}[D])\\&
=\sum_{u'\in\cu^\fa}b_{\fD(A),u'}^v\tr(C^{\ss'}_{\uD\l'}\tT_{\uD},E_{u'}^v)
=\sum_{u'\in\cu^\fa}b_{A,u'}^v\tr(C^{\ss'}_{\uD\l'{}\i}\tT_{\uD},E_{u'}^v).\tag c
\endalign$$
Introducing (b),(c) in (a) we obtain
$$\align&\sum_{A\in\ci_n}
(\sum_{u'\in\cu^\fa}b_{A,u'}^v\tr(C^{\ss'}_{\uD\l'{}\i}\tT_{\uD},E_{u'}^v))
(\sum_{u\in\cu^\fa}(b_{A,u}^v)^\sp\tr(C^\ss_{\uD\l}\tT_{\uD},E_u^v)^\sp)
\\&=\sum_{u\in\cu^\fa}\tr(C^{\ss'}_{\uD\l'}\tT_{\uD},E_u^v)
\tr(C^\ss_{\uD\l\i}\tT_{\uD},E_u^v)^\sp,\endalign$$
that is,
$$\sum_{u,u'\in\cu^\fa}(\sum_{A\in\ci_n}b_{A,u'}^v(b_{A,u}^v)^\sp-\d_{u,u'})
\tr(C^{\ss'}_{\uD\l'{}\i}\tT_{\uD},E_{u'}^v)\tr(C^\ss_{\uD\l}\tT_{\uD},E_u^v)^\sp=0.
\tag d$$
Recall that here $\ss$ is assumed to satisfy $s_1s_2\do s_r\uD\l=\l$. The
$\ca$-submodule of $H^D_n$ spanned by the elements $C^\ss_{\uD\l}\tT_{\uD}$ with $\ss$ 
as above (and $\l$ fixed) is just $1_\l H_n\tT_{uD}1_\l$ hence in (d) we may replace 
$C^\ss_{\uD\l}\tT_{\uD}$ by any element in $1_\l H_n\tT_{uD}1_\l$ and the equality
remains true. Similarly, we may replace $C^{\ss'}_{\uD\l'{}\i}\tT_{\uD}$ by any element
in $1_{\l'{}\i}H_n\tT_{\uD}1_{\l'{}\i}$ and the equality remains true. Thus we have
$$\sum_{u,u'\in\cu^\fa}(\sum_{A\in\ci_n}b_{A,u'}^v(b_{A,u}^v)^\sp-\d_{u,u'})
\tr(\tT_{x'}1_{\uD\l'{}\i}\tT_{\uD},E_{u'}^v)\tr(\tT_x1_{\uD\l}\tT_{\uD},E_u^v)^\sp)
=0\tag e$$
for any $x,x'\in\WW$ such that 

(f) $x\uD\l=\l,x'\uD\l'{}\i=\l'{}\i$.
\nl
Now (e) remains true even if (f) does not hold. (If for example, $x\uD\l=\l_1\ne\l$
then for any $u\in\cu^\fa$ we have
$$\tr(\tT_x1_{\uD\l}\tT_{\uD},E_u^v)=\tr(1_{\l_1}\tT_x1_{\uD\l}\tT_{\uD},E_u^v)=
\tr(\tT_x\tT_{\uD}1_\l1_{\l_1},E_u^v)=0,$$
hence the left hand side of (e) is zero.) We now multiply both sides of (e) by

$\tr(\tT_{x'}1_{\uD\l'{}\i}\tT_{\uD},E_{u'_1}^v)^\sp
\tr(\tT_x1_{\uD\l}\tT_{\uD},E_{u_1}^v)$
\nl
where $u_1,u'_1\in\cu^\fa$ and sum the resulting equalities over all $x,x'\in\WW$ and 
$\l,\l'\in\ufs_n$. We obtain
$$\align&\sum_{u,u'\in\cu^\fa}(\sum_{A\in\ci_n}b_{A,u'}^v(b_{A,u}^v)^\sp-\d_{u,u'})
(\sum_{x,\l}\tr(\tT_x1_{\uD\l}\tT_{\uD},E_{u_1}^v)
\tr(\tT_x1_{\uD\l}\tT_{\uD},E_u^v)^\sp))\\&
\T(\sum_{x',\l'}\tr(\tT_{x'}1_{\uD\l'{}\i}\tT_{\uD},E_{u'}^v)
\tr(\tT_{x'}1_{\uD\l'{}\i}\tT_{\uD},E_{u'_1}^v)^\sp)=0.\endalign$$
Using 34.18(a) we deduce
$$\sum_{u,u'\in\cu^\fa}(\sum_{A\in\ci_n}b_{A,u'}^v(b_{A,u}^v)^\sp-\d_{u,u'})
\d_{u',u'_1}f_{u'_1}^v\dim E_{u'_1}\d_{u,u_1}f_{u_1}^v\dim E_{u_1}=0,$$
that is,
$$(\sum_{A\in\ci_n}b_{A,u'_1}^v(b_{A,u_1}^v)^\sp-\d_{u_1,u'_1})
f_{u'_1}^v\dim E_{u'_1}f_{u_1}^v\dim E_{u_1}=0.$$
Since $f_{u'_1}^v\dim E_{u'_1}f_{u_1}^v\dim E_{u_1}\ne 0$, we deduce
$$\sum_{A\in\ci_n}b_{A,u'_1}^v(b_{A,u_1}^v)^\sp=\d_{u_1,u'_1}\tag g$$
for any $u_1,u'_1\in\cu^\fa$.

\subhead 35.20\endsubhead 
Let $A_0\in\ci_n$. Then $\supp(A_0)$ is the closure of a stratum $Y=Y_{L,S}$ of $D$. 
Let $e=\dim Y$. Replacing $\FF_q$ by a finite extension, we may assume that the 
following holds: 

(a) {\it there exists a finite group $\G$ and a sequence of maps $\G@>>>Y^{F^m}$,
$\g\m y_{\g,m}$ ($m=1,2,3,\do$) such that for any $m\ge 1$ and any $A\in\ci_n$ we 
have}
$$q^{-(\dim D-e)m}|\G|\i\sum_{\g\in\G}\c_{A,\k_A^{(m)}}(y_{\g,m})
\c_{\fD(A_0),\k'_{A_0}{}^{(m)}}(y_{\g,m})=\d_{A_0,A'}.$$
(We apply 33.7 with $\ci=\ci_n$.)

\subhead 35.21\endsubhead 
We identify $H_n$ with a subalgebra of $H_{q-1}$ by $\tT_w1_\l\m\tT_w1_\l$ for 
$w\in\WW$, $\l\in\ufs_n$. (We have $\ufs_n\sub\ufs_{q-1}$ since $n$ divides $q-1$.) 
Similarly, for any $\k\in\fU$ we identify $H_n^\k$ with a subalgebra of $H_{q-1}^\k$. 

Let $\{\tE_u;u\in\ti\cu\}$ be a set of representatives for the isomorphism classes of 
simple modules for $H_{q-1}^1$. We may assume that $\cu\sub\ti\cu$, that for $u\in\cu$ 
we have $\tE_u=E_u$ as an $H_n^1$-module with $1_\l$ acting as $0$ for 
$\l\in\ufs_{q-1}-\ufs_n$ and for $u\in\ti\cu-\cu$, $H_n^1$ acts on $\tE_u$ as zero. For
$u\in\ti\cu$, the $H_{q-1}^{\sqq}$-module $\tE_u^{\sqq}$ is defined as in 34.15. If
$u\in\cu$ we have again $\tE_u^{\sqq}=E_u^{\sqq}$ as an $H_n^{\sqq}$-module with $1_\l$
acting as $0$ for $\l\in\ufs_{q-1}-\ufs_n$. For $u\in\ti\cu$ we set 
$V_u=\Hom_{H_{q-1}^{\sqq}}(\tE_u^{\sqq},\fP)$ where $\fP$ is regarded as an 
$H_{q-1}^{\sqq}=\fT$-module as in 35.4, 35.5. Since $\fT=\End_{G^{0F}}(\fP)$ (see 
35.5), the $G^{0F}$-module structure on $\fP$ makes $V_u$ into an irreducible 
$G^{0F}$-module and we have an isomorphism 
$$\vt:\op_{u\in\ti\cu}(\tE_u^{\sqq}\ot V_u)@>\si>>\fP,e\ot x\m x(e),e\in\tE_u^{\sqq},
x\in V_u.$$ 
Hence $\fP=\op_{u\in\ti\cu}\fP_u$ where $\fP_u=\vt(\tE_u^{\sqq}\ot V_u)$. For 
$u\in\cu^\fa$ and $x\in V_u$ we define $R_{u,x}:\tE_u^{\sqq}@>>>\fP$ by 
$R_{u,x}(e)=\r_{d,d}(x(\tT_{\uD}\i e))$. We show that 

(a) $R_{u,x}(he)=hR_{u,x}(e)$ for $h\in H_{q-1}^{\sqq}$, $e\in\tE_u^{\sqq}$. 
\nl
If $h\in H_{q-1}^{\sqq}1_\l$ with $\l\in\ufs_{q-1}-\ufs_n,w\in\WW$ then both sides of
(a) are zero (we use 35.12(a)). Hence we may assume that $h\in H_n^{\sqq}$. Recall that
$\r_{d,d}(\fa\i(h))=h\r_{d,d}:\fP@>>>\fP$ (see 35.12(a) and 35.14(b)). Hence
$$\align R_{u,x}(he)&=\r_{d,d}(x(\tT_{\uD}\i(he)))
=\r_{d,d}(x((\fa\i(h))(\tT_{\uD}\i e))=\r_{d,d}(\fa\i(h))(x(\tT_{\uD}\i e))\\&
=h\r_{d,d}(x(\tT_{\uD}\i e))=hR_u(x)(e),\endalign$$
as required. We see that $R_{u,x}\in V_u$. Thus, we have a map $R_u:V_u@>>>V_u$,
$x\m R_{u,x}$. From the definitions we have

$\vt(\tT_{\uD}e\ot R_{u,x})=\r_{d,d}\vt(e\ot x)$ for $e\in E_u^{\sqq},x\in V_u$.
\nl
In particular, $\r_{d,d}:\fP@>>>\fP$ maps $\fP_u$ into itself. On the 
other hand it maps $\fP_u$ for $u\in\cu-\cu^\fa$ into $\fP_{u'}$ for some 
$u'\in\cu,u'\ne u$. It also maps $\op_{u\in\ti\cu-\cu}\fP_u$ into itself. Hence if 
$h\in H_n^{\sqq}$ and $g_0\in G^{0F}$ then $h\r_{d,d}\r_{g_0,1}$ acts as zero on
$\op_{u\in\ti\cu-\cu}\vt(E_u^{\sqq}\ot V_u)$ and we have 
$$\tr(h\r_{d,d}\r_{g_0,1},\fP)=\sum_{u\in\cu^\fa}
\tr(h\tT_{\uD},E_u^{\sqq})\tr(R_ug_0,V_u).$$
Now, from Lemma 35.10 (with $J=\II$) we have 
$\c_{\bK,\bao}(g)=\tr(C^\ss_{\uD\l}\r_{g,d},\fP)$ for any $g\in D^F$. (Recall that 
$a_{\uD\l,F,\ss}=1$, see 35.14(a).) Hence
$$\c_{\bK,\bao}(g)
=\sum_{u\in\cu^\fa}\tr(C^\ss_{\uD\l}\tT_{\uD},E_u^{\sqq})\tr(R_ud\i g,V_u).\tag b$$
We have
$$\align&\c_{\bK,\bao}(g)=\sum_{A\in\ci_n}\c_{A,\k_A}(g)\sum_i(-\sqq)^i\sqq^{-\dim G}
\x_A\dim V_{A,i,\ss,\cl}\\&=\sum_{A\in\ci_n}\c_{A,\k_A}(g)\x_A\sum_{u\in\cu^\fa}
b_{A,u}^{\sqq}\tr(C^\ss_{\uD\l}\tT_{\uD},E_u^{\sqq}).\tag c\endalign$$
(The first equality follows from 35.18(b),(c). The second equality is obtained from the
identity 

$\sum_i(-v)^iv^{-\dim G}\dim V_{A,i,\ss,\cl}=\sum_{u\in\cu^\fa}
b_{A,u}^v\tr(C^\ss_{\uD\l}\tT_{\uD},E_u^v)$, 
\nl
see 34.19, under the specialization $v\m\sqq$. Note that the rational function 
$b_{A,u}^v$ does not have a pole at $v=\sqq$; indeed, in 34.19(b), we have
$\z^A(\tT_w1_{\uD\l}[D])\in\ca$, $\tr(\tT_w1_{\uD\l}\tT_{\uD},E_u^v)^\sp\in\fU[v,v\i]$ 
by 34.15(a), while $f_u^v|_{v=\sqq}=f_u^{\sqq}\ne 0$ by 34.15(b). Hence $b_{A,u}^v$ can
be specialized for $v=\sqq$ and yields a value $b_{A,u}^{\sqq}\in\fU$.) From (b),(c) we
deduce
$$\sum_{u\in\cu^\fa}\tr(C^\ss_{\uD\l}\tT_{\uD},E_u^{\sqq})
(\tr(R_ud\i g,V_u)-\sum_{A\in\ci_n}\c_{A,\k_A}(g)\x_Ab_{A,u}^{\sqq})=0.$$
As in the paragraph preceding 35.19(e) we see that here we may replace
$C^\ss_{\uD\l}\tT_{\uD}$ by any element in $1_\l H_n\tT_{\uD}1_\l$. Thus we have
$$\sum_{u\in\cu^\fa}\tr(\tT_x1_{\uD\l}\tT_{\uD},E_u^{\sqq})
(\tr(R_ud\i g,V_u)-\sum_{A\in\ci_n}\c_{A,\k_A}(g)\x_Ab_{A,u}^{\sqq})=0\tag d$$
for any $x\in\WW$ such that $x\uD\l=\l$; moreover, if $x\uD\l\ne\l$ then \lb
$\tr(\tT_x1_{\uD\l}\tT_{\uD},E_u^{\sqq})=0$ so that (d) holds again. We see that (d)
holds for any $x\in\WW$ and any $\l\in\ufs_n$. We now multiply both sides of (d) by
$\tr(\tT_x1_{\uD\l}\tT_{\uD},E_{u_1}^{\sqq})^\sp=0$ where $u_1\in\cu^\fa$ and sum the
resulting equalities over all $x\in\WW$ and $\l\in\ufs_n$. We obtain
$$\align&\sum_{u\in\cu^\fa}(\sum_{x,\l}
\tr(\tT_x1_{\uD\l}\tT_{\uD},E_u^{\sqq})\tr(\tT_x1_{\uD\l}\tT_{\uD},E_{u_1}^{\sqq})^\sp)
\\&\T(\tr(R_ud\i g,V_u)-\sum_{A\in\ci_n}\c_{A,\k_A}(g)\x_Ab_{A,u}^{\sqq})=0.\endalign$$
Using 34.18(b) with $\k=\sqq=\sqq^\sp$ we deduce that
$$\sum_{u\in\cu^\fa}\d_{u,u_1}f_u^{\sqq}\dim E_u
(\tr(R_ud\i g,V_u)-\sum_{A\in\ci_n}\c_{A,\k_A}(g)\x_Ab_{A,u}^{\sqq})=0$$
that is, 
$$f_{u_1}^{\sqq}\dim E_{u_1}(\tr(R_{u_1}d\i g,V_{u_1})
-\sum_{A\in\ci_n}\c_{A,\k_A}(g)\x_Ab_{A,u_1}^{\sqq})=0.$$
Since 
$f_{u_1}^{\sqq}\dim E_{u_1}\ne 0$, we deduce
$$\tr(R_{u_1}d\i g,V_{u_1})=\sum_{A\in\ci_n}\c_{A,\k_A}(g)\x_Ab_{A,u_1}^{\sqq}=0.\tag e
$$
We show:

(f) {\it $R_ud\i g:V_u@>>>V_u$ has finite order.}
\nl
For any $t\ge 1$ we have $R_u^t(x)=\r_{d,d}^t(x(\tT_{\uD}^{-t}e))$, hence $R_u^c=1$ for
some $c\ge 1$. For $g_0\in G^{0F}$ we have

$R_ug_0(x)(e)=\r_{d,d}\r_{g_0,1}x(\tT_{\uD}\i e)=\r_{dg_0d\i,1}\r_{d,d}x\tT_{\uD}\i e$ 
\nl
hence $R_ug_0=(dg_0d\i)R_u$. From this we deduce

$(R_ug_0)^t=(dg_0d\i)(d^2g_0d^{-2})\do(d^tg_0d^{-t})R_u^t$
\nl
for $t\ge 1$. We have $d^a=1$ for some $a\ge 1$. Let 

$g_1=(dg_0d\i)(d^2g_0d^{-2})\do(d^ag_0d^{-a})$.
\nl
We have $g_1^b=1$ for some $b\ge 1$. We have $(R_ug_0)^{ab}=g_1^bR_u^{ab}=R_u^{ab}$. 
Thus $(R_ug_0)^{abc}=R_u^{abc}=1$. Taking $g_0=d\i g$ we see that (f) holds.

From (f) we see that $\tr(R_ud\i g,V_u)$ is a cyclotomic integer. Introducing
this in (e) we see that

(g) {\it $\sum_{A\in\ci_n}\c_{A,\k_A}(g)\x_Ab_{A,u}^{\sqq}$ is a cyclotomic integer for
any $u\in\cu^\fa,g\in D^F$}.
 
\proclaim{Lemma 35.22} In the setup of 35.20, let $u\in\cu^\fa$. Assume that $E_u$ is 
quasi-rational (see 34.20). Then $b_{A_0,u}^v\in\et\QQ[v,v\i]$ where $\et$ is a root of
$1$.
\endproclaim
Note that 35.21(g) remains true if $\FF_q$ is replaced by $\FF_{q^m}$ where $m\ge 1$.
Thus:
$$\sum_{A\in\ci_n}\c_{A,\k_A^{(m)}}(g)\x_A^mb_{A,u}^{\sqrt{q^m}}\tag a$$
is a cyclotomic integer for any $m\ge 1,g\in D^{F^m}$. We use notation of 35.20. Assume
that $g\in Y^{F_m}$. We multiply (a) by 
$\sqq^{m(e-\dim D)}\c_{\fD(A_0),\k'_{A_0}{}^{(m)}}(g)$
which is a cyclotomic integer by 35.16(ii). We obtain again a cyclotomic integer. Now 
take $g=y_{\g,m}$ (see 35.20) and sum over all $\g\in\G$ (see 35.20). We see that
$$\sqq^{m(e-\dim D)}\sum_{\g\in\G}\sum_{A\in\ci_n}\c_{A,\k_A^{(m)}}(y_{\g,m})
\c_{\fD(A_0),\k'_{A_0}{}^{(m)}}(y_{\g,m})\x_A^mb_{A,u}^{\sqrt{q^m}}\tag b$$
is a cyclotomic integer for any $m\ge 1$. Using 35.20(a) we see that (b) equals
$$\sqq^{m(\dim D-e)}|\G|\x_{A_0}^mb_{A_0,u}^{\sqrt{q^m}}\tag c$$
which is therefore a cyclotomic integer for any $m\ge 1$. By 35.22 we have 
$b_{A_0,u}^v=\et Q(v)$ where $\et\in\fU$ is a root of $1$ and $Q(v)\in\QQ(v)$. Let $K$ 
be a subfield of $\bbq$ such that $K$ is a finite Galois extension of $\QQ$ of degree 
$a$ which contains $\et$ and $\x_{A_0}$. Let $N:K@>>>\QQ$ be the norm map. Since all 
complex conjugates of $\x_{A_0}$ have absolute value $1$ we see that 
$N(\x_{A_0})=\pm 1$. We have also $N(\et)=\pm 1$. Hence applying $N$ to (c) (with 
$m=2m'$) we see that $|\G|^aq^{am'(\dim D-e)}Q(q^{m'})^a$ is a cyclotomic integer. This
being also a rational number, is an ordinary integer. Let 

$R(v)=|\G|^av^{a(\dim D-e)}Q(v)^a\in\QQ(v)$. 
\nl
We see that $R(q^{m'})\in\ZZ$ for any integer $m'\ge 1$. This forces $R(v)\in\QQ[v]$. 
Thus, $(v^{\dim D-e}Q(v))^a\in\QQ[v]$. It follows that $v^{\dim D-e}Q(v)\in\QQ[v]$ 
hence $Q(v)\in\QQ[v,v\i]$. The lemma is proved.

\proclaim{Theorem 35.23}Assume that $D$ is clean (see 33.4(b)). Let $A\in\ci_n$. Let 
$u\in\cu^\fa$ be such that $E_u$ is quasi-rational (see 34.20). Then
$b_{A,u}^v\in\et\QQ$ for some $\et$, a root of $1$.
\endproclaim
Using 35.19(g) with $u_1=u'_1=u$ we see that 
$\sum_{A\in\ci_n}b_{A,u}^v(b_{A,u}^v)^\sp=1$. Using 35.22 we write $b_{A,u}^v=\et_AQ_A$
for $A\in\ci_n$ where $\et_A$ is a root of $1$ and $Q_A\in\QQ[v,v\i]$. Then 
$(b_{A,u}^v)^\sp=\et_A\i Q_A$ so that $\sum_{A\in\ci_n}Q_A^2=1$. Since 
$Q_A\in\QQ[v,v\i]$, this forces each $Q_A$ to be a constant. The theorem is proved.

\enddocument